\documentstyle[amssymb,10pt]{amsart}

\newtheorem{theorem}{Theorem}[section]
\newtheorem{lemma}[theorem]{Lemma}
\newtheorem{proposition}[theorem]{Proposition}
\newtheorem{corollary}[theorem]{Corollary} 
\theoremstyle{definition}  
\newtheorem{definition}[theorem]{Definition}

\newtheorem{example}[theorem]{Example}

\newtheorem{remark}[theorem]{Remark}

\newcommand{\id}{\operatorname{id}}

\newcommand{\Ad}{\operatorname{Ad}}

\newcommand{\nc}{\newcommand}
\nc{\Symm}{{\on{Sym}}}

\newcommand{\on}{\operatorname}   
\newcommand{\eps}{\varepsilon}
\newcommand{\Ups}{\Upsilon}
 \nc{\cE}{{\cal E}}

\nc{\D}{{\mathfrak d}}
\nc{\SL}{{\mathfrak sl}}

\nc{\gt}{{\mathfrak gt}}
\nc{\grt}{{\mathfrak grt}}
\nc{\gtm}{{\mathfrak gtm}}
\nc{\grtm}{{\mathfrak grtm}}
\nc{\gtmd}{{\mathfrak gtmd}}
\nc{\grtmd}{{\mathfrak grtmd}}

\nc{\HH}{{\mathfrak h}}
\newcommand{\g}{{\mathfrak{g}}}

\renewcommand{\u}{{\mathfrak{u}}}
\newcommand{\p}{{\mathfrak{p}}}

\newcommand{\q}{{\mathfrak{q}}}

\renewcommand{\l}{{\mathfrak{l}}}

\renewcommand{\t}{{\mathfrak{t}}}
\newcommand{\SG}{{\mathfrak{S}}}
\newcommand{\f}{{\mathfrak{f}}}
\newcommand{\ul}{\underline}

\nc{\wh}{\widehat}\nc{\wt}{\widetilde}

\newcommand{\kk}{{\bf k}}
\newcommand{\Pseudo}{{\bf Pseudo}}
\newcommand{\Psdist}{{\bf Psdist}}
\newcommand{\Assoc}{{\bf Assoc}}

\newcommand{\ben}{\begin{enumerate}}
\newcommand{\een}{\end{enumerate}}
\newcommand{\ad}{{\rm ad}}
\renewcommand{\i}{{\rm i}}

\newcommand{\CC}{{\mathbb{C}}}
\newcommand{\QQ}{{\mathbb{Q}}}
\newcommand{\RR}{{\mathbb{R}}}
\newcommand{\PP}{{\mathbb{P}}}

\newcommand{\ZZ}{{\mathbb{Z}}}
\newcommand{\cC}{{\mathcal C}}

\newcommand{\cF}{{\mathcal F}}
\newcommand{\cS}{{\mathcal S}}

\newcommand{\cM}{{\mathcal M}}

\newcommand{\NN}{{\mathbb N}}

\begin{document}

\title[Quasi-reflection algebras and cyclotomic associators]{Quasi-reflection
algebras and cyclotomic associators}

\begin{abstract} 
We develop a cyclotomic analogue of the theory of associators. 
Using a trigonometric version of the universal KZ equations, 
we prove the formality of a morphism $B_n^1 \to (\ZZ/N\ZZ)^n \rtimes \SG_n$, 
where $B_n^1$ is a braid group of type $B$. The formality isomorphism
depends algebraically on a series $\Psi_{\on{KZ}}$, the ``KZ pseudotwist''. 
We study the scheme of pseudotwists and show that it is a torsor under 
a group $\on{GTM}(N,\kk)$, mapping to Drinfeld's group $\on{GT}(\kk)$, 
and whose Lie algebra is isomorphic to its associated graded 
$\grtm(N,\kk)$. We prove that Ihara's subgroup $\on{GTK}$ of the 
Grothendieck-Teichm\"uller group, defined using distribution relations,
in fact coincides with it. We show that the subscheme of pseudotwists 
satisfying distribution relations is a subtorsor. We study 
the corresponding analogue $\grtmd(N,\kk)$ of $\grtm(N,\kk)$; 
it is a graded Lie algebra with an action of  $(\ZZ/N\ZZ)^\times$, and 
we give a lower bound for the character of its space of generators. 
\end{abstract}

\author{Benjamin Enriquez}
\address{IRMA (CNRS), rue Ren\'e Descartes, F-67084 Strasbourg, France}
\email{enriquez@@math.u-strasbg.fr}

\maketitle

\section*{Introduction}

\subsection{} 
The purpose of this paper is to develop a cyclotomic analogue of
Drinfeld's theory of associators. 
This theory can be described as follows. In \cite{Ko}, Kohno 
proves that the pure braid group of order $n$ is formal, 
i.e., its Malcev Lie algebra is isomorphic to its associated
graded Lie algebra, which has an explicit presentation; this result 
relies on the universal Knizhnik-Zamolodchikov (KZ) equations, 
but it can also be proved in the context of minimal model theory
(\cite{Su}). Kohno's isomorphisms involve the periods of 
$\PP^1 - \{0,1,\infty\}$ and are therefore defined over $\CC$. 
One of the consequences of \cite{Dr2} is that these isomorphisms are
defined over $\QQ$. To prove this,  
Drinfeld shows that the isomorphisms algebraically depend on a 
particular element in a prounipotent group $F_2(\CC)$, the KZ associator
$\Phi_{\on{KZ}}$. 
He then defines the scheme of associators over a field $\kk$ with 
$\on{char}(\kk)=0$ as the set of elements of $\kk\times F_2(\kk)$ sharing 
some properties of\footnote{We set $\i:=\sqrt{-1}$.} 
$(2\pi\i,\Phi_{\on{KZ}})$ and proves that a rational 
associator exists. The main point of this proof is the introduction of a 
group $\on{GT}(\kk)$, over which the scheme of associators is a torsor. 
Its definition relies on the categorical meaning of associators: 
associators allow to construct quasitriangular quasibialgebras (QTQBA's)
and therefore $\kk[[\hbar]]$-linear braided monoidal categories, in which 
the image of the pure braid groups are close to 1; $\on{GT}(\kk)$
may be viewed as the group of automorphisms of such structures. 
This group admits profinite and pro-$l$ variants $\wh{\on{GT}}$ and 
$\on{GT}_l$, where these categories
are replaced by braided monoidal categories, in which the image of the 
pure braid groups are finite or $l$-groups. These variants play a role in 
arithmetics; $\on{GT}(\kk)$ may be viewed as an explicit
approximation of the motivic fundamental group of $\PP^1 - \{0,1,\infty\}$. 

\subsection{} In this paper, we use universal versions of the ``cyclotomic''
KZ system introduced in \cite{EE} (see also \cite{GL} for $N=2$) 
to prove (Section \ref{sec:1}) the formality of the morphism 
$B_n^1\to (\ZZ/N\ZZ)^n \rtimes \SG_n$ ($n,N\geq 1$; $B_n^1$ is a braid group
of type $B$). We show (Section \ref{sec:2}) that the 
formality isomorphisms algebraically depend on a pair
$(\Phi_{\on{KZ}},\Psi_{\on{KZ}})$, where $\Psi_{\on{KZ}}$ is a generating 
series for periods of $\CC^\times - \mu_N(\CC)$. We then introduce a 
scheme $\ul\Pseudo(N,\kk)$ of all quadruples $(a,\lambda,\Phi,\Psi)
\in\ZZ/N\ZZ\times \kk\times F_{2}(\kk)\times F_{N+1}(\kk)$ 
sharing some properties of $(-\bar 1,2\pi\i,\Phi_{\on{KZ}},\Psi_{\on{KZ}})$; 
these are called pseudotwists. 
A pseudotwist may be viewed as a universal object allowing to 
construct quasi-reflection algebras over quasitriangular quasi-bialgebras 
(Sections \ref{sec:ps:1}, \ref{sec:QRA}), and the elements of the preimage 
$\Pseudo(N,\kk)$ of $(\ZZ/N\ZZ)^\times\times \kk^\times$
by $\ul\Pseudo(N,\kk)\to (\ZZ/N\ZZ)\times\kk$ give rise to 
formality isomorphisms for $B_n^1\to (\ZZ/N\ZZ)^n \rtimes \SG_n$. 
To show that a $\Pseudo(N,\QQ)$ is nonempty, we show 
(Section \ref{sect:7}) that $\Pseudo(N,\kk)$ is a torsor over a group 
$\on{GTM}(N,\kk)$, 
which can be defined as the group of automorphisms of braided 
module categories over $\kk[[\hbar]]$-linear braided monoidal categories, 
in which the image of $K_{n,N} := \on{Ker}(B_n^1\to (\ZZ/N\ZZ)^n \rtimes 
\SG_n)$ is close to 1 (Section \ref{sec:GTM}); this group is defined 
by imposing the defining conditions of a module version GTM of the group 
GT in products of partially completed groups. The group $\on{GTM}(N,\kk)$ 
admits a pro-$l$ variant
$\on{GTM}(N)_l$ of automorphisms of braided module 
categories over braided monoidal categories, in which the image of $K_{n,N}$ 
is a $l$-group. We have natural morphisms 
$\pi_{NN'}:\on{GTM}(N)_l \to \on{GTM}(N')_l$ for $N'|N$
(this is an isomorphism if $N' = N/l^\alpha$, 
where $\alpha$ is the $l$-adic valuation of $N$), 
and morphisms $\on{GTM}(N)_{l}\to \on{GTM}(N,\QQ_{l})$
and $\wh{\on{GT}}\to \on{GTM}(N)_{l}$ (Section \ref{sec:GTMl}). 
The map $\on{Gal}(\bar\QQ/\QQ)\to \on{GT}_{l}$ studied by 
Ihara \cite{Ih1} factors through a sequence of maps
$\on{Gal}(\bar\QQ/\QQ)\hookrightarrow
\wh{\on{GT}}\to \on{GTM}(N)_{l}\to \on{GT}_{l}$. 
One can show that the natural morphism 
$\on{Gal}(\bar\QQ/\QQ)\to \on{Out}(\pi_1^{orb}
(({\Bbb G}_m - \mu_N)/D_N)^{(l)})$ 
factors through $\on{Gal}(\bar\QQ/\QQ)\to \on{GTM}(N)_{l}$ 
(Section \ref{sect:galois}). 

We study the group $\on{GRTM}(N,\kk)$ of $\on{GTM}(N,\kk)$-automorphisms of 
$\Pseudo(N,\kk)$; it is non-canonically isomorphic to $\on{GTM}(N,\kk)$
and is equipped with a morphism $\on{GRTM}(N,\kk)\to (\ZZ/N\ZZ)^\times
\times \kk^\times$. The Lie algebra $\grtm_{(\bar 1,1)}(N,\kk)$ of the 
preimage of $(\bar 1,1)$ is graded, with a $(\ZZ/N\ZZ)^\times$-action. 

We then study distribution relations. In Section \ref{sec:ih} (which can 
be read independently),  
we prove that the inclusion $\on{GTK} \subset \wh{\on{GT}}$ (\cite{Ih2})
is in fact an equality. We then construct distribution morphisms 
$\delta_{NN'}:\on{GTM}(N)_l\to \on{GTM}(N')_l$, for $N'|N$. Using 
$\pi_{NN'}$ and $\delta_{NN'}$, one then constructs
a subgroup $\on{GTMD}(N)_l \subset \on{GTM}(N)_l$. 
A consequence of Section \ref{sec:ih} is that we have a factorization
$\wh{\on{GT}}\to \on{GTMD}(N)_l \subset \on{GTM}(N)_l$. 
We also construct a field version $\on{GTMD}(N,\kk)\subset \on{GTM}(N,\kk)$. 

We then show that the level $N$ version $\Psi^N_{\on{KZ}}$ of $\Psi_{\on{KZ}}$
satisfies distribution relations. This means that as in \cite{Go,Rac}, the pseudotwists 
$\Psi_{\on{KZ}}^N$ and $\Psi_{\on{KZ}}^{N'}$ are related by two 
transformations for $N'|N$ related to $\pi_{NN'}$ and $\delta_{NN'}$; 
this implies a system of equations satisfied 
by $\Psi^N_{\on{KZ}}$. Adjoining these conditions
to the defining equations of $\ul\Pseudo(N,\kk)$, we obtain a diagram 
$\kk^{\nu} \leftarrow \ul\Psdist(N,\kk) \hookrightarrow \ul\Pseudo(N,\kk)$
(here $\nu$ is the number of distinct prime factors of $N$). 
We show that ${\Psdist}(N,{\kk}):= {\ul\Psdist}(N,{\kk}) \cap  
\Pseudo(N,\kk)$ is a torsor under $\on{GTMD}(N,\kk) \subset 
\on{GTM}(N,\kk)$ and $\on{GRTMD}(N,\kk) \subset \on{GRTM}(N,\kk)$. 
We denote the Lie algebra of the kernel of the composed morphism 
$\on{GRTMD}(N,\kk) \subset \on{GRTM}(N,\kk)
\to (\ZZ/N\ZZ)^\times \times \kk^\times$ by $\grtmd_{(\bar 1,1)}(N,\kk)$. 
This is a $\ZZ_{\geq 0}$-graded Lie 
algebra, equipped with an action of $(\ZZ/N\ZZ)^\times$ and a morphism 
of graded Lie algebras $\grtmd_{(\bar 1,1)}(N,\kk) \to \grt_1(\kk)$.  
We use Drinfeld's method (\cite{Dr2}, Proposition 6.3) to 
construct elements of $\grtmd_{(\bar 1,1)}(N,\CC)$; using the nonvanishing of 
Dirichlet $L$-functions, and the proof of \cite{W} of the 
theorem of Bass on cyclotomic units, we derive lower bounds for the character
(as a $(\ZZ/N\ZZ)^\times$-module) of the space of generators ${\mathfrak g}/
[{\mathfrak g},{\mathfrak g}]$ of ${\mathfrak g} 
:= \grtmd_{(\bar 1,1)}(N,\CC)$ (Section \ref{sect:9}). 
We derive from this the surjectivity of a group morphism 
$\on{GRTMD}(N,\kk)\to (\ZZ/N\ZZ)^{\times}\times (\kk^{\nu}
\rtimes \kk^{\times})$ and therefore of the map $\Psdist(N,\kk) \to 
(\ZZ/N\ZZ)^{\times} \times (\kk^\nu \rtimes \kk^\times)$ (Section \ref{sect:rat}). 

Finally, we relate the Lie algebras $\grtm_{(\bar 1,1)}(N,\kk)$ and 
$\grtmd_{(\bar 1,1)}(N,\kk)$ with a 
generalization of Ihara's algebra, defined using special derivations
(Section \ref{sect:Ihara}). 

\subsection{}

This paper is related to the following works:  

(a) in \cite{Rac}, Racinet introduces a torsor $\on{DMRD}(N,\kk)$
over a prounipotent Lie group, based on the shuffle and distribution 
relations satisfied by special values of multiple polylogarithms (which coincide 
with the periods involved in $\Psi_{\on{KZ}}$, see Appendix
\ref{app:a}). 
He conjectures that for $N=1$, there is a morphism of torsors 
$\on{DMRD}(1,\kk) \to \ul\Assoc(\kk)$ taking the analogue of $\Phi_{\on{KZ}}$
in $\on{DMRD}(1,\CC)$ to $\Phi_{\on{KZ}}$, and that this is an isomorphism. 
A natural generalization of this conjecture would be to 
expect a morphism of torsors $\on{DMRD}(N,\kk) \to \ul\Psdist(N,\kk)$ when 
$N\neq 1$.  

(b) in \cite{DG}, Section 5, Deligne and Goncharov studied the motivic 
fundamental group of ${\mathbb G}_m - \mu_N$, of which 
$\on{GTMD}(N,\kk)$ should be an explicit approximation.
The element $\Psi_{\on{KZ}}$ appears in \cite{DG}, Proposition 5.17 and 
``motivic" conditions satisfied by $\Psi_{\on{KZ}}$ are given in Section 
5.19. The lower bound found in Section \ref{sect:9} for the character of 
the space of generators of $\grtmd_{(\bar 1,1)}(N,\CC)$ seems related 
to the upper bound of \cite{DG}, Cor. 5.25 for the dimension of the 
$\QQ$-linear span of the coefficients of $\Psi_{\on{KZ}}$ of fixed degree. 

(c) in \cite{Ih1}, Ihara studies the kernel of
$\on{Gal}(\bar\QQ/\QQ(\mu_{l^\infty}))\to \on{GT}_l$; 
this leads to a graded Lie algebra over $\ZZ_l$, in which 
he constructed analogues of the Drinfeld generators (using 
Soul\'e's cyclotomic characters) and conjectural integral
versions of these generators. These constructions 
should generalize to the morphism 
$\on{Gal}(\bar\QQ/\QQ(\mu_{N l^\infty}))\to \on{GTMD}(N)_l$. 

(d) using twists of QRA's over QTQBA's, one can hope for a 
generalization of the Kohno-Drinfeld theorem to representations of 
$B_n^1$ arising from specialized versions of the cyclotomic KZ
equations.  This study was started in \cite{GL} for $N=2$. 

\subsection*{Acknowledgements} This project was started in collaboration 
with P. Etingof in September 2003. I wish to thank him heartily 
for collaboration in \cite{EE} and for many discussions 
at several stages of this work. I also would like to thank D. Calaque, 
F. Lecomte, V. Leksin, I. Marin and V. Toledano Laredo for 
comments.

\section{The formality of $B_n^1\to (\ZZ/N\ZZ)^n\rtimes \SG_n$} \label{sec:1}

We fix a base field $\kk$ of characteristic $0$. 
If $\Gamma$ is a finitely generated group, we denote by $\Gamma(\kk)$
the corresponding Malcev pro-unipotent group and by $\on{Lie}(\Gamma(\kk))$
the Lie algebra of $\Gamma(\kk)$; if $\kk/\QQ$ is a finite field extension, then 
$\on{Lie}(\Gamma(\kk)) = \on{Lie}(\Gamma(\QQ)) \otimes\kk$. 

\subsection{Partial prounipotent completions} \label{partial:prounip}

Let $\varphi : \Gamma\to \Gamma_0$ be a surjective group morphism. Assume that 
$\Gamma_1:= \on{Ker}\varphi$ is finitely generated. Then one can construct a 
non-connected pro-algebraic group $\Gamma(\varphi,\kk)$, fitting in an exact 
sequence $1\to \Gamma_1(\kk)\to \Gamma(\varphi,\kk)\to\Gamma_0\to 1$, 
and a group morphism $\Gamma\to\Gamma(\varphi,\kk)$, such that the diagram 
$$
\begin{matrix}
1 \to & \Gamma_1 & \to & \Gamma & \to & \Gamma_0 & \to 1 \\  
      & \downarrow & & \downarrow & & \parallel & \\ 
1 \to & \Gamma_1(\kk) & \to & \Gamma(\varphi,\kk) & \to & \Gamma_0 & \to 1  
\end{matrix}
$$ 
commutes, and with the following properties. 

The commuting triangle $\begin{matrix} \Gamma & \to & \Gamma(\varphi,\kk) \\
\searrow & & \swarrow \\  & \Gamma_{0}&  \end{matrix}$ is an initial object 
in the category where objects are similar commuting triangles, where 
$\Gamma(\varphi,\kk)$ is replaced by a non-connected prounipotent group $U$, 
and morphisms are morphisms $U\to U'$, such that the resulting tetrahedron 
commutes. The group $\Gamma(\varphi,\kk)$ is then called the partial 
prounipotent completion of $\Gamma\stackrel{\varphi}{\to}\Gamma_{0}$. 
A commuting square $\begin{matrix} \Gamma & \stackrel{\varphi}{\to} & 
\Gamma_{0} \\ \downarrow & & \downarrow \\ \Gamma' & \stackrel{\varphi'}{\to} 
& \Gamma'_{0}\end{matrix}$ gives rise to a commuting square 
$\begin{matrix} \Gamma(\varphi,\kk) & \to & 
\Gamma_{0} \\ \downarrow & & \downarrow \\ \Gamma'(\varphi',\kk) & 
\to & \Gamma'_{0}\end{matrix}$. 

\subsection{Formality of a group morphism}

Recall that the finitely generated group $\Gamma$ is called formal iff there exists a 
Lie algebra isomorphism $\on{Lie}(\Gamma(\kk)) \to \hat{\on{gr}}\on{Lie}(\Gamma(\kk))$, 
whose associated graded morphism is the identity (here 
$\on{gr}({\mathfrak g})$ is the graded Lie algebra associated to a 
pronilpotent Lie algebra ${\mathfrak g}$ and $\hat{\on{gr}}({\mathfrak g})$ 
is the degree completion of $\on{gr}({\mathfrak g})$). 

In the above situation, the choice of a section $\gamma_0\mapsto 
\tilde\gamma_0$ of $\Gamma(\kk)\to \Gamma_0$ 
induces a map $\Gamma_0\to \on{Aut}(\on{Lie}(\Gamma_1(\kk)))$, 
$\gamma_0\mapsto \theta_{\gamma_0}$, where $\theta_{\gamma_0}(x)
= \tilde\gamma_0 x \tilde\gamma_0^{-1}$; the map $\gamma_0\mapsto 
\theta_{\gamma_0}$ may fail to be a group morphism (by inner automorphisms).
Each $\theta_{\gamma_0}$ induces a graded Lie algebra automorphism 
$\on{gr}\theta_{\gamma_0}\in \on{Aut}(\on{gr}\on{Lie}(\Gamma_1(\kk)))$, 
which is independent of the choice of the lifting 
$\gamma_0\mapsto \tilde\gamma_0$, and $\gamma_0\mapsto 
\on{gr}\theta_{\gamma_0}$ is a group morphism. We can then form the semidirect 
product $\on{exp}(\hat{\on{gr}}\on{Lie}\Gamma_1(\kk))\rtimes\Gamma_0$. 

We say that the group morphism $\varphi : \Gamma\to \Gamma_0$ is formal 
if there exists a group isomorphism $\Gamma(\kk,\varphi) \simeq 
\on{exp}(\on{gr}\on{Lie}\Gamma_1(\kk))\rtimes\Gamma_0$, restricting to a 
formality isomorphism for $\Gamma_1$, and such that the diagram 
$$
\begin{matrix}
1 \to & \Gamma_1(\kk) & \to & \Gamma(\kk,\varphi) & \to & \Gamma_0 & \to 1 \\  
      & \wr\mid & & \wr\mid & & \parallel & \\ 
1 \to & \on{exp}(\hat{\on{gr}}\on{Lie}\Gamma_1(\kk)) & \to & 
\on{exp}(\hat{\on{gr}}\on{Lie}\Gamma_1(\kk))\rtimes\Gamma_0 & \to & \Gamma_0 & \to 1  
\end{matrix}
$$ 
commutes. 

For example, the morphism 
$B_n \to \SG_n$ is formal, where $B_n$ is the braid group of order $n$ 
(when $\kk = \CC$, this follows from \cite{Ko}, and for $\kk = \QQ$, 
from \cite{Dr2}). 

\subsection{The morphism $\varphi_{n,N} : B_n^1\to (\ZZ/N\ZZ)^n\rtimes S_n$} 
\label{sec:phi:nN}

Let $B_n^1:= B_{n+1}\times_{\SG_{n+1}}\SG_n$, where $\SG_n\subset \SG_{n+1}$
is the subgroup of permutations of $[0,N]$ such that $\sigma(0)=0$. We have  
$B_n^1 = \pi_1(X_n,\SG_n b)$, where $X_n := [(\CC^\times)^n
-\cup_{1\leq i<j\leq n}D_{ij}]/\SG_n$, $D_{ij} = \{(Z_1,...,Z_n)|Z_i=Z_j\}$
and $b:= \{(Z_{1},...,Z_{n})
|0<Z_{1}<...<Z_{n}\}$. 

Let $W_{n,N} := (\CC^\times)^n - \cup_{1\leq i<j\leq n}\cup_{\zeta\in 
\mu_N(\CC)} D_{ij,\zeta}$, where $D_{ij,\zeta} = \{(z_1,...,z_n)|z_i=
\zeta z_j\}$. Then the map $W_{n,N}\to X_n$ induced by $(z_1,...,z_n)
\to [(z_1^N,...,z_n^N)]$ is a covering with group $(\ZZ/N\ZZ)^n\rtimes \SG_n$; 
it induces a group morphism 
$\varphi_{n,N} : B_n^1\to (\ZZ/N\ZZ)^n\rtimes \SG_n$. 
We have then $K_{n,N} := \on{Ker}\varphi_{n,N} \simeq 
\pi_1(W_{n,N},b)$. 

We now introduce the universal ``cyclotomic'' KZ system and use it to 
prove the formality of $\varphi_{n,N}$. 

\subsection{The universal cyclotomic KZ system}

Let $\t_{n,N}(\kk)$ (or $\t_{n,N}$) be the $\kk$-Lie algebra with generators
$t_0^{0i}$, ($i\in [1,n]$), and $t(a)^{ij}$, 
($i\neq j\in [1,n]$, $a\in \ZZ/N\ZZ$), and relations: 
$$
t(a)^{ij} = t(-a)^{ji}, \quad 
[t(a)^{ij},t(a+b)^{ik} + t(b)^{jk}] =0, 
\quad 
[t_0^{0i} + t_0^{0j} + \sum_{b\in \ZZ/N\ZZ}t(b)^{ij},
t(a)^{ij}]=0, 
$$
$$
[t_0^{0i} ,t_0^{0j} + \sum_{b\in\ZZ/N\ZZ}t(b)^{ij}]=0,
\quad 
[t_0^{0i},t(a)^{jk}]=0, 
\quad 
[t(a)^{ij},t(b)^{kl}]=0
$$
($i,j,k,l\in[1,n]$ are all distinct and $a,b\in\ZZ/N\ZZ$).
The group $(\ZZ/N\ZZ)^n \rtimes \SG_n$ acts on $\t_{n,N}$
as follows: let $s_i=(\bar 0,...,\bar 1,...,\bar 0)$ be the $i$th 
generator of $(\ZZ/N\ZZ)^n$, then $s_i \cdot t_{0}^{0j} = t_0^{0j}$
and $s_i \cdot t(a)^{jk} = t(a+\delta_{ij} - \delta_{ik})^{jk}$;
for $\sigma\in\SG_n$, $\sigma \cdot t_0^{0i} = t_0^{0\sigma(i)}$, 
$\sigma\cdot t(a)^{ij} = t(a)^{\sigma(i) \sigma(j)}$. 

When $\kk = \CC$, we define $t[\zeta]^{ij}$ for $\zeta\in\mu_N(\CC)$
by $t[\zeta_N^a]^{ij} = t(a)^{ij}$, where 
$\zeta_N := e^{2\pi\i/N}\in \mu_N(\CC)$. 

Set 
$$
K_i(z_1,...,z_n):= {{t_0^{0i}}\over {z_i}} + \sum_{j\in [1,n] - \{i\}}
\sum_{\zeta\in\mu_N(\CC)} {{t[\zeta]^{ij}}\over{z_i - \zeta z_j}}, 
$$
then the connection 
$$
\on{d} - \sum_{i=1}^n K_i(z_1,...,z_n)\on{d}z_i 
$$
on the trivial 
$\on{exp}(\hat\t_{n,N}(\CC))$-bundle\footnote{$\hat{\mathfrak g}$ 
denotes the degree completion of a 
$\NN$-graded Lie algebra; here $\on{deg}(t_0^{0i})=\on{deg}(t(a)^{ij})=1$.} 
over $W_{n,N}$ is flat.  
This connection is also equivariant w.r.t. the action of 
$\Gamma_0 := (\ZZ/N\ZZ)^n \rtimes \SG_n$, therefore its descends to a 
flat connection on the $\on{exp}(\hat\t_{n,N}(\CC))\rtimes\Gamma_0$-bundle 
over $X_{N}$, 
such that the set of sections of ${\cal U} \subset X_{N}$
is the set of functions $f:\pi^{-1}({\cal U})\to \on{exp}(\hat\t_{n,N}(\CC))\rtimes 
\Gamma_0$, such that $f(\gamma_0 \cdot x) = f(x) \gamma_0^{-1}$ for 
$\gamma_0 \in \Gamma_0$
(here $\pi : W_{n,N} \to X_n$ is the canonical projection). 

\subsection{Monodromy morphisms}

There exists a unique solution $F : \tilde X_{N} \to \on{exp}(\hat\t_{n,N}(\CC))$
of the system $\partial F/\partial z_i = K_i F$, with asymptotic behavior
$F(z_1,...,z_n) \simeq z_1^{t_0^{0i}}z_2^{t_0^{01,2}}...z_n^{t_0^{0...n-1,n}}$
when $z_1 \ll z_2 ...\ll z_n$; 
here $t_0^{0...r-1,r}=t_0^{0r}+\sum_{i=1}^{r-1}\sum_{\zeta\in\mu_N(\CC)}
t[\zeta]^{ir}$. Here $\tilde X_{N}\to X_{N}$ 
is the universal cover map.  

If now $\gamma\in B_n^1$, the map $z\mapsto F(z)^{-1}\varphi_{n,N}(\gamma)
F(\gamma^{-1}z)\in \on{exp}(\hat\t_{n,N}(\CC))\rtimes\Gamma_0$ is a constant, 
called $\mu_{n}(\gamma)$. 
This defines a group morphism 
$$
\mu_n : B_n^1 \to \on{exp}(\hat\t_{n,N}(\CC))\rtimes \Gamma_0,  
$$
which factors through a morphism $B_n^1(\varphi_{n,N},\CC) \to 
\on{exp}(\hat\t_{n,N}(\CC))\rtimes \Gamma_0$. This morphism commutes with the 
projections to $\Gamma_0$, and therefore restricts to a morphism 
$\mu_n(\CC) : K_{n,N}(\CC)\to\on{exp}(\hat\t_{n,N}(\CC))$. 

\subsection{Formality of $\varphi_{n,N}$}

To prove that $\varphi_{n,N}$ is formal, it then suffices to prove 
that $\on{Lie}\mu_n(\CC) : \on{Lie}(K_{n,N}(\CC)) \to \hat\t_{n,N}(\CC)$
is an isomorphism. 

$K_{n,N}$ may be be viewed as the kernel of the morphism $P_{n+1}\to 
(\ZZ/N\ZZ)^n$ given by $x_{ij}\mapsto 0$ if $i,j\neq 0$,  
$x_{0i}\mapsto s_i = (\bar 0...\bar 1...\bar 0)$ if $i\neq 0$
(here the $x_{ij}$, $0\leq i<j\leq n$ are the Artin generators of the 
pure braid group $P_{n+1}$). 
Set $X_{0i}:= x_{0i}^N$, $x_{ij}(\alpha) := x_{0i}^\alpha x_{ij} x_{0i}^{-\alpha}$ 
for $\alpha\in\ZZ$. 

One can show that the $X_{0i}$ and $x_{ij}(\alpha)$ ($i\neq j\in [1,n]$, 
$\alpha\in [0,N-1]$) are generators of $K_{n,N}$. We have: 
(a) $\on{Lie}\mu_n(\CC)$ takes $[\on{log}X_{0i}]$ and 
$[\on{log}x_{ij}(\alpha)]$ to multiples to $t_0^{0i}$ and $t(\bar\alpha)^{ij}$
(here $X_{0i},x_{ij}(\alpha)$ are identified with their images in 
$K_{n,N}(\CC)$, $[x]$ is the class of $x\in \on{Lie}(K_{n,N}(\CC))$
in $\on{gr}_1\on{Lie}(K_{n,N}(\CC))$) and $\alpha\mapsto\bar\alpha$
is the canonical projection $\ZZ/\to\ZZ/N\ZZ$; (b) one can exhibit relations 
between these generators, whose logarithms have degree $\geq 2$, 
and whose degree 2 parts coincide with the defining relations of 
$\t_{n,N}(\CC)$. These relations are explicitly given as follows
(recall that $x_{ij}(\alpha+N)=X_{0i}x_{ij}(\alpha)X_{0i}^{-1}$, so 
$[\on{log}x_{ij}(\alpha)]$ depends only on $\bar\alpha$): 

\begin{proposition} \label{beginning}
$$
(X_{0j}x_{ij}(0)...x_{ij}(N-1)X_{0i},X_{0j}) = 
(X_{0j}x_{ij}(0)...x_{ij}(N-1)X_{0i},x_{ij}(\alpha)) = 1,  
$$
\begin{align} \label{new:rel:1}
& ( x_{ij}(\alpha)x_{ik}(\alpha+\beta|\alpha)u_{ik}(\alpha,\beta)
x_{jk}(\beta|0)u_{ik}(\alpha,\beta)^{-1},x_{ij}(\alpha)) 
\\ & \nonumber = ( x_{ij}(\alpha)x_{ik}(\alpha+\beta|\alpha)u_{ik}(\alpha,\beta)
x_{jk}(\beta|0)u_{ik}(\alpha,\beta)^{-1},
x_{ik}(\alpha+\beta|\alpha)) = 1,  
\end{align}
\begin{align} \label{new:rel:2}
& (x_{ij}(\alpha),u_{ik}(\alpha,-\beta)x_{kl}(\beta)u_{ik}(\alpha,-\beta)^{-1}) =
(u_{ij}(\alpha+1,-\beta)^{-1}x_{il}(\alpha) u_{ij}(\alpha+1,-\beta),x_{jk}(\beta)) 
\\ & \nonumber = (\on{Ad}[u_{ij}(\alpha+1,-\beta-1)^{-1}][x_{ik}(\alpha)], 
\on{Ad}[u_{ij}(\alpha,-\beta-1)^{-1} u_{ij}(\alpha,-\beta)][x_{jl}(\beta)])= 1, 
\end{align}
\begin{equation} \label{new:rel:3}
(u_{ij}(N,-\alpha)^{-1}X_{0i},x_{jk}(\alpha)) = (u_{ik}(\alpha,-N)X_{0k},x_{ij}(\alpha)) = 
(u_{ij}(\alpha+1,-N)X_{0j},x_{ik}(\alpha)) = 1
\end{equation}
where $i<j<k<l$, $\alpha,\beta\in\ZZ$,  
$x_{ik}(\alpha+\beta|\alpha) := \on{Ad}[x_{ik}(\alpha)x_{ik}(\alpha+1)...
x_{ik}(\alpha+\beta-1)](x_{ik}(\alpha+\beta))$
and $u_{ik}(\alpha,\beta)=x_{ik}(\alpha)...x_{ik}(\alpha+\beta-1)
x_{ik}(\beta-1)^{-1}...x_{ik}(0)^{-1}$; 
when $\beta\leq 0$, $x_{ik}(\alpha)x_{ik}(\alpha+1)...x_{ik}(\alpha+\beta-1)$ means
$x_{ik}(\alpha-1)^{-1}x_{ik}(\alpha-2)^{-1}...x_{ik}(\alpha+\beta)^{-1}$ (this is $1$ if 
$\beta=0$). 
\end{proposition}

{\em Proof.} These relations are consequences of 
the pure braid group relations $(x_{ij}x_{ik}x_{jk},x_{ij})
= (x_{ij}x_{ik}x_{jk},x_{jk})=1$ for $i<j<k$, $(x_{ij},x_{kl})=(x_{il},x_{jk})=
(x_{ik},x_{ij}^{-1}x_{jl}x_{ij})=1$ for $i<j<k<l$ (\cite{Art}); for example, 
(\ref{new:rel:1}) is the conjugation of $(x_{ij}x_{ik}x_{jk},x_{ij})
=(x_{ij}x_{ik}x_{jk},x_{ik})=1$ by $x_{0i}^{\alpha}x_{0k}^{-\beta}$, using the 
relations $\on{Ad}(x_{0i}^\alpha x_{0k}^{-\beta})(x_{ik}) 
= x_{ik}(\alpha+\beta|\beta)$ and 
$(x_{0i}^\alpha,x_{0k}^{-\beta}) = u_{ik}(\alpha,b)$; (\ref{new:rel:2}) is the conjugation 
of $(x_{ij},x_{jk})=1$ by $x_{0i}^\alpha x_{0j}^\beta$, of 
$(x_{il},x_{jk})=1$ by $x_{0j}^\beta x_{0i}^\alpha$ and of $(x_{ik},x_{jl}(-1))=1$ 
by $x_{0j}^{\beta+1}x_{0i}^\alpha$; for example, in this last case, we have 
$\on{Ad}[x_{0j}^{\beta+1}x_{0i}^\alpha][x_{ik}] = \on{Ad}[(x_{0j}^{\beta+1},
x_{0i}^{\alpha+1})
x_{0i}^{\alpha+1}x_{0j}^{\beta+1}][x_{ik}(-1)] = \on{Ad}[u_{ij}(\alpha+1,-\beta-1)^{-1}]
[x_{ik}(\alpha)]$ and $\on{Ad}[x_{0j}^{\beta+1}x_{0i}^\alpha][x_{jl}(-1)] = 
\on{Ad}[(x_{0j}^{\beta+1},x_{0i}^\alpha)(x_{0i}^\alpha,x_{0j}^\beta)x_{0j}^\beta 
x_{0i}^\alpha][x_{jl}] 
= \on{Ad}[u_{ij}(\alpha,-\beta-1)^{-1} u_{ij}(\alpha,-\beta)][x_{jl}(\beta)]$; 
(\ref{new:rel:3})
is the conjugation of $(x_{0i}^N,x_{jk})=1$ by $x_{0j}^\alpha$, of 
$(x_{0k}^N,x_{ij})=1$ by $x_{0i}^\alpha$, and of $(x_{0j}^N,x_{ik}(-1))=1$
by $x_{0i}^{\alpha+1}$. \hfill \qed \medskip 

It then follows from (b) that we have a morphism $\alpha_n : \t_{n,N}(\CC)\to 
\on{gr}\on{Lie}(K_{n,N}(\CC))$ given by $t_0^{0i}\mapsto [\on{log}X_{0i}]$, 
$t(a)^{ij}\mapsto [\on{log}x_{ij}(\tilde a)]$ ($\tilde a$ is the lift of $a$ to 
$[0,N-1]$). This morphism is surjective, 
since $\on{gr}\on{Lie}(K_{n,N}(\CC))$ is generated in degree 1, as the 
associated graded Lie algebra of the quotient of a free Lie algebra (recall that 
$\on{Lie}K_{n,N}(\CC)$ is the quotient of the topologically free Lie algebra generated
by the $\on{log}X_{0i}$, $\on{log}x_{ij}(\alpha)$ by the logarithm of their relations). 
It follows from (a) that 
$\on{gr}\mu_n(\CC) : \on{gr}\on{Lie}(K_{n,N}(\CC)) \to \t_{n,N}(\CC)$
is surjective, and one checks that $\on{gr}\mu_n(\CC) \circ \alpha_n$ 
is bijective. So $\alpha_n$ is injective, hence an isomorphism, hence so 
it $\on{gr}\mu_n(\CC)$. So $\mu_n(\CC)$ is an isomorphism, as wanted. 

\section{Algebraic construction of formality isomorphisms} \label{sec:2}

\subsection{Insertion-coproduct morphisms and structure of the $\t_{n,N}$}

\subsubsection{} Let $\t_n(\kk)$ (or $\t_{n}$) be the $\kk$-Lie 
algebra of infinitesimal pure braids 
introduced in \cite{Dr2}; it has generators $t^{ij}$, $i\neq j\in [1,n]$,
and relations $t^{ij}=t^{ji}$, $[t^{ij},t^{ik}+t^{jk}]=0$ and 
$[t^{ij},t^{kl}]=0$ for $i,j,k,l$ distinct. 

Let $f : [0,m]\to [0,n]$ be a partially defined function such that $f(0)=0$, 
then there is a unique Lie algebra morphism $\t_{n,N} \to \t_{m,N}$, 
$x\mapsto x^f$, such that 
$$
(t(a)^{ij})^f = \sum_{i'\in f^{-1}(i), j'\in f^{-1}(j)}
t(a)^{i'j'}
\quad
(i\neq j\in [1,n]) 
$$ 
$$
(t_0^{0j})^f = \sum_{j'\in f^{-1}(j)} t_0^{0j'} + 
\sum_{j',j''\in f^{-1}(j), j'<j''} \sum_{b\in\ZZ/N\ZZ}t(b)^{j'j''}
+ \sum_{i'\in f^{-1}(0) \setminus \{0\}, j'\in f^{-1}(j)} 
\sum_{b\in \ZZ/N\ZZ}t(b)^{i'j'}
\quad 
(j\neq 0). 
$$
Let $f : [1,m]\to [1,n]$ be a partially defined function. 
There is a unique Lie algebra morphism $\t_n \to \t_{m,N}$, 
$x\mapsto x^f$, such that 
$$
(t^{ij})^f = \sum_{i'\in f^{-1}(i), j'\in f^{-1}(j)} 
t(0)^{i'j'}
\quad (i\neq j\in [1,n]).  
$$
We denote $x^f$ as $x^{f^{-1}(0),\ldots,f^{-1}(n)}$ 
($x^{f^{-1}(1),\ldots,f^{-1}(n)}$ in the second case).  

We also denote by $x\mapsto x^f$ the induced morphisms 
$U(\t_{n,N})\to U(\t_{m,N})$, $U(\t_n)\to U(\t_{m,N})$. 

We denote by $T_{n,N}$ the semidirect product 
$U(\t_{n,N}) \rtimes (\ZZ/N\ZZ)^n$ and by $s_i\in T_{n,N}$ 
the image of the $i$th generator of $(\ZZ/N\ZZ)^n$. 

If $f:[0,m]\to[0,n]$ is a partially defined function with $f(0)=0$, 
then the algebra morphism $U(\t_{n,N})\to U(\t_{m,N})$, $x\mapsto x^f$
extends to an algebra morphism $T_{n,N}\to T_{m,N}$ by $s_i \mapsto 
\prod_{i'\in f^{-1}(i)} s_{i'}$.

\subsubsection{}

Let us describe the structure of $\t_{n,N}$. Let $\f_{n,N} \subset \t_{n,N}$
be its Lie subalgebra generated by $t_0^{0n}$ and the $t(a)^{in}$ 
($i\in [1,n-1]$, $a\in \ZZ/N\ZZ$). These generators freely generate
$\f_{n,N}$, which is an ideal of $\t_{n,N}$. The map $x\mapsto x^{0,\dots,n-1}$
is an injection $\t_{n-1,N} \hookrightarrow \t_{n,N}$. This map induces an 
action of $\t_{n-1,N}$ by derivations on $\f_{n,N}$, and $\t_{n,N}$ is the
corresponding semidirect product. 

When $n=2$, $\t_{2,N}$ can be described more simply as the direct sum of its
center, generated by $t_0^{01} + t_0^{02} 
+ \sum_{a\in \ZZ/N\ZZ} t(a)^{12}$, with 
the free Lie algebra $\t_{2,N}^0$ generated by $t_0^{01}$ 
and the $t(a)^{12}$ ($a\in \ZZ/N\ZZ$). We will use the following notation 
for generators of $\t_{2,N}^0$: $A:= t_0^{01}$, $b(a) := 
t(a)^{12}$, $C:=$ the image $t_0^{02}$ by the projection  $\t_{2,N}\to 
\t_{2,N}^0$, so $A+C + \sum_{a\in\ZZ/N\ZZ}b(a) = 0$. 

In general, the element $\sum_{0\leq i<j \leq n} t^{ij}_0$ is central in 
$\t_{n,N}$, where $t_0^{ij} := \sum_{a\in \ZZ/N\ZZ} t(a)^{ij}$ for 
$j>i\geq 1$. If $\t_{n,N}^0$ is the Lie algebra with the same generators 
(except $t_0^{0n}$) and relations as $\t_{n,N}$ (except those involving 
$t_0^{0n}$), then $\t_{n,N}^0 \hookrightarrow \t_{n,N}$ and 
$\t_{n,N} = \t_{n,N}^0 \oplus \kk (\sum_{0\leq i<j \leq n} t^{ij}_0)$.  

\subsection{Definition and properties of $\Psi_{\on{KZ}}$} \label{sect:PsiKZ}

We define $\Psi_{\on{KZ}}$ as the renormalized holonomy from $0$ to $1$ 
of the differential equation 
\begin{equation} \label{eq:H}
 { {\on{d}} \over {\on{d}z} } H(z) 
= \big( {{t_0^{01}}\over z}
+ \sum_{\zeta\in \mu_N(\CC)} {{t[\zeta]^{12}}\over{z-\zeta}}\big) H(z), 
\end{equation}
i.e., $\Psi_{\on{KZ}} = H_{1^{-}}^{-1}H_{0^{+}}$, where 
$H_{0^{+}},H_{1^{-}}$ are the solutions such that $H_{0^{+}}(z)
\sim z^{t_{0}^{01}}$ when $z\to 0^{+}$ and $H_{1^{-}}(z)
\sim z^{t[1]^{12}}$ when $z\to 1^{-}$. 
Then $\Psi_{\on{KZ}}$ belongs to $\on{exp}(\hat\t_{2,N}^{0}(\CC))$. 

Recall that $\Phi_{\on{KZ}}\in\on{exp}(\hat\t_{3}(\CC))$ is  the renormalized holonomy
from $0$ to $1$ of $G'(z) = ({{t^{12}}\over z} + {{t^{12}}\over{z-1}})G(z)$, 
i.e., $\Phi_{\on{KZ}} = G_{0^{+}}^{-1}G_{1^{-}}$, where 
$G_{0^{+}},G_{1^{-}}$ are the solutions such that $G_{0^{+}}(z)\sim z^{t^{12}}$ 
when $z\to 0^{+}$ and $G_{1^{-}}(z)\sim (1-z)^{t^{23}}$ when $z\to 1^{-}$. 
Recall that 
$$
\Phi_{\on{KZ}}(V,U)=\Phi_{\on{KZ}}(U,V)^{-1}, \; 
e^{\pi\i U}\Phi_{\on{KZ}}(W,U)
e^{\pi\i W}\Phi_{\on{KZ}}(V,W)
e^{\pi\i V}\Phi_{\on{KZ}}(U,V)=1, 
$$
where $U=\overline{t^{12}}\in \t_3^0 := \t_3/(t^{12}+t^{13}+t^{23})$, 
$V=\overline{t^{23}}\in \t_3^0$ and $U+V+W=0$, and  
$$
\Phi^{12,3,4}_{\on{KZ}}\Phi^{1,2,34}_{\on{KZ}}=\Phi^{2,3,4}_{\on{KZ}}
\Phi^{1,23,4}_{\on{KZ}}\Phi^{1,2,3}_{\on{KZ}} 
$$
(relation in $\on{exp}(\hat\t_4(\CC))$). 

\begin{proposition}
The pair $(\Phi_{\on{KZ}},\Psi_{\on{KZ}})$ satisfies the mixed pentagon relation
$$
\Psi^{12,3,4}_{\on{KZ}}\Psi^{1,2,34}_{\on{KZ}}=\Phi^{2,3,4}_{\on{KZ}}
\Psi^{1,23,4}_{\on{KZ}}\Psi^{1,2,3}_{\on{KZ}}. 
$$ 
in $\on{exp}(\hat\t_{3,N}^0(\CC))$
and the octogon relation
\begin{align*}
& \Psi_{\on{KZ}}(A|b(0),b(1),...)^{-1} e^{\pi\i b(0)} 
\Psi_{\on{KZ}}(C|b(0),b(-1),...)
e^{(2\pi\i/N)C} 
\\ & \Psi_{\on{KZ}}(C|b(1),b(0),b(-1),...)^{-1}e^{\pi\i b(1)}
\Psi_{\on{KZ}}(A|b(1),b(2),...)e^{(2\pi\i/N)A}=1
\end{align*}
in $\on{exp}(\hat\t_{2,N}^0(\CC))$. 
\end{proposition}

{\em Proof.} The first equation follows as in \cite{EE} from the 
the flatness of the system 
$$
{{\partial}\over{\partial z}} H = \big( {{t_0^{12}}\over z}
+ \sum_{\zeta\in \mu_N(\CC)} {{t[\zeta]^{23}}\over{z-\zeta w}}
+ \sum_{\zeta\in \mu_N(\CC)} {{t[\zeta]^{24}}\over{z-\zeta}}
\big) H , 
$$
$$
{{\partial}\over{\partial w}} H = \big( {{t_0^{13}}\over w}
+ \sum_{\zeta\in \mu_N(\CC)} {{t[\zeta]^{32}}\over{w-\zeta z}}
+ \sum_{\zeta\in \mu_N(\CC)} {{t[\zeta]^{34}}\over{w-\zeta}}
\big) H .  
$$
Let us prove the second equation. Let us denote by ${\bf D}$ the domain of 
$\CC$ of all complex numbers $z$, such that $0\leq \on{arg}(z) \leq 2\pi/N$, 
and $z\neq 0,1,e^{2\pi\i/N}$.
We set $\zeta_N = e^{2\pi\i/N}$. 

Recall that (\ref{eq:H}) may be rewritten $H'(z) = ({A\over z} + 
\sum_{a\in\ZZ/N\ZZ} {{b(a)}\over{z-\zeta_N^a}})H(z)$. 
There are unique solutions $H_{0^+}$, $H_{1^-}$, $H_{1^+}$, $H_\infty$, 
$H_{\zeta_N \infty}$, $H_{\zeta_N^+}$, $H_{\zeta_N^-}$, $H_{\zeta_N 0^+}$
of (\ref{eq:H}) in ${\bf D}$, with the following asymptotic behaviors: 
$H_{0^+}(z) \sim z^{A}$ when $z\to 0^+$, 
$H_{1^-}(z) \sim (1-z)^{b(0)}$ when $z\to 1^-$, 
$H_{1^+}(z) \sim (z-1)^{b(0)}$ when $z\to 1^+$, 
$H_\infty(z) \sim z^{-C}$ when $z\to +\infty$ on the real axis, 
$H_{\zeta_N\infty}(z) \sim (z/\zeta_N)^{-C}$ when $z/\zeta\to +\infty$
on the real axis, 
$H_{\zeta_N^+}(z) \sim ({z\over{\zeta_N}} -1)^{b(1)}$
when $z\to \zeta_N$ in such a way that $z/\zeta_N > 1$, 
$H_{\zeta_N^-}(z) \sim (1-{z\over{\zeta_N}})^{b(1)}$
when $z\to \zeta_N$ in such a way that $z/\zeta_N < 1$, 
$H_{\zeta_N 0^+}(z) \sim (z/\zeta_N)^{A}$ as $z\to 0$ 
in such a way that $z/\zeta_N >0$. 

We have: 
$$
\forall z\in ]0,1[, \; H_{0^+}(z) = H_0(z), \quad 
H_{1^-}(z) = H_1(z),  
$$
$$
\forall z\in ]1,+\infty[, \; 
H_{1^+}(z) = H_1(1/z)^{1,3,2}, 
\quad   
H_{\infty}(z) = H_0(1/z)^{1,3,2}, 
$$
$$
\forall z\in ]\zeta_N,\zeta_N\infty[, \; 
H_{\zeta_N\infty}(z) = \on{Ad}(s_1)(H_0(\zeta_N/z)^{1,3,2}), 
\quad   
H_{\zeta_N^+}(z) = 
\on{Ad}(s_1) ( H_1(\zeta_N/z)^{1,3,2}), 
$$
$$
\forall z\in ]0,\zeta_N[, \; 
H_{\zeta_N^-}(z) = \on{Ad}(s_1)(H_1(z/\zeta_N)), \quad  
H_{\zeta_N 0^+}(z) = \on{Ad}(s_1)(H_0(z/\zeta_N)).  
$$
where $x\mapsto x^{1,3,2}$ and $\on{Ad}(s_1)$ are 
the automorphisms of $\hat\t_{2,N}^0$ given by 
$(A,C,b(a)) \mapsto (C,A,b(-a))$ and $(A,C,b(0),b(1),...)
\mapsto (A,C,b(1),b(2),...)$. 

It follows that $H_{1^-} = H_{0^+} \Psi_{\on{KZ}}^{-1}$, 
$H_\infty = H_{1^+} \Psi_{\on{KZ}}^{1,3,2}$, $H_{\zeta_N^+} 
= H_{\zeta_N \infty} \on{Ad}(s_1)((\Psi_{\on{KZ}}^{-1})^{1,3,2})$, 
$H_{\zeta_N 0^+} = H_{\zeta_N^-} \on{Ad}(s_1)(\Psi_{\on{KZ}})$. 

Moreover, local study at the points $0,1,\infty,\zeta_N$ implies that
$H_{1^+}=H_{1^-}e^{\pi\i b(0)}$, $H_{\zeta_N\infty}=H_\infty e^{(2\pi\i/N)C}$, 
$H_{\zeta_N^-}=H_{\zeta_N^+}e^{\pi\i b(1)}$, $H_{0^+}=H_{\zeta_N 0^+}
e^{(2\pi\i/N)A}$.  

Therefore 
$$
H_{0^+} = H_{0^+} \cdot 
\Psi_{\on{KZ}}^{-1} e^{\pi\i b(0)} \Psi_{\on{KZ}}^{1,3,2} e^{(2\pi\i/N)C}
\on{Ad}(s_1)((\Psi_{\on{KZ}}^{-1})^{1,3,2}) e^{\pi\i b(1)}
\on{Ad}(s_1)(\Psi_{\on{KZ}}) e^{(2\pi\i/N)A}=1.  
$$
Since $H_{0^+}$ is invertible, this implies that the right factor of 
$H_{0^+}$ in this expression in $1$, i.e., the second relation. 
\hfill \qed \medskip 

Moreover, the formality isomorphism of Section \ref{sec:1} may be expressed 
algebraically in terms of $(\Phi_{\on{KZ}},\Psi_{\on{KZ}})$ as follows. 
The group $B_n^1$ has the following presentation: generators 
are $\tau,\sigma_1,...,\sigma_{n-1}$ and relations are 
\begin{equation} \label{std:braid}
\sigma_i\sigma_j = \sigma_j \sigma_i \ (|i-j|\geq 2), \; \; \; 
\sigma_i \sigma_{i+1} \sigma_i = \sigma_{i+1} \sigma_i \sigma_{i+1}
\ (i\in [1,n-1]), 
\end{equation}
$$
\tau\sigma_i = \sigma_i\tau\; (i\in [2,n]), \quad (\tau\sigma_1)^2 =
(\sigma_1\tau)^2.  
$$
The inclusion $B_n^1\subset B_{n+1}$ is then 
induced by $\sigma_i\mapsto \sigma_i$, $\tau\mapsto \sigma_0^2$
(the generators of $B_{n+1}$ are called $\sigma_{0},...,\sigma_{n-1}$). 

The monodromy of the solution $F(z_1,...,z_n)$ at order $n$ can 
be related to to that of the solution at order $2$. It follows that the 
monodromy morphism $\mu_n : B_n^1 \to \on{exp}(\hat\t_{n,N}(\CC))
\rtimes (\ZZ/N\ZZ)^n \rtimes \SG_n$ is given by 
$$
\tau\mapsto e^{(2\pi\i/N)t_0^{01}}s_1^{-1}, \quad 
\sigma_i \mapsto (\Psi_{\on{KZ}}^{-1})^{0...i-1,i,i+1}
(i,i+1) e^{\pi\i t(0)^{i,i+1}} \Psi_{\on{KZ}}^{0...i-1,i,i+1}
$$
for $i=1,...,n-1$. 

\subsection{The schemes $\Pseudo(N,\kk)\subset \ul{\Pseudo}(N,\kk)$ 
and formality isomorphisms}

Let $\kk$ be a ring over $\QQ$. Recall that $\ul{\Assoc}(\kk)$ is defined as 
the set of pairs $(\lambda,\Phi)$, where $\lambda\in \kk$ and 
$\Phi\in \on{exp}(\hat\t_{3}^{0}(\kk))$ satisfies 
\begin{equation} \label{cond:Phi:1}
\Phi(V,U) = \Phi(U,V)^{-1}, \quad e^{\lambda U/2} \Phi(W,U) 
e^{\lambda W/2} \Phi(V,W) e^{\lambda V/2} \Phi(U,V) = 1
\end{equation}
where $U,V$ generate a free Lie algebra and $U+V+W=0$, and 
\begin{equation} \label{cond:Phi:2}
\Phi^{1,2,34} \Phi^{12,3,4} = \Phi^{2,3,4} \Phi^{1,23,4} \Phi^{1,2,3}. 
\end{equation}

We then define the set of cyclotomic associators $\ul{\Pseudo}(N,\kk)$ 
as the set of quadruples 
$(a,\lambda,\Phi,\Psi)\in (\ZZ/N\ZZ \times \kk) \times 
\on{exp}(\hat\t_3^0(\kk)) \times \on{exp}(\hat\t_{2,N}^0(\kk))$, 
such that $(\lambda,\Phi)\in \ul{\Assoc}(\kk)$ and $\Psi$ satisfies 
\begin{align} \label{cond:Psi:1}
& \Psi(A | b(a),b(a+1),\ldots,b(a+N-1))^{-1} e^{(\lambda/2) b(a)} 
\Psi(C | b(a),b(a-1),\ldots,b(a+1-N)) e^{(\lambda/N) C}
\\ & \nonumber 
\Psi(C | b(0),b(-1),\ldots,b(1-N))^{-1} 
e^{(\lambda/2) b(0)}
\Psi(A | b(0),b(1),\ldots,b(N-1)) e^{(\lambda/N)A}=1, 
\end{align}  
\begin{equation} \label{cond:Psi:2}
\Psi^{1,2,34} \Psi^{12,3,4} = \Phi^{2,3,4} \Psi^{1,23,4} \Psi^{1,2,3}. 
\end{equation}

We also set 
$\Assoc(\kk) := \{(\lambda,\Phi) \in \ul{\Assoc}(\kk) | \lambda\in
\kk^\times\}$, and  
$\Pseudo(N,\kk) := \{(a,\lambda,\Phi,\Psi) \in
\ul{\Pseudo}(N,\kk) | (a,\lambda)\in 
(\ZZ/N\ZZ)^\times\times \kk^\times\}$. If $(a,\lambda)\in
(\ZZ/N\ZZ)\times \kk$, 
we also set $\ul{\Pseudo}_{(a,\lambda)}(N,\kk)
:= \{(\Phi,\Psi)  | (a,\lambda,\Phi,\Psi) \in \ul{\Pseudo}(N,\kk)\}$ 
(we write ${\Pseudo}_{(a,\lambda)}(N,\kk)$ instead of 
$\ul{\Pseudo}_{(a,\lambda)}(N,\kk)$ if $(a,\lambda)\in
(\ZZ/N\ZZ)^\times\times\kk^{\times}$). 

Then $\Phi_{\on{KZ}} \in \Assoc_{2\pi\i}(\CC)$, and $(\Phi_{\on{KZ}},
\Psi_{\on{KZ}}) \in \Pseudo_{(-\bar 1,2\pi\i)}(N,\CC)$. 

\begin{proposition}
There is a unique map 
$$
\ul{\Pseudo}(N,\kk) \to \on{Mor}(B_{n}^{1},
\on{exp}(\hat\t_{n,N}(\kk))\rtimes (\ZZ/N\ZZ)^{n}\rtimes\SG_{n}), 
$$
taking $(a,\lambda,\Phi,\Psi)$ to the morphism
$$
\tau\mapsto e^{(\lambda/N)t_0^{01}}s_1^a, \quad 
\sigma_i \mapsto (\Psi^{-1})^{0...i-1,i,i+1}
(i,i+1) e^{(\lambda/2)t(0)^{i,i+1}} \Psi_{\on{KZ}}^{0...i-1,i,i+1}
$$
for $i=1,...,n-1$. When $(a,\lambda)\in (\ZZ/N\ZZ)^\times \times \kk^\times$, 
this morphism induces an isomorphism $B_{n}^{1}(\kk,\varphi_N) \to 
\on{exp}(\hat\t_{n,N}(\kk))\rtimes (\ZZ/N\ZZ)^{n}\rtimes\SG_{n}$, hence 
we get a map 
$$
{\Pseudo}(N,\kk) \to \on{Iso}(B_{n}^{1}(\varphi_{N},\kk),
\on{exp}(\hat\t_{n,N}(\kk))\rtimes (\ZZ/N\ZZ)^{n}\rtimes\SG_{n}). 
$$
\end{proposition}

{\em Proof.} The proof of the fact that the above formulas define a group morphism 
follows the construction of representations of $B_{n}^{1}$ in Section \ref{sec:QRA}. 
One studies the image of the generators $X_{0i}$, $x_{ij}(a)$ to prove that we get 
an isomorphism when $(a,\lambda)\in(\ZZ/N\ZZ)^{\times}\times
\kk^{\times}$. \hfill \qed \medskip 

We will see later that ${\Pseudo}(N,\kk)$ is a torsor; one
can prove that the above map is a morphism of torsors. To prove that a
formality isomorphism $B_{n}^{1}(\varphi_{N},\QQ)\simeq 
\on{exp}(\hat\t_{n,N}(\QQ))\rtimes (\ZZ/N\ZZ)^{n}\rtimes \SG_{n}$
defined over $\QQ$ exists, we will show that 
$\Pseudo_{(a,\lambda)}(N,\QQ)\neq\emptyset$
for any $(a,\lambda)\in(\ZZ/N\ZZ)^{\times}\times\kk^{\times}$, 
which will be based on the action of a group $\on{GTM}(N,\kk)$. 
The two next sections serve as a motivation for the definition of $\on{GTM}(N,\kk)$, 
which is the group of automorphisms of certain braided module categories. 

\section{Pseudotwists over quasibialgebras and module categories} 
\label{sect:2} \label{sec:ps:1}

In this section, we recall the notion of a pseudotwist over a quasi-bialgebra\footnote{All 
(co)algebras are understood to be with (co)unit. 
We use the standard notation for multiple coproducts in coalgebras, 
so $x^{12,34} = (\Delta\otimes \Delta)(x)$, etc. 
If $A$ is an algebra, we denote by $A^\times$ the group of invertible 
elements of $A$.}
and relate it to the notion of module category over a monoidal category. 

\subsection{Pseudotwists}

Let $(A,\Delta_A,\Phi_A)$ be a quasibialgebra (QBA). Recall that this 
means that $A$ is an algebra, $\Delta_A : A \to A^{\otimes 2}$
is a morphism, $\Phi_A \in A^{\otimes 3}$ is invertible, and 
\begin{equation} \label{Delta:Phi}
\forall a\in A, \; (\id_A\otimes\Delta_A)(\Delta_A(a)) = 
\Phi_A (\Delta_A\otimes\id_A)(\Delta_A(a)) \Phi_A^{-1} 
\end{equation}
\begin{equation} \label{pentagon}
\Phi_A^{1,2,34} \Phi_A^{12,3,4} = \Phi_A^{2,3,4} \Phi_A^{1,23,4} \Phi_A^{1,2,3}. 
\end{equation}
The counit $\eps_A$ is also required to satisfy
$(\eps_A \otimes \id_A) \circ \Delta_A
= (\id_A\otimes \eps_A) \circ \Delta_A = \id_A$, 
$\eps_A^{i}(\Phi_A) = 1^{\otimes 2}$ if $i=1,2,3$
(here $\eps_A^1 = \eps_A \otimes \id_A^{\otimes 2}$, etc.; see \cite{Dr1}). 

\begin{definition}
{\it A pseudotwist over $A$ is a triple $(B,\Delta_B,\Psi_B)$ of an 
algebra $B$, an algebra morphism $\Delta_B : B \to B \otimes A$
and an invertible $\Psi_B\in B \otimes A^{\otimes 2}$, such that 
$(\id_{B}\otimes \eps_{A})\circ \Delta_{B}=\id_{B}$ and 
\begin{equation} \label{Delta:Psi}
\forall b\in B, \; 
(\id_B\otimes \Delta_A)(\Delta_{B}(b)) = 
\Psi_B (\Delta_B\otimes \id_A)(\Delta_{B}(b)) \Psi_B^{-1},  
\end{equation}
\begin{equation} \label{pseudotwist} \label{mixed:pent}
\Psi_B^{1,2,34}\Psi_B^{12,3,4} = \Phi_A^{2,3,4} \Psi_B^{1,23,4} \Psi_B^{1,2,3}.      
\end{equation}
Here $\Psi_B^{12,3,4}=(\Delta_B\otimes \on{id}_A^{\otimes 2})(\Psi_B)$, 
$\Psi_B^{1,23,4}=(\on{id}_B \otimes \Delta_A\otimes \on{id}_A)(\Psi_B)$, 
etc.}
\end{definition}

So the notion of pseudotwist is the ``quasi'' version of that of 
comodule-algebra over a bialgebra $A$ (which corresponds to $\Phi_A=1$, 
$\Psi_B=1$). 

\begin{remark}
Similarly to \cite{EN}, Propositions 2.4 and 2.6, one proves that  
if $B \subset A$ is a subalgebra such that $\Delta_{A}(B)\subset B\otimes A$, 
there is a bijection between the following sets:  
(a) invertible elements $\Psi_B \in B\otimes A^{\otimes 2}$, satisfying 
(\ref{Delta:Psi}) and (\ref{pseudotwist}), where 
$\Delta_{B}:= (\Delta_A)_{|B}$; and (b) invertible elements 
$\Ups \in A^{\otimes 2}$, such that $\forall b\in B$, 
$\Upsilon\Delta_{A}(b) = \Delta_{A}(b)\Upsilon$, and $\Upsilon * \Phi_A 
:= (\Ups^{1,23})^{-1} \Phi_A \Ups^{12,3} \Ups^{1,2}
\in B\otimes A^{\otimes 2}$. 
This bijection takes $\Psi_B$ to $\Upsilon := 
(\eps_{A}\otimes \id_{A}^{\otimes 2})(\Psi_B)$, and its inverse takes 
$\Upsilon$ to $\Upsilon * \Phi_A$. 
\end{remark}

\subsection{Module categories over monoidal categories}

Let $(\cC,\otimes,\Phi_\cC,{\bf 1},l,r)$ be a monoidal category. 
Recall that this means that $\cC$ is $\kk$-linear, $\otimes:\cC^{2}\to
\cC$ is a bifunctor, $\Phi_\cC$ is a functorial and additive 
assignment\footnote{If $\cC$ is a category, we denote by $\cC(X,Y)$ the set of 
morphisms from $X$ to $Y$, by $\on{Iso}_{\cC}(X,Y)\subset \cC(X,Y)$ 
the subset of isomorphisms and by $\on{Aut}_{\cC}(X) = \on{Iso}_{\cC}(X,X)$
the set of automorphisms of $X$.} 
$\Phi_{X,Y,Z}\in \on{Iso}_{\cC}(X\otimes (Y\otimes Z),(X\otimes Y)\otimes Z)$, 
satisfying the pentagon identity, ${\bf 1}\in\on{Ob}(\cC)$ and $l,r$ are 
additive and functorial isomorphisms $l_{X}\in 
\on{Iso}_{\cC}({\bf 1}\otimes X,X)$,
$r_{X}\in \on{Iso}_{\cC}(X\otimes {\bf 1},X)$, compatible with $\Phi_\cC$. 

Let $(\cM,\otimes,\Psi_\cM,l')$ be a right module category over $\cC$ 
(\cite{O,CF}). This means that $\cM$ is $\kk$-linear, 
$\otimes : \cM\times\cC\to\cM$ is a bifunctor, $\Psi_\cM$ is a functorial 
and additive assignment $\Psi_{M,X,Y}
\in \on{Iso}_{\cM}(M\otimes (X\otimes Y),(M\otimes X)\otimes Y)$
for $M\in\on{Ob}(M)$, $X,Y\in \on{Ob}(\cC)$, and $r'$ is a functorial and 
additive assignment $r'_{M}\in \on{Iso}_{\cM}(M\otimes {\bf 1},M)$, such that
$\Psi_{M\otimes X,Y,Z} \circ \Psi_{M,X,Y\otimes Z} = 
(\Psi_{M,X,Y}\otimes \id_{Z}) \circ \Psi_{M,X\otimes Y,Z}
\circ (\id_{M}\otimes \Phi_{X,Y,Z})$ (mixed pentagon), and 
$(r_{M\otimes X} \otimes\on{id}_{\bf 1}) \circ \Psi_{M,X,{\bf 1}} 
= \id_{M}\otimes r_{X}$, 
$(r'_{M}\otimes\id_{X})\circ \Psi_{M,{\bf 1},X}=\on{id}_{M}
\otimes l_{X}$. 

If $(B,\Delta_{B},\Psi_B)$ is a pseudotwist over the QBA $(A,\Delta_A,\Phi_A)$,
then $\cM:= \on{Rep}(B)$ is a right module category over $\cC:= \on{Rep}(A)$. 

\subsection{A coherence theorem for module categories}

Let us define the groupoid ${\bf Pa}_n$ of parenthesizations of 
$n$ letters as follows  (here $n\geq 1$). We set 
$\on{Ob}({\bf Pa}_n) := \{$parenthesizations of the word $x_{1}...x_{n}\}$. 
Let us define the morphisms of ${\bf Pa}_n$. If $w',w'',w'''$ are parenthesizations 
of $x_{a}...x_{b}$, $x_{b+1}...x_{c}$ and $x_{c+1}...x_{d}$, if 
$w_{0}$ contains $w'(w''w''')$ and $w_{1}$ is obtained from $w_{0}$
by the replacement $w'(w''w''')\to (w'w'')w''$, then 
there is an invertible morphism $a(w',w'',w''')\in {\bf Pa}_{n}(w_{0},w_{1})$; 
then ${\bf Pa}_n$ is the free groupoid generated by these morphisms and their
inverses. 

Let $\cC$ be a monoidal category. To $X_{1},...,X_{n}\in\on{Ob}(\cC)$,
we attach a representation  of ${\bf Pa}_{n}$; the 
space corresponding to $w\in\on{Ob}({\bf Pa}_{n})$ is
$\otimes_{1\leq i\leq n}^{w} X_{i}$ (meaning that the 
order in which tensor products are taken follows $w$). The image of 
$a(w',w'',w''')$ is $\id_{...}\otimes\Phi_{X',X'',X'''}\otimes \id_{...}$, where 
$X' := \otimes^{w'}_{1\leq i\leq a}X_{i}$, 
$X'': = \otimes^{w''}_{a+1\leq i\leq b}X_{i}$, 
$X''' := \otimes^{w'''}_{b+1\leq i\leq n}X_{i}$.  
Then McLane's theorem says that this representation is trivial 
(i.e., its restriction to each ${\bf Pa}_{n}(w,w)$ is trivial). 

Let now $\cM$ be a module category over $\cC$.
To $M\in\on{Ob}(\cM)$ and $X_{1},...,X_{n}\in\on{Ob}(\cC)$, we attach 
a representation of ${\bf Pa}_{n+1}$ as above (the $\Phi_{...}$ are replaced by 
$\Psi_{...}$ when $X'$ contains $M$). One can check following 
the proof of McLane's theorem that this representation is trivial. 
 
\section{Quasi-reflection algebras and braided module categories} 
\label{sect:3} \label{sec:QRA}

In this section, we introduce the notion of quasi-reflection algebras (QRA's) over a 
quasitriangular quasi-bialgebra (QTQBA). The corresponding notion is that of 
braided module category over a monoidal category. We show that such 
objects give rise to representations of $B_{n}^{1}$. 

\subsection{QRA's over QTQBA's}

Let $(A,\Delta_A,R,\Phi_A,\eps_{A})$ be a quasitriangular quasibialgebra. 
This means that $A$ is an algebra, $\Delta_A : A \to A^{\otimes 2}$
is a morphism, $R\in A^{\otimes 2}$ and $\Phi_A \in A^{\otimes 3}$
are invertible, and 
$$
\forall a\in A, \; 
\Delta_A^{2,1}(a) = R \Delta_A(a) R^{-1}, 
\; (\id_A\otimes\Delta_A) \circ \Delta_A (a) = 
\Phi_A (\Delta_A\otimes\id_A) \circ \Delta_A(a) \Phi_A^{-1} 
$$  
$$
\Phi_A^{1,2,34} \Phi_A^{12,3,4} = \Phi_A^{2,3,4} \Phi_A^{1,23,4} \Phi_A^{1,2,3}
$$
$$
R^{12,3} = \Phi_A^{3,1,2} R^{1,3} (\Phi_A^{1,3,2})^{-1}
 R^{2,3} \Phi_A, \; 
R^{1,23} = (\Phi_A^{2,3,1})^{-1} R^{1,3} \Phi_A^{2,1,3} 
R^{1,2} \Phi_A^{-1}. 
$$
These conditions are the pentagon and the two hexagon equations. 
The counit $\eps_{A}$ satisfies $(\eps_{A} \otimes \id_{A}) \circ \Delta_A
= (\id_{A}\otimes \eps_{A}) \circ \Delta_A = \id_A$, 
$\eps_{A}^{i}(R) = 1$ for $i=1,2$, 
$\eps_{A}^{i}(\Phi_A) = 1^{\otimes 2}$ if $i=1,2,3$. 

\begin{definition}
{\it A quasi-reflection algebra (QRA) over $A$ is a quadruple 
$(B,\Delta_B,E,\Psi_B)$, where $(B,\Delta_B,\Psi_B)$ is a pseudotwist 
over $A$, and $E$ satisfies 
\begin{equation} \label{octogon}
(\Delta_{B} \otimes \id_{A})(E) = \Psi_B^{-1} R^{3,2} \Psi_B^{1,3,2}
E^{1,3} (\Psi_B^{1,3,2})^{-1} R^{2,3} \Psi_B.  
\end{equation}}
\end{definition}

We will call the relations (\ref{mixed:pent}) and (\ref{octogon})
the mixed pentagon and octogon 
relations.  

When $\Phi_A=1$, $A$ is a quasitriangular bialgebra, and when $\Psi_B=1$, 
we say that $B$ is a reflection algebra (RA) over $A$. 

\begin{remark} If $A$ is a QTQBA and $B\subset A$ be a subalgebra 
such that $\Delta_A(B) \subset B\otimes A$, 
and set $\Delta_B:= (\Delta_A)_{|B}$. One can show that there is a bijection 
between the following sets: 
(a) pairs $(E,\Psi_B)$, where $E\in B\otimes A$ and $\Psi_B \in B \otimes 
A^{\otimes 2}$ are invertible, s.t. $(B,\Delta_B,E,\Psi_B)$ is a QRA over $A$; 
(b) pairs $(\Sigma,\Upsilon)$, where $\Sigma\in A$ and $\Upsilon
\in A^{\otimes 2}$ are invertible, such that 
$\forall b\in B, \; \Sigma b=b\Sigma$, $\Upsilon \Delta_A(b)
=\Delta_A(b)\Upsilon$, 
$(\Upsilon^{1,23})^{-1} \Phi_A \Upsilon^{12,3} \Upsilon^{1,2}
\in B \otimes A^{\otimes 2}$ and $\Ups^{-1} R^{2,1} \Ups^{2,1} \Sigma^2 
(\Ups^{2,1})^{-1} R \Ups \in B \otimes A$.  
This bijection takes $(E,\Psi_B)$ to $(\Sigma(E),\Upsilon(\Psi_B))$, 
where $\Sigma(E) := (\eps_{A}\otimes \id_{A})(E)$ and $\Upsilon(\Psi_B) := 
(\eps_{A}\otimes \id_{A}^{\otimes 2})(\Psi_B)$. The inverse bijection takes 
$(\Sigma,\Upsilon)$ to $(E(\Sigma,\Ups),\Psi(\Upsilon))$, where 
$$
E(\Sigma,\Ups) := \Ups^{-1} R^{2,1} \Ups^{2,1} \Sigma^2 
(\Ups^{2,1})^{-1} R \Ups, \; 
\Psi(\Upsilon) := 
(\Upsilon^{1,23})^{-1} \Phi_A \Upsilon^{12,3} \Upsilon^{1,2}.
$$ 
When $\Phi_A=1$, this restricts to a bijection between the sets of pairs
$(E,\Psi_B)$ and $(\Sigma,\Upsilon)$ as above, where $\Psi_B=1$ and $\Upsilon=1$. 
\end{remark}

\begin{remark} \label{rem:1:5} If $B$ is a RA over the QTBA $A$, then 
$(E,R)$ satisfy the reflection equation
\begin{equation} \label{RE}
E^{1,2}R^{3,2} E^{1,3} R^{2,3}
= R^{3,2} E^{1,3} R^{2,3} E^{1,2} . 
\end{equation}
In particular, if $(\rho,V)$ is a $B$-module, then 
$K:= (\rho\otimes \id)(E)$ satisfies the reflection equation
$R^{2,1}K^2 R^{1,2} K^1 = K^1 R^{2,1} K^2 R^{1,2}$.
Reflection equation algebras (REA) are defined using this 
equation (\cite{FRT}), $R$ being fixed; REA's yield 
representations of $B_n^1$. Applications of REA's include the topology
of handlebodies (\cite{KS}) and scattering theory with reflection 
(\cite{Ch}). The analogy between both situations is that winding around 
the first fixed strand is analogous to the reflection of a particle on a wall. 
In \cite{DKM}, universal versions of REA's were studied; they rely on the 
datum of a ``semiuniversal" object $K\in \on{End}(V) \otimes B$, where 
$V$ is a $A$-module, $A$ is a Hopf algebra and $B$ is an algebra. The data 
of a RA $B$ over a QTBA $A$ and a $A$-module $V$ yield such a $K$.  
\end{remark} 

\begin{example} If we set $B := A$, $E:= R^{2,1}R$, 
$\Psi_B := \Phi_A$, then $(B,E,\Psi_B)$ is a QRA over $A$. 
It corresponds to $\Sigma = 1$, $\Upsilon = 1^{\otimes 2}$.  
\end{example}

\begin{example} \label{ex:g:Gamma}
Let $\g$ be a Lie algebra, let $\Gamma$ be a group 
acting on $\g$ by automorphisms. 
Set $A := U(\g) \rtimes \Gamma$. Define $\Delta_A : A \to A^{\otimes 2}$
as the unique algebra morphism such that $\Delta_A(x) = 
x\otimes 1 + 1 \otimes x$ for $x\in \g$, $\Delta_A(\gamma)
= \gamma \otimes \gamma$ for $\gamma\in\Gamma$. Set $R :=1^{\otimes 2}$, 
then $(A,\Delta_A,R)$ is a QTBA. 
Fix $\sigma\in \Gamma$ and a Lie subalgebra $\l \subset \g^\sigma$. 
Set $B:= U(\l)$, $E := 1\otimes \sigma$. Then $(B,E)$ is a RA 
over $A$. We will later introduce a deformation of 
this example, where $\Gamma = \ZZ / N\ZZ$. 
\end{example}

\subsection{Twists}

Let $(A,\Delta_{A},R,\Phi_A)$ be a QTQBA. To an invertible 
$F \in A^{\otimes 2}$, one associates the twisted QTQBA 
$^F\!A := (A, ^F\!\Delta, ^F\!R, ^F\!\Phi_A)$, where 
$$
^F\!\Delta_{A}(a) = F\Delta_{A}(a)F^{-1}, \; 
^F\!R = F^{2,1} R F^{-1}, \; 
^F\!\Phi_A = F^{2,3} F^{1,23} \Phi_A (F^{1,2} F^{12,3})^{-1}. 
$$

If now $(B,E,\Psi_B)$ is a QRA over $A$ and $G\in B\otimes A$
is invertible, set 
$$
^{G}\!\Delta_{B}(b)=G\Delta_{B}(b)G^{-1}, \; 
^{G}E := GEG^{-1}, \; ^{F,G}\Psi_B:= F^{2,3}
(\on{id}_{B}\otimes\Delta_{A})(G) \Psi_B 
(G^{1,2}(\Delta_{B}\otimes\on{id}_{A})(G))^{-1}.
$$

Then $(B,^{G}\!\Delta_{B},^{G}\!E,^{F,G}\!\Psi_B)$ is a QRA over 
$(A,^{F}\!\Delta_{A},^{F}\!R,^{F}\!\Phi_A)$; we call it 
the twist of $B$ by $(F,G)$. 

\subsection{Braided module categories over braided monoidal categories}

Let $\cC$ be a braided monoidal category. This means that 
$\cC$ is a monoidal category, equipped with a braiding
$\beta_{X,Y}\in \on{Iso}_{\cC}(X\otimes Y,Y\otimes X)$
(where 
$X,Y\in\on{Ob}(\cC)$), $\beta_{X\otimes Y,Z} = (\beta_{X,Z}\otimes 
\id_{Y}) \circ (\id_{X}\otimes \beta_{Y,Z})$, 
$\beta_{X,Y\otimes Z} = (\id_{Y}\otimes\beta_{X,Z})
\circ (\beta_{X,Y}\otimes \id_{Z})$, $\beta_{{\bf 1},X}
=\beta_{X,{\bf 1}}=\id_{X}$ (we omitted the 
associativity constraints). 

We define a braided module category over $\cC$ to be a 
module category $\cM$ over $\cC$, equipped with an additive and
functorial $\eps_{M,X}\in \on{Aut}_{\cM}(M\otimes X)$, such that 
$$
\eps_{M\otimes X,Y} = \Psi^{-1}_{M,X,Y} \circ 
(\id_{M}\otimes\beta_{Y,X}) \circ \Psi_{M,Y,X} \circ
(\eps_{M,X}\otimes\id_{Y}) \circ \Psi_{M,Y,X}^{-1} \circ 
(\id_{M}\otimes\beta_{X,Y}) \circ \Psi_{M,X,Y}.
$$
Then if $A$ is a QTQBA and $B$ is a QRA over $A$, then $\cM:= \on{Rep}(B)$
is a braided module category over $\cC:= \on{Rep}(A)$; the braided monoidal 
structure on $\cC$ is  
$\beta_{X,Y}:= \sigma_{X,Y} \circ (\rho_{X}\otimes \rho_{Y})(R)$
and $\Phi_{X,Y,Z}:= (\rho_X\otimes \rho_Y \otimes \rho_Z)(\Phi_A)$
for $(X,\rho_{X}),...\in\on{Ob}(\on{Rep}(A))$ (here $\sigma_{X,Y}
\in\on{Hom}_{\kk}(X\otimes Y,Y\otimes X)$ is the exchange operator)
and the braided module structure on $\cM$ is $\eps_{M,X}:= (\rho_{M}\otimes 
\rho_{X})(E)$, $\Psi_{M,X,Y}:= (\rho_M \otimes \rho_X \otimes \rho_Y)(\Psi_B)$
for $(M,\rho_M)\in\on{Ob}(\on{Rep}(B))$. 

\subsection{Braided module categories and representations of $B_{n}^{1}$}

Recall that if $\cC$ is a strict braided tensor category, and 
$X\in\on{Ob}(\cC)$, then $\beta$ gives rise to a group morphism 
$B_{n}\to \on{Aut}_{\cC}(X^{\otimes n})$ for any 
$X\in\on{Ob}(\cC)$. If $\sigma_1,...,\sigma_{n-1}$ are the standard 
generators of $B_n$, this morphism is $\sigma_i\mapsto 
\on{id}_{X^{\otimes i-1}}\otimes \beta_{X,X} \otimes 
\on{id}_{X^{\otimes n-i-1}}$. 

If now $\cM$ is a strict braided module category over the strict braided 
monoidal category $\cC$ (i.e., the two natural bifunctors $\cC^{3}\to\cC$ 
coincides, as well as the two bifunctors $\cM\times\cC^{2}\to \cM$, 
and the constraints $\Phi_\cC$ and $\Psi_\cM$ are the identity), if 
$M\in\on{Ob}(\cM)$ and $X\in\on{Ob}(\cC)$, then 
$\eps,\beta$ give rise to a compatible 
group morphism $B_{n}^{1}\to \on{Aut}_{\cM}(M\otimes X^{\otimes n})$ 
for $M\in\on{Ob}(\cM)$, 
$X\in\on{Ob}(\cC)$, by $\tau\mapsto \eps_{M,X}\otimes 
\on{id}_{X^{\otimes n-1}}$, $\sigma_i\mapsto 
\on{id}_{M\otimes X^{\otimes i-1}} \otimes \beta_{X,X} 
\otimes \on{id}_{X^{\otimes n-i-1}}$. 

The fact that $(\tau\sigma_{1})^{2}=(\sigma_{1}\tau)^{2}$ is preserved
follows from the identity $(\eps_{M,X}\otimes \on{id}_X) \circ 
\eps_{M\otimes X,X}=\eps_{M\otimes X,X}\circ
(\eps_{M,X}\otimes\on{id}_X)$, 
which follows from the functoriality of $\eps$.  

In \cite{BN}, the groupoid ${\bf PaB}_n$ was defined by 
$\on{Ob}({\bf PaB}_{n}) = \{(\sigma,w) | \sigma\in\SG_{n}, w$ is a 
parenthesization of $x_{\sigma(1)}... x_{\sigma(n)}\}$ and 
${\bf PaB}_{n}((\sigma,w),(\sigma',w')) = \{$braids relating 
$x_{\sigma(1)}...x_{\sigma(n)}$ and $x_{\sigma'(1)}...x_{\sigma'(n)}\}$.  
Combining the proof of braid representations in the strict case with 
McLane's theorem, one shows that to any braided monoidal category $\cC$
and any $X_{1},...,X_{n}\in\on{Ob}(\cC)$, there corresponds a 
representation of ${\bf PaB}_{n}$, where the object corresponding to 
$(\sigma,w)$ is $\otimes^{w}_{1\leq i\leq n}X_{i}$. 
In particular, we get a morphism $B_{n}\to \on{Aut}_{\cC}(X^{\otimes n})$, 
where $X^{\otimes n} = ((X\otimes X)\otimes ...)\otimes X$. 
When $\cC = \on{Rep}(A)$, where $A$ is a QTQBA, this morphism 
is induced by a morphism $B_{n}\to (A^{\otimes n})^{\times}\rtimes\SG_{n}$, 
such that 
$$
\begin{matrix}
B_{n}^{1} & \to & (A^{\otimes n})^{\times}\rtimes \SG_{n}
\\
  &\searrow & \downarrow\\
  & & \SG_{n}
\end{matrix}
$$
commutes; it is given by 
$$\sigma_{i}\mapsto (
\Phi_A^{((12)..).i-1,i,i+1})^{-1} 
\cdot (i,i+1)R^{i,i+1}\cdot 
\Phi_A^{((12)...)i-1,i,i+1}
$$
(the generators of $B_{n}$ are $\sigma_{1},...,\sigma_{n-1}$); here
$\Phi_A^{((12)..).i-1,i,i+1}
=(\Delta_{A}^{((12)...)i-1}
\otimes\id_{A}^{\otimes 2})(\Phi_A)$, 
and $\Delta_{A}^{((12)...)i}
=(\Delta_{A}\otimes\id_{A}^{\otimes i-2})\circ
... \circ \Delta_{A}$. 

We define similarly the groupoid ${\bf PaB}_{n}^{1}$ by 
$\on{Ob}({\bf PaB}_{n}) = \{(\sigma,w) | \sigma\in\SG_{n}, w$ is a 
parenthesization of $x_{0}x_{\sigma(1)}... x_{\sigma(n)}\}$ and 
${\bf PaB}_{n}((\sigma,w),(\sigma',w')) = \{$braids relating 
$x_{0}x_{\sigma(1)}...x_{\sigma(n)}$ and 
$x_{0}x_{\sigma'(1)}...x_{\sigma'(n)}\}$. As above, to the data of 
a braided module category $\cM$ over the braided monoidal category 
$\cC$, and to $M\in\on{Ob}(\cM)$, $X_{1},...,X_{n}\in\on{Ob}(\cC)$, 
there corresponds a representation of ${\bf PaB}_{n}^{1}$, where 
to $(\sigma,w)$, we attach $\otimes^{w}_{0\leq i\leq n}X_{\sigma(i)}$, 
where  $X_{0}:= M$ and $\sigma(0)=0$. In particular, we get a 
morphism $B_{n}^{1}\to \on{Aut}_{\cM}(M\otimes X^{\otimes n})$, 
where $M\otimes X^{\otimes n} = ((M\otimes X)\otimes ... )\otimes X$. 
When $\cC = \on{Rep}(A)$, $\cM = \on{Rep}(B)$, $A$ is a QTQBA and 
$B$ is a QRA over $A$, this morphism is induced by a morphism 
$B_{n}^{1}\to (B\otimes A)^{\otimes n}\rtimes \SG_{n}$, 
such that 
$$
\begin{matrix}
B_{n}^{1} & \to & (B\otimes A^{\otimes n})^{\times}\rtimes \SG_{n}
\\
  &\searrow & \downarrow\\
  & & \SG_{n}
\end{matrix}
$$
commutes; this morphism is given by 
$$
\tau\mapsto E^{0,1}, \quad  \sigma_{i}\mapsto 
 (\Psi_B^{((12)...)i-1,i,i+1})^{-1}  \cdot (i,i+1)R^{i,i+1} \cdot
 (\Psi_B^{((12)...)i-1,i,i+1})^{-1}, 
$$ 
where
$\Psi_B^{((12)...)i-1,i,i+1}=(\Delta_{B}^{((12)...)i-1}
\otimes\id_{A}^{\otimes 2})(\Psi_B)$, 
and $\Delta_{B}^{((12)...)i}
=(\Delta_{B}\otimes\id_{A}^{\otimes i-2})\circ
... \circ \Delta_{B}$. 

\begin{remark}
If $g\geq 1$, let $B_n^g:= B_{g+n} \times_{\SG_{g+n}}\SG_{n}$, where 
$\SG_n\subset \SG_{g+n}$ is the subgroup of permutations of the 
last $n$ elements. A presentation of $B_n^g$ can be found in 
\cite{Sos,V}; generators are $\tau_1,...,\tau_g,\sigma_1,...,\sigma_g$
and relations are (\ref{std:braid}), 
$\tau_k \sigma_i = \sigma_i \tau_k$ ($k\in [1,g]$, $i\in[2,n]$), 
$\tau_k\sigma_1\tau_k\sigma_1 =\sigma_1 \tau_k \sigma_1\tau_k$
($k \in [1,g]$), $\sigma_1\tau_{\ell} \sigma_1^{-1} \tau_k = 
\tau_k\sigma_1 \tau_{\ell} \sigma_1^{-1}$ ($1\leq k < l \leq g$). 
The inclusion $B_n^g \subset B_{n+g}$ is induced by $\sigma_i \mapsto
\sigma_{g+i}$, $\tau_i \mapsto \sigma_g^{-1} \cdots \sigma_{i+1}^{-1} 
\sigma_i^2 \sigma_{i+1} \cdots \sigma_g$. 

If $\cC$ is a strict braided monoidal category, $\cM$ is a 
strict monoidal category and $F : \cC \to \cM$ is a monoidal functor, 
then $\cM$ is a right strict module category over $\cC$ (the 
bifunctor $\otimes : \cM\times \cC\to \cM$ is $(M,X)\mapsto 
M \otimes F(X)$). Then a braided module category structure of  
$\cM$ over $\cC$ gives rise to group morphisms 
$B_n^g\to \on{Aut}_{\cM}(\otimes_{i=1}^g M_i 
\otimes X^{\otimes n})$. 
We describe these morphisms when $\cC = \on{Rep}(A)$, 
$\cM = \on{Rep}(B)$, $A$ is a QTBA, $B\subset A$ is a subbialgebra. 
The morphism $B_n^g \to (B^{\otimes g} \otimes A^{\otimes n})^\times
\rtimes \SG_n$ is given by $\sigma_i\mapsto (\overline i, \overline{i+1}) 
\circ R^{\overline i,\overline{i+1}}$, $\tau_i\mapsto 
(E^{i+1...g,\bar 1})^{-1} E^{i...g,\bar 1}$ (we set 
$E^{\emptyset,\bar 1} = 1$ and $\overline{i} = i+g$). 
\end{remark}

\subsection{Specializations of the $\t_{n,N}$ and examples for QRA's over 
QTQBA's} \label{rem:5:1} 

If $\g$ is a Lie algebra over $\kk$, $t_\g\in S^2(\g)^\g$, and $\sigma \in 
\on{Aut}(\g,t_\g)$ is such that $\sigma^N = \id$, then $(\g,t_\g,\sigma)$ 
gives rise to a representation $T_{n,N} \to U(\l) \otimes 
(U(\g)\rtimes \ZZ/N\ZZ)^{\otimes n}$
(here $\l = \g^\sigma$), as follows. Set $\u = \on{Im}(\sigma-\id)$, then 
$\g = \l\oplus \u$ and $t_\g=t_\l\oplus t_\u$, where $t_\l\in S^2(\l)$
and $t_\u\in S^2(\u)$. The representation is $t_0^{0i}\mapsto 
\hbar N (t_\l^{0i} + {1\over 2} t_\l^{ii})$,  
$t(a)^{ij} \mapsto \hbar (\sigma^a\otimes \id)(t_\g)^{ij}$, and 
$s_i\mapsto \ul\sigma^i$ (here $\hbar$ is a formal parameter and 
$\ul\sigma$ is the generator of $\ZZ/N\ZZ\subset 
U(\g)\rtimes\ZZ/N\ZZ$). 
$(\g,t_\g,\sigma)$ also gives rise to a representation 
$\t_n \to U(\g)^{\otimes n}$, $t^{ij} \mapsto \hbar t_\g^{ij}$. 
This collection of representations is compatible with (a) the 
insertion-coproduct maps $U(\l) \otimes (U(\g)\rtimes \ZZ/N\ZZ)^{\otimes n}
\to U(\l) \otimes (U(\g)\rtimes \ZZ/N\ZZ)^{\otimes m}$ induced 
by a partially defined function $f : [0,m]\to [0,n]$ with $f(0)=0$, 
and (b) the insertion-coproduct maps $U(\g)^{\otimes n}
\to U(\l) \otimes U(\g)^{\otimes m}$ induced by a partially defined 
function $f : [1,m]\to [1,n]$.    

If $(a,\lambda,\Phi,\Psi) \in \ul\Pseudo(N,\kk)$, then $(
U(\g) \rtimes \ZZ/N\ZZ[[\hbar]], R = e^{\hbar (\lambda/2) t_\g}, 
\Phi(\hbar t_\g^{12},\hbar t_\g^{23}))$ is a QTQBA and 
$$
(U(\l)[[\hbar]], E = e^{\hbar (\lambda/N) (t_\l^{12} + {1\over 2} t_\l^{22})}
(1\otimes \underline\sigma^a),
\Psi(\hbar (t_\l^{12} + {1\over 2} t_\l^{22}), 
\hbar t_\g^{23},\ldots, \hbar (1\otimes \sigma^{N-1})(t_\g)^{23}))
$$
is a QRA over it. 

\section{The semigroups $\ul{\on{GTM}}$ and $\wh{\ul{\on{GTM}}}$} 
\label{sect:6} \label{sec:GTM}

\subsection{The semigroup $\ul{\on{GTM}}$}

Let $(A,\Delta_{A},R,\Phi_A)$ be a QTQBA and let $(B,\Delta_{B},E,\Psi_B)$ be a 
QRA over $A$. Let $P_n \subset B_n$ be the pure braid group.
We have group morphisms $P_n \to (A^{\otimes n})^\times$ and 
$P_n\to (B\otimes A^{\otimes n-1})^\times$. Let us describe them 
explicitly when $n=3$. 

We have an isomorphism $P_3\simeq P_3^0\times\ZZ\simeq F_2 \times\ZZ$, 
where the center $\ZZ$ is generated by $(\sigma_1\sigma_2)^3 
= (\sigma_1^2\sigma_2)^2=x_{12}x_{13}x_{23}$ 
($x_{12}=\sigma_1^2$, $x_{23}=\sigma_2^2$, $x_{13} = 
\sigma_2\sigma_1^2\sigma_2^{-1}$) and $P_3^0 = \on{Ker}(P_3 \to 
P_2 \simeq \ZZ)$ where $P_3\to P_2$ is the erasing of the strand $2$, 
i.e., $x_{12},x_{23}\mapsto 0$, $x_{13}\mapsto 1$; $P_3^0\simeq F_2$ 
is freely generated by $x:= x_{12}$ and $y:= x_{23}$. 

The morphism $\pi_A : P_3 \to (A^{\otimes 3})^\times$ is 
given by 
$$
\sigma_1^2 \mapsto R^{2,1}R^{1,2}, \quad 
\sigma_2^2 \mapsto \Phi_A^{-1} R^{3,2}R^{2,3} \Phi_A, \quad 
(\sigma_1\sigma_2)^3 \mapsto 
R^{2,1} R^{1,2} \Phi_A R^{3,2} (\Phi_A^{1,3,2})^{-1}
R^{3,1} R^{1,3} \Phi_A^{1,3,2} R^{2,3} \Phi_A^{-1}. 
$$
The morphism $\pi_{BA} : P_3 \to (B\otimes A^{\otimes 2})^\times$ is given by 
$$
\sigma_1^2 \mapsto E^{1,2}, \quad 
\sigma_2^2 \mapsto \Psi_B^{-1} R^{3,2} R^{2,3} \Psi_B, \quad 
(\sigma_1\sigma_2)^3 \mapsto 
E^{1,2} \Psi_B^{1,2,3} R^{3,2} (\Psi_B^{1,3,2})^{-1}
E^{1,3} \Psi_B^{1,3,2} R^{2,3} (\Psi_B^{1,2,3})^{-1}. 
$$

Let $\lambda\in 2\ZZ+1$, $\mu\in\ZZ$ and 
$f:= f(\sigma_1^2,\sigma_2^2)(\sigma_1\sigma_2)^{3a}$, 
$g:= g(\sigma_1^2,\sigma_2^2)(\sigma_1\sigma_2)^{3a'}$ be elements of 
$P_3$. Set 
$$
\tilde R := R(R^{2,1}R)^{(\lambda-1)/2}, \quad 
\tilde\Phi_A := \Phi_A \pi_A(f), \quad 
\tilde E := E^{\mu}, \quad 
\tilde\Psi_B := \Psi_B \pi_{BA}(g). 
$$
Let us write under which conditions   
$(A,\Delta_{A},\tilde R,\tilde\Phi_A)$ is a QTQBA and 
$(B,\Delta_{B},\tilde E,\tilde\Psi_B)$
is a QRA over it. According to \cite{Dr2}, $(A,\Delta_{A},\tilde R,\tilde\Phi_A)$ 
is a QTQBA if $a=0$ and  
\begin{equation} \label{cond:f:1}
f(y,x) = f(x,y)^{-1}, \quad 
f(x_3,x_1) x_3^{(\lambda-1)/2} f(x_2,x_3) x_2^{(\lambda-1)/2} f(x_1,x_2) 
x_1^{(\lambda-1)/2} = 1
\end{equation}
where $x_1,x_2,x_3$ are three variables subject to the only relation 
$x_1x_2x_3=1$,  and  
\begin{equation} \label{cond:f:2}
f(x_{12},x_{23}x_{24}) f(x_{13}x_{23},x_{34}) = 
f(x_{23},x_{34}) f(x_{12}x_{13},x_{24}x_{34}) f(x_{12},x_{23})
\end{equation}
holds in $P_4$, where $x_{ij}$ are the Artin generators
$$
x_{ij} = (\sigma_{j-2}\cdots \sigma_i)^{-1}
\sigma_{j-1}^2 (\sigma_{j-2}\cdots \sigma_i)
= (\sigma_{j-1} \cdots \sigma_{i+1}) \sigma_i^2 
(\sigma_{j-1} \cdots \sigma_{i+1})^{-1}
\quad (1\leq i < j \leq 4). 
$$

\begin{lemma}
If these conditions are met, $(B,\Delta_{B},\tilde E,\tilde\Psi_B)$ is a QRA over 
$(A,\Delta_{A},\tilde R,\tilde\Phi_A)$ if $a'=0$, 
\begin{equation} \label{cond:g:1}
g(x,y)^{-1}y^{{{\lambda-1}\over 2}} g(x^{-1}y^{-1},y)(x^{-1}y^{-1})^{\mu}
g(x^{-1}y^{-1})^{-1}y^{{{\lambda+1}\over 2}}g(x,y)x^{\mu}=1
\end{equation} 
(identity in the free group generated by $x,y$), and 
\begin{equation} \label{cond:g:2}
g(x_{01},x_{12}x_{13}) g(x_{02}x_{12},x_{23}) = 
f(x_{12},x_{23}) g(x_{01}x_{02},x_{13}x_{23}) 
g(x_{01},x_{12}) 
\end{equation}
holds in $P_4$ (with Artin generators $x_{ij}$, $0\leq i<j\leq 3$). 
\end{lemma}

{\em Proof.} The identity implying that 
$(\tilde\Phi_A,\tilde\Psi_B)$ satisfies (\ref{mixed:pent}) is 
\begin{align} \label{aux:id}
& g(x_{12},x_{23}x_{24}) z_{1,2,34}^{a'} g(x_{13}x_{23},x_{34}) z_{12,3,4}^{a'}
\\ & \nonumber 
= f(x_{23},x_{24}) g(x_{12}x_{13},x_{24}x_{34}) z_{1,23,4}^{a'}
g(x_{12},x_{13}) z_{1,2,3}^{a'}, 
\end{align}
where $z_{1,2,3} = (\sigma_1\sigma_2)^3 = x_{23}x_{12}x_{13}$, 
$z_{2,3,4} = (\sigma_2\sigma_3)^3 = x_{24}x_{34}x_{23}$, 
$z_{12,3,4} = x_{12}^{-1}z$, $x_{1,23,4}=x_{23}^{-1}z$, $x_{1,2,34} = x_{34}^{-1}z$, 
where $z := (\sigma_3\sigma_2\sigma_1)^3 = x_{14} x_{24} x_{34}
x_{13} x_{23} x_{12}$ generates the center of $P_4$. There is a unique 
morphism $P_{4}\to \ZZ$ with $x_{14}\mapsto 1$, $x_{ij}\mapsto 0$ for 
$1\leq i<j\leq 4$, $(i,j)\neq (1,4)$ (it corresponds to the removal of strands $2,3$). 
The image of (\ref{aux:id}) by this morphism yields $2a'=a'$, hence $a'=0$. 
Then (\ref{aux:id}) implies (\ref{cond:g:2}).  (\ref{cond:g:1}) is then the 
condition for $(\tilde E,\tilde\Psi_B)$ to satisfy (\ref{octogon}). 
\hfill \qed \medskip 

\begin{remark}
(\ref{cond:g:1}) 
implies that 
$\theta_{(\lambda,\mu,g)}\in \on{End}(F_{2})$ ($F_{2}$ is the free group 
generated by $x,y$) defined by $x\mapsto g(x,y)x^{\mu}g(x,y)^{-1}$, 
$y\mapsto y^{\lambda}$ also takes $x^{-1}y^{-1}$ to a conjugate of a power of 
this element, namely $x^{-1}y^{-1}\mapsto y^{{{\lambda-1}\over 2}}
g(x^{-1}y^{-1},y)\cdot (x^{-1}y^{-1})^{\mu} \cdot
(y^{{{\lambda-1}\over 2}} g(x^{-1}y^{-1},y))^{-1}$. 
 \hfill \qed \medskip 
\end{remark} 
 
Define $\ul{\on{GTM}}$ as the set of $(\lambda,\mu,f,g)\in (2\ZZ+1)
\times \ZZ \times F_2\times F_2$, satisfying (\ref{cond:f:1}), 
(\ref{cond:f:2}), (\ref{cond:g:1}) and (\ref{cond:g:2}). Then $\ul{\on{GTM}}$ 
has a semigroup structure, defined by 
$$
(\lambda_1,\mu_1,f_1,g_1)* (\lambda_2,\mu_2,f_2,g_2)
= (\lambda,\mu,f,g), 
$$
where
\begin{equation} \label{comp:f}
\lambda = \lambda_{1}\lambda_{2}, \quad 
f(x,y) = f_1\big( f_2(x,y) x^{\lambda_{2}} f_2(x,y)^{-1},y^{\lambda_2}\big)
f_2(x,y),  
\end{equation}
$$
\mu = \mu_{1}\mu_{2}, \quad 
g(x,y) = g_1\big( g_2(x,y) x^{\mu_2} g_2(x,y)^{-1},y^{\lambda_2}\big)
g_2(x,y).   
$$
Recall that $\ul{\on{GT}}$ was defined in \cite{Dr2} as the semigroup 
of all pairs $(f,\lambda)$ satisfying (\ref{cond:f:1}), (\ref{cond:f:2}). 

We have a natural morphism $\ul{\on{GTM}} \to \ul{\on{GT}}$, taking 
$(\lambda,\mu,f,g)$ to $(\lambda,f)$. 
The semigroup $\ul{\on{GTM}}$ acts on the set of pairs $(\cC,\cM)$
of a braided monoidal category $\cC$ and a braided module category $\cM$
over $\cC$ (GTM stands for ``module category version of GT"). 

According to \cite{Dr2}, Proposition 4.1, $\ul{\on{GT}}$ 
consists of the pairs $(\lambda,f)$, where $\lambda = \pm 1$ 
and $f(x,y) = 1$; so $\ul{\on{GT}} \simeq \ZZ/2\ZZ$, where $\ZZ/2\ZZ$
acts by $(A,\Delta_{A},R,\Phi_A) \mapsto (A,\Delta_{A},(R^{2,1})^{-1},\Phi_A)$. 
In the same way, we have: 

\begin{proposition} \label{prop:GTM}
$\ul{\on{GTM}}$ consists of the  
$(\lambda,\mu,f,g)$, where $\lambda = \mu = \pm 1$, 
$f(x,y) = 1$, $g(x,y) = y^s$ ($s\in \ZZ$), so $\ul{\on{GTM}}\simeq
\ZZ\rtimes \ZZ/2\ZZ$ (the action of $\ZZ/2\ZZ$ on $\ZZ$ is nontrivial); 
$\ZZ/2\ZZ$ acts by $(A,\Delta_{A},R,\Phi_A;B,\Delta_{B},E,\Psi_B)
\mapsto (A,\Delta_{A},(R^{2,1})^{-1},\Phi_A;B,\Delta_{B},E^{-1},\Psi_B)$ . 
\end{proposition}

{\em Proof.} We already have $\lambda = \pm 1$ and $f(x,y) = 1$. 
Apply the morphism $P_{4}\to P_{3}$, $x_{ij}\mapsto x_{ij}$ if $i\geq 1$, 
$x_{0j}\mapsto 1$ corresponding to removing the first strand, 
to (\ref{cond:g:2}). We get $g(x_{12},x_{23}) = (x_{12}x_{13})^{-s}
f(x_{12},x_{23})(x_{13}x_{23})^{s}x_{12}^{s}$, where $s\in\ZZ$ is such that 
$g(1,y)=y^{s}$. Since $x_{12}x_{13}x_{23}$ is central in $P_{3}$, 
we get $g(x,y)=y^{s}f(x,y)=y^s$. Then $g(x,y)$ satisfies (\ref{cond:g:2})
and (\ref{cond:g:1}) is satisfied iff 
\begin{equation} \label{temp}
(xyx^{-1})^{(\lambda-1)/2} (y^{-1} x^{-1})^{\mu-1} y^{(\lambda-1)/2} 
x^{\mu-1}=1,    
\end{equation}
which implies (as the degree in $y$ of the l.h.s. should vanish) 
$\lambda=\mu$, and since $\lambda=\pm 1$, we have necessarily 
$\lambda=\mu=\pm1$. Conversely, (\ref{temp}) is satisfied when 
$\lambda=\mu=\pm 1$.  \hfill \qed \medskip 

While the solutions of (\ref{cond:f:1}), (\ref{cond:f:2}), (\ref{cond:g:1}), 
(\ref{cond:g:2}) act on (general) braided module
categories, solutions of these equations in completions of the pure
braid groups act on braided module categories with additional properties. 

\subsection{The semigroup $\wh{\ul{\on{GTM}}}$}

If $G$ is a group (other that $\on{GT}$ or $\on{GTM}$), we denote by 
$\wh G$ (or $G^\wedge$) the profinite completion of $G$; if $\phi:G\to H$ 
is a group morphism, we denote by $\wh\phi:\wh G\to \wh H$ the induced 
morphism. 

Recall the definition of $\wh{\ul{\on{GT}}}$ (\cite{Dr2}). 
The conditions expressing that $(\lambda,f)\in\ul{\on{GT}}$
are expressed via group morphisms, namely $\kappa_{21}(f)f=1$, 
$\kappa_{31}(f)\kappa_3({{\lambda-1}\over 2}) 
\kappa_{23}(f)\kappa_2({{\lambda-1}\over 2}) 
f\kappa_1({{\lambda-1}\over 2}) =1$, 
$\partial_{3}(f)\partial_{1}(f)=\partial_{0}(f)\partial_{2}(f)
\partial_{4}(f)$, where $\kappa_{ij} : F_2 \to F_2$ is 
$f(x_1,x_2) \to f(x_i,x_j)$, $\kappa_i:\ZZ\to F_2$ is $\mu\mapsto y_i^\mu$, 
and $\partial_{i} : F_2 \subset P_{3}\to P_{4}$ are the morphisms from 
\cite{Dr2}; explicitly, $\partial_{0}:x,y\mapsto x_{12},x_{23}$, 
$\partial_{4}:x,y\mapsto x_{01},x_{12}$, 
$\partial_{1}:x,y\mapsto x_{02}x_{12},x_{23}$, 
and $\partial_{2}:x,y\mapsto x_{01}x_{02},x_{13}x_{23}$. 
On the other hand, the composition is $f:= 
\theta_{(\lambda_2,f_2)}(f_1)f_2$, where $\theta_{(\lambda,f)} : 
F_2\to F_2$ is $(x,y)\mapsto (f x^\lambda f^{-1},y^\lambda)$. 
$\wh{\ul{\on{GT}}}$ is then defined as the set of $(\lambda,f)
\in (2\wh\ZZ+1)\times\wh F_2$, satisfying the profinite analogues of the 
above conditions. Any $(\lambda,f)\in\wh\ZZ\times\wh F_2$ gives rise to 
$\theta_{(\lambda,f)}\in\on{End}(\wh F_2)$, which serves to define the 
composition in $\wh{\ul{\on{GT}}}$. 

In the same way, we define $\ul{\wh{\on{GTM}}}$ as the set of all 
$(\lambda,\mu,f,g) \in (2\wh\ZZ+1)\times \wh\ZZ \times \wh F_{2}
\times \wh F_{2}$ satisfying the profinite analogues of the
defining conditions of $\ul{\on{GTM}}$. The composition formulas 
for $\ul{\on{GTM}}$ extend to $\wh{\ul{\on{GTM}}}$ and define a semigroup 
structure. $\wh{\ul{\on{GTM}}}$ acts on the pairs $(\cC,\cM)$, such that 
the image of the braid groups are finite. We have a semigroup morphism 
$\wh{\ul{\on{GTM}}}\to \wh{\ul{\on{GT}}}$, compatible with the action 
of $\wh{\ul{\on{GT}}}$ on similar categories $\cC$. 

\begin{proposition} \label{prop:str:GTM}
$\ul{\wh{\on{GTM}}} = \{(\lambda,\lambda,f,y^sf)|(\lambda,f)\in 
\ul{\wh{\on{GT}}}$ and $s\in\wh\ZZ\}$. So 
$\ul{\wh{\on{GTM}}} \simeq \ul{\wh{\on{GT}}} \ltimes \wh\ZZ$, 
where $\wh\ZZ$ is equipped with its additive structure and the 
action of $\ul{\wh{\on{GT}}}$ on it is via the semigroup 
character $\chi : \ul{\wh{\on{GT}}}\to \wh\ZZ$, $(f,\lambda)\mapsto \lambda$
(where $\wh\ZZ$ is equipped with its mutiplicative semigroup structure). 
\end{proposition}

{\em Proof.} Let $(f,\lambda)\in\wh{\ul{\on{GT}}}$ and let us find the set 
of all $(g,\mu)$ such that $(\lambda,\mu,f,g)\in\wh{\ul{\on{GTM}}}$. 
Let $s\in\wh\ZZ$ be such that $g(1,y)=y^{s}$. As in Proposition 
\ref{prop:GTM}, applying the morphism $\wh P_{4}\to \wh P_{3}$ of removal
of the first strand to (\ref{cond:g:2}), we get $g(x,y)=y^{s}f(x,y)$.  
Abelianizing the profinite version of (\ref{cond:g:1}), we find that 
$\mu=\lambda$. Conversely, one shows that for any $s\in\wh\ZZ$, 
$(\lambda,\lambda,f,y^sf)\in\wh{\ul{\on{GTM}}}$. So we have a 
bijection $\wh{\ul{\on{GTM}}}\simeq \wh{\ul{\on{GT}}}\times\wh\ZZ$. One 
then checks that the semigroup structure is the semidirect product 
$(\ul f_{1},s_{1})(\ul f_{2},s_{2}) = (\ul f_{1}\ul f_{2},s_{1}\chi(\ul f_{2})
+s_{2})$. \hfill \qed\medskip 

We now identify $(\lambda,\lambda,f,y^s f)\in \ul{\wh{\on{GTM}}}$
with $((f,\lambda),s)\in\ul{\wh{\on{GT}}}\ltimes\wh\ZZ$. If $G$ is a 
semigroup, we denote by $G^{inv}$ the group of its invertible elements. 
Recall that $\ul{\wh{\on{GT}}}^{inv}$ is denoted $\wh{\on{GT}}$ and that 
$\wh{\on{GT}} \subset \wh{\ul{\on{GT}}} \times_{\wh\ZZ}\wh\ZZ^{\times}$. 
If we set $\wh{\on{GTM}}:= \wh{\on{\ul{GTM}}}^{inv}$, then $\wh{\on{GTM}}
= \wh{\on{GT}}\ltimes\wh\ZZ = \{((f,\lambda),s) | (f,\lambda)\in\wh{\on{GT}}\}$. 

\section{The groups $\on{GTM}(N)_l$ and $\on{GTM}(N,\kk)$}\label{sec:GTMl}

If $G$ is a group and $l$ is a prime number, we denote by $G_{l}$ its pro-$l$ completion. 

\subsection{Partial profinite and pro-$l$ completions}

Let $\varphi : \Gamma\to\Gamma_{0}$ be a surjective group morphism. 
There exists a group $\Gamma(\varphi,l)$, fitting in an exact sequence 
$1\to (\Gamma_{1})_{l}\to \Gamma(\varphi,l)\to \Gamma_{0}\to 1$, 
and a group morphism $\Gamma\to \Gamma(\varphi,l)$, such that 
the diagram 
$$
\begin{matrix}
1 \to & \Gamma_1 & \to & \Gamma & \to & \Gamma_0 & \to 1 \\  
      & \downarrow & & \downarrow & & \parallel & \\ 
1 \to & (\Gamma_1)_{l} & \to & \Gamma(\varphi,l) & \to & \Gamma_0 & \to 1  
\end{matrix}
$$ 
commutes and with universal properties similar to those of partial 
prounipotent completions (see Section \ref{partial:prounip}). 

We also have a group $\wh\Gamma(\varphi)$ with the same properties, 
where $(\Gamma_{1})_{l}$ is replaced by $\wh\Gamma_{1}$; we then have a 
morphism $\wh\Gamma(\varphi)\to \Gamma(\varphi,l)$, fitting in a
commutative diagram $\begin{matrix}  & \nearrow& \wh\Gamma(\varphi)& \searrow 
& \\
\Gamma & & \downarrow & & \Gamma_{0}\\
& \searrow & \Gamma(\varphi,l)& \nearrow & \end{matrix}$

We have $\wh\Gamma(\varphi) = \lim_{\leftarrow N \lhd \Gamma_{1}\on{\ 
cofinite\ and\ } N \lhd \Gamma}\Gamma/N$, $\Gamma(\varphi,l) = 
\lim_{\leftarrow N \lhd \Gamma_{1}\ \on{co-}l
\on{\ and\ } N \lhd \Gamma}\Gamma/N$, where "$N \lhd \Gamma_{1}$ cofinite"
(resp., co-$l$) means that $\Gamma_{1}/N$ is finite (resp., a $l$-group). 

By universal properties, the construction of $\Gamma(\varphi,l)$ and $\wh\Gamma(\varphi)$
is functorial w.r.t. commuting squares $\begin{matrix}\Gamma & \to & \Gamma_{0} \\
\downarrow & & \downarrow \\ \Gamma' & \to & \Gamma'_{0}\end{matrix}$. 

\begin{lemma}
1) If $\varphi:
\Gamma \to \Gamma_{0}$ and $\pi : \Gamma_{0}\to \Gamma'_{0}$ are surjective 
group morphisms, such that $\on{Ker}\pi$ is a $l$-group, then we have an 
isomorphism $\Gamma(\varphi,l)\simeq \Gamma(\pi\circ\varphi,l)$, such that 
the diagram $\begin{matrix} \Gamma(\varphi,l) & \simeq & 
\Gamma(\pi\circ\varphi,l) \\ \downarrow & & \downarrow \\ \Gamma_{0} & 
\stackrel{\pi}{\to} & \Gamma'_{0}\end{matrix}$ commutes. 

2) Same statements replacing "$l$-group" by "finite" and $\Gamma(\psi,l)$ by 
$\wh\Gamma(\psi)$ ($\psi = \varphi$, $\pi\circ \varphi$). 
\end{lemma}

{\em Proof.} Let us prove 1) (2) is proved in the same way). 
By universal properties, we have a morphism 
$\Gamma(\varphi,l)\to \Gamma(\pi\circ\varphi,l)$ commuting with the inclusions 
of $\Gamma$ and the morphisms to $\Gamma'_{0}$. To show that this is an isomorphism, 
reinterpret this map as follows. 
 We have $\Gamma(\varphi,l) = \lim_{\leftarrow N\in \cS}\Gamma/N$
and $\Gamma(\pi\circ\varphi,l)=\lim_{\leftarrow N'\in \cS'}\Gamma/N'$, 
where $\cS = \{N|N\lhd \Gamma_{1}$ co-$l, N\lhd\Gamma\}$, 
$\cS' = \{N'|N'\lhd\Gamma'_{1}$ co-$l, N'\lhd\Gamma\}$
(recall that $\Gamma_{1}= \on{Ker}\varphi$; we set $\Gamma'_{1}:= 
\on{Ker}\pi\circ\varphi$). The morphism is then induced by the map 
$\cS'\to \cS$, $N'\mapsto N'$ (which is contained in $N'$). A morphism 
$\Gamma(\pi\circ \varphi,l)\to 
\Gamma(\varphi,l)$ is then induced by $\cS\to \cS'$, $N\mapsto 
\cap_{\gamma\in\Gamma_{1}/\Gamma'_{1}}\tilde\gamma N \tilde\gamma^{-1}$
(which is contained in $N$), 
where $\gamma\mapsto\tilde\gamma$ is a section of the projection $\Gamma_{1}
\to \Gamma_{1}/\Gamma'_{1}$.  \hfill \qed \medskip 

If $\wh \Gamma$ is the profinite completion of $\Gamma$ and $\Gamma_{0}$ is finite, 
then we have a 
commuting triangle $\begin{matrix} \wh\Gamma & \to & \Gamma(\varphi,l) \\
\searrow & & \swarrow \\ & \Gamma_{0}& \end{matrix}$. 
The morphism $\wh\Gamma\to \Gamma(\varphi,l)$ is the composition  
$\wh\Gamma\simeq \wh\Gamma(\varphi)\to \Gamma(\varphi,l)$. 


Recall that if $F_{n}$ is the free group with $n$ generators, we have an 
injection $(F_{n})_{l} \hookrightarrow F_{n}(\QQ_{l})$. This implies
that if $G$ is finitely generated, we have a functorial morphism 
$G_{l}\to G(\QQ_{l})$. More generally, if $\Gamma_{1}$ is finitely generated, then 
one has a morphism 
$\Gamma(\varphi,l)\to \Gamma(\varphi,\QQ_{l})$, 
and a commuting square  $\begin{matrix} \Gamma & \stackrel{\varphi}{\to} 
& \Gamma_{0} \\ {\scriptstyle f}\downarrow & & \downarrow \\
\Gamma' & \stackrel{\varphi'}{\to} & \Gamma'_{0}\end{matrix}$ gives rise
to a commuting square $\begin{matrix} \Gamma(\varphi,l) &\to & 
\Gamma(\varphi,\QQ_{l}) \\ {\scriptstyle f(\varphi,\varphi',l)}\downarrow 
& & \downarrow 
\\ \Gamma'(\varphi',l) &\to & \Gamma'(\varphi',\QQ_{l})
\end{matrix}$. 

\subsection{The rings $R(N)$} \label{R(N)}

Let $R$ be a ring with unit. We attach to it a ring $R(N)$, equipped with a surjective
ring morphism $R(N) \twoheadrightarrow \ZZ/N\ZZ$ and an exact sequence 
$0 \to R \to R(N) \to \ZZ/N\ZZ \to 0$ of additive groups. 

The ring $R(N)$ is defined as follows: $R(N) = (\ZZ/N\ZZ) \times R$, with
addition defined by $(a,r) + (a',r') = (a+a',r+r'+\sigma(a,a')N)$, 
where $\sigma(a,a')\in \{0,1\}$ is defined by $\tilde a + \tilde a' = 
\wt{a+a'} + \sigma(a,a')N$ (we denote by $\tilde x\in [0,N-1]$ the 
representative of $x\in \ZZ/N\ZZ$). The product is defined by 
$(a,r)(a',r') = (aa',
ra'+ar' + Nrr' + \pi(a,a')N)$, where $\pi(a,a')\in \NN$ is defined by 
$\tilde a \tilde a' = \wt{aa'} + \pi(a,a')N$. 

The assignment $R\mapsto R(N)$ is an endofunctor of the category of
rings with unit with the following properties: 

(a) if $N_{1}$ and $N_{2}$ are coprime, then $R(N_{1})(N_{2})\simeq R(N_{1}N_{2})$; 

(b) we have a ring morphism $R(N)\to (\ZZ/N\ZZ)\times R$, $\lambda\mapsto
(\bar\lambda,[\lambda])$, where $\bar \lambda = a$ and $[\lambda] = 
\tilde a + Nr$ for $\lambda=(a,r)$; this is an isomorphism if $N$ is invertible in $R$; 

(b') more generally, if $N'|N$ and $d=N/N'$, we have a ring morphism 
$[[-]]:R(N)\to R(N')$, $(a,r)\mapsto (\bar{\bar a},\tilde{\tilde a}+dr)$, where 
$a\mapsto \bar{\bar a}$ is the morphism $\ZZ/N\ZZ\to \ZZ/d\ZZ$, $\bar 1\mapsto \bar 1$; 
$\tilde{\tilde a}\in [0,N'-1]$ is the integral part of $\tilde a/d$; if $d$ is invertible in $R$, then 
$R(N)\to (\ZZ/d\ZZ)\times R(N')$, $(a,r)\mapsto (\bar a,[[(a,r)]])$ is an isomorphism 
(where $a\mapsto \bar a$ is $\ZZ/N\ZZ\to \ZZ/d\ZZ$, $\bar 1\mapsto \bar 1$); 

(c) if $R,R'$ are rings with unit, then $(R\times R')(N) = R(N)\times_{\ZZ/N\ZZ}R'(N)$. 

If $l$ is a prime number, then $(\ZZ/l^{n}\ZZ)(l)\simeq \ZZ/l^{n+1}\ZZ$; 
it follows that $\ZZ_{l}(l)\simeq \ZZ_{l}$.  Using (a) and (b), we then get 
$\ZZ_{l}(N) \simeq \ZZ_{l}\times (\ZZ/N'\ZZ)$, where $N' = N/l^{\alpha}$
and $\alpha$ is the $l$-adic valuation of $N$. (c) then implies that $\wh\ZZ(N)
\simeq\wh\ZZ$. One also checks directly that the map $(a,r)\mapsto 
\tilde a + Nr$ induces an isomorphism $\ZZ(N)\simeq\ZZ$. 

If $l$ is a prime number, then $\ZZ(\on{can},l) \simeq \ZZ_{l}(N)$, 
where $\on{can}:\ZZ\to \ZZ/N\ZZ$ is the canonical projection; on the 
other hand, $\ZZ(0,l)=\ZZ_{l}$. If $N'|N$, one checks that the morphism 
$\ZZ(N)_{l}\to \ZZ(N')_{l}$ induced by $\begin{matrix} \ZZ & \to & \ZZ/N\ZZ \\
 & \searrow & \downarrow \\  & & \ZZ/N'\ZZ\end{matrix}$ coincides with the morphism 
 from (b') for $R=\ZZ_{l}$. 

\subsection{The semigroup $\ul{\on{GTM}}(N)_{l}$}

Recall that $\ul{\on{GT}}_{l}$ is the set of pairs $(\lambda,f)
\in (2\ZZ_{l}+1)\times (F_{2})_{l}$, satisfying the pro-$l$ versions of
the defining conditions of $\ul{\on{GT}}$; the duality and hexagon relations 
take place in $(F_{2})_{l}$ and the pentagon relation in $(P_{4})_{l}$. 
We then have a semigroup morphism $\wh{\ul{\on{GT}}}\to \ul{\on{GT}}_l$; 
$\ul{\on{GT}}_l$ acts on the braided monoidal categories $\cC$, such that 
the image of the pure braid groups are $l$-groups.

We now introduce a semigroup $\ul{\on{GTM}}(N)_l$ acting on pairs 
$(\cC,\cM)$, such that 
the images of $K_{n,N}\subset B_n^1$ are $l$-groups. The setup of the 
definition of $\ul{\on{GTM}}$ has to be modified as follows: 
$\lambda\in 2\ZZ_l+1$, $f\in (P_3^0)_l = (F_2)_l$;  
$\mu\in \ZZ(\on{can},l)=\ZZ_l(N)$, and 
$g\in P_3^0(\varphi_{3,N},l)\simeq F_2(\varphi_N,l)$; 
here\footnote{The isomorphism 
$P_3^0(\varphi_{3,N},l)\simeq F_2(\varphi_N,l)$ is a consequence of 
the commutativity of 
$\begin{matrix} F_2 & \to  & P_3 \\ {\scriptstyle \varphi_N}\downarrow & & 
\downarrow {\scriptstyle \varphi_{2,N}} 
\\ \ZZ/N\ZZ& \stackrel{a\mapsto(a,\bar 0)}{\to}  & (\ZZ/N\ZZ)^2  \end{matrix}$} 
$\varphi_{N} : F_2 \to \ZZ/N\ZZ$ is defined by 
$x\mapsto\bar 1$, $y\mapsto \bar 0$. 

We now describe the conditions defining $\ul{\on{GTM}}(N)_l$. 
Recall that the octogon condition (\ref{cond:g:1}) can be written 
\begin{equation} \label{oct'}
g^{-1}\kappa_{y}({{\lambda-1}\over 2})\kappa(g)\kappa_{z}(\mu)
\kappa(g)^{-1}\kappa_{y}({{\lambda+1}\over 2})g\kappa_{x}(\mu)=1, 
\end{equation}
where the morphisms
$\kappa_x,\kappa_y,\kappa_z : \ZZ\to F_2$ are $a\mapsto x^a$, 
$a\mapsto y^a$, $a\mapsto (y^{-1}x^{-1})^a$, and 
$\kappa : F_2\to F_2$ is $g(x,y)\mapsto g(x^{-1}y^{-1},y)$. 

Now $\varphi_N \circ \kappa_y=0$, so $\kappa_y$ extends to a group 
morphism $\ZZ_l\to F_2(\varphi_N,l)$, and the diagrams 
$\begin{matrix} \ZZ & \stackrel{\kappa_{x,z}}{\to} & F_2 \\  
{\scriptstyle\on{can}}\downarrow & & \downarrow{\scriptstyle\varphi_N} \\  
\ZZ/N\ZZ & \stackrel{\pm\on{id}}{\to}& \ZZ/N\ZZ \end{matrix}$
commute, which implies that $\kappa_{x,z}$ extend to a morphism 
$\ZZ_l(N) \to F_2(\varphi_N,l)$. We also have a commuting diagram 
$\begin{matrix} F_2 & \stackrel{\kappa}{\to} & F_2 \\  
{\scriptstyle\varphi_N}\downarrow & & \downarrow{\scriptstyle\varphi_N} \\  
\ZZ/N\ZZ & \stackrel{-\on{id}}{\to}& \ZZ/N\ZZ \end{matrix}$, so 
$\kappa$ extends to an endomorphism of $F_2(\varphi_N,l)$. 
All this implies that $(\lambda,\mu,g)\mapsto$ l.h.s. of 
(\ref{cond:g:1})
extends to a map $\on{oct}_{N,l} : (2\ZZ_l+1)\times \ZZ_l(N) \times F_2(\varphi_N,l)\to
F_2(\varphi_N,l)$. 

The mixed pentagon condition (\ref{cond:g:2}) is written 
\begin{equation} \label{mixed:pentagon}
\partial_3(g)
\partial_1(g) = \partial_0(f)\partial_2(g)\partial_4(g).
\end{equation} 
For $i=1,...,4$, we have commuting diagrams 
$\begin{matrix} F_{2} & \stackrel{\partial_{i}}{\to}
& P_{4} \\ {\scriptstyle \varphi_{N}}\downarrow & & \downarrow{\scriptstyle
\varphi_{3,N}} \\ \ZZ/N\ZZ & \stackrel{\alpha_{i}}{\to} & (\ZZ/N\ZZ)^{3} 
 \end{matrix}$
where $\alpha_{i}(\bar 1) = (\bar 0,\bar 1,\bar 0)$ for $i=1$, 
$(\bar 1,\bar 1,\bar 0)$ for $i=2$, $(\bar 1,\bar 0,\bar 0)$ for $i=3,4$; 
moreover, $\varphi_{3,N}\circ\partial_0=0$. So the morphisms 
$\partial_i : F_{2} \to P_{4}$ induce, for $i=1,...,4$, group 
morphisms $F_{2}(\varphi_{N},l)\to P_{4}(\varphi_{3,N},l)$, 
and the morphism $\partial_0 : F_{2} \to P_{4}$ induces a morphism 
$(F_{2})_{l} \to P_{4}(\varphi_{3,N},l)$.  It follows that the 
map $(f,g)\mapsto$ l.h.s. of (\ref{cond:g:2}) extends to a map 
$\on{pent}_{N,l} : (F_2)_l \times F_2(\varphi_N,l)\to P_4(\varphi_{3,N},l)$. 

\begin{definition} {\it 
$\ul{\on{GTM}}(N)_l$ is the set of all $(\lambda,\mu,f,g)
\in (2\ZZ_l+1)\times \ZZ_l(N) \times (F_2)_l \times F_2(\varphi_N,l)$, 
such that $(\lambda,f)\in\ul{\on{GT}}_l$, $\on{oct}_{N,l}(\lambda,\mu,g) 
=1$ and $\on{pent}_{N,l}(f,g)=1$.}
\end{definition}

We now define the composition in $\ul{\on{GTM}}(N)_l$. 

We define $\lambda\mapsto x^{\lambda},y^{\lambda}$ as the morphisms
$\ZZ_{l}\to (F_{2})_{l}$ induced by $\ZZ\to F_{2}$, $a\mapsto x^{a},y^{a}$. 
Then for $(\lambda,f)\in\ZZ_l\times(F_2)_l$, there is a unique 
$\theta_{(\lambda,f)}\in \on{End}((F_2)_l)$ given by 
$x\mapsto f x^\lambda f^{-1}$, $y\mapsto y^\lambda$. 

Define $\mu\mapsto x^{\mu}$ as the morphism $\ZZ_{l}(N)\to F_{2}(\varphi_{N},l)$
induced by $\ZZ\to F_{2}$, $a\mapsto x^{a}$ and $\lambda\mapsto y^{\lambda}$
as the morphism $\ZZ_{l}\to F_{2}(\varphi_{N},l)$ induced by $\ZZ\to F_{2}$, 
$a\mapsto y^{a}$. For $(\lambda,\mu,g)\in 
\ZZ_l \times \ZZ_l(N) \times F_2(\varphi_N,l)$, 
there is a unique $\theta_{(\lambda,\mu,g)}\in 
\on{End}(F_2(\varphi_N,l))$ given by $x\mapsto g x^{\mu} g^{-1}$, 
$y\mapsto y^\lambda$. 

The composition law of $\ul{\on{GTM}}(N)_l$ is defined by the same formulas
as for $\ul{\on{GTM}}$, where the composition for $f$ is understood as 
$f := \theta_{(\lambda_2,f_2)}(f_1)f_2$ and the composition for $g$ as 
$g:= \theta_{(\lambda_2,\mu_2,g_2)}(g_1)g_2$. One then checks that 
$\ul{\on{GTM}}(N)_l$ is a semigroup and 
$\ul{\on{GTM}}(N)_l\to \ul{\on{GT}}_l$, $(\lambda,\mu,f,g)
\mapsto (\lambda,f)$ is a semigroup morphism. 

\begin{lemma}
If $(\lambda,\mu,f,g)\in \ul{\on{GTM}}(N)_l$, then $\lambda=[\mu]$, 
where $[-]$ is the morphism $\ZZ_l(N)\to \ZZ_l$ from \ref{R(N)}, b). 
\end{lemma}

{\em Proof.} The abelianization morphism $F_2\to\ZZ^2$ fits in a 
commutative diagram $\begin{matrix} F_2 & \to & \ZZ^2 \\
 {\scriptstyle\varphi_N}\downarrow & & \downarrow{
  \begin{matrix} {\scriptstyle (a,b)\mapsto} \\ \scriptstyle{\on{can}(a)}
  \end{matrix} }
 \\ \ZZ/N\ZZ & \stackrel{\id}{\to}& \ZZ/N\ZZ\end{matrix}$ so it induces a
morphism $\on{ab} : F_2(\varphi_N,l)\to \ZZ_l(N) \times \ZZ_l$. 
Let us apply $\on{ab}$ to (\ref{oct'}). We have $\on{ab}(\kappa_y(\nu)) = (0,\nu)$ for 
$\nu\in\ZZ_l$, and for $\nu'\in\ZZ_l(N)$, $\on{ab}(\kappa_x(\nu'))=(\nu',0)$  
and $\on{ab}(\kappa_z(\nu'))=-(\nu',[\nu'])$; this follows from the 
fact that the morphism $\ZZ_l(N)\to\ZZ_l$ induced by the diagram 
$\begin{matrix} \ZZ & \stackrel{\id}{\to}& \ZZ \\ {\scriptstyle \on{can}}
\downarrow & & \downarrow{\scriptstyle 0} \\ \ZZ/N\ZZ & \stackrel{0}{\to} & 
\ZZ/N\ZZ \end{matrix}$ is $[-]$. The image of (\ref{oct'}) by ab is then 
$(0,[\mu]-\lambda)$, so $[\mu]=\lambda$. \hfill \qed
\medskip

Let $b:F_2\to\ZZ$ be the morphism $x\mapsto 1$, $y\mapsto 0$. 
We have a commutative diagram $\begin{matrix} F_2 &\stackrel{b}{\to}& \ZZ 
\\ \scriptstyle{\varphi_N}\downarrow && \downarrow{\scriptstyle{\on{can}}}
\\ \ZZ/N\ZZ &\stackrel{\id}{\to}& \ZZ/N\ZZ \end{matrix}$ hence 
$b$ induces a morphism $b_l : F_2(\varphi_N,l)\to \ZZ_l(N)$. 
In particular, the morphism $F_2(\varphi_N,l)\stackrel{(\varphi_N)_l}{\to}
\ZZ/N\ZZ$ induced by $\varphi_N$ factors as $F_2(\varphi_N,l)
\stackrel{b_l}{\to}\ZZ_l(N) \to \ZZ/N\ZZ$. 

\begin{lemma} \label{lemma:b0}
If $(\lambda,\mu,f,g)\in \ul{\on{GTM}}(N)_{l}$, then 
$g\in \on{Ker}[b_l : F_2(\varphi_N,l) \to \ZZ_l(N)]$; in 
particular, $g\in \on{Ker}(\varphi_N)_l = (\on{Ker}\varphi_N)_l$. 
\end{lemma}

{\em Proof.} It is well-known that $f$ belongs to the kernel of the
abelianization morphism $(F_2)_l \to \ZZ_l^{\oplus 2}$ (this can be checked 
by abelianizing the pentagon identity). The mixed pentagon identity can be 
similarly abelianized: the abelianizations of the maps 
$\partial_i$ ($i=0,...,4$) fit in commuting squares 
$\begin{matrix} F_2 &\stackrel{\partial_i}{\to} & P_4 \\ \scriptstyle{\on{ab}}
\downarrow && \downarrow\scriptstyle{\on{ab}} \\ \ZZ x^{ab}\oplus \ZZ y^{ab} &
\stackrel{\partial_i^{ab}}{\to} 
& \oplus_{(i,j)|0\leq i<j\leq 3} \ZZ x_{ij}^{ab} \end{matrix}$ 
which give rise to 
squares
$$
\begin{matrix} F_2(\varphi_N,l) &\stackrel{\partial_i}{\to} &
P_4(\varphi_{3,N},l) \\ \scriptstyle{\on{ab}}
\downarrow && \downarrow\scriptstyle{\on{ab}} \\ 
\begin{matrix} \ZZ_l(N)x^{ab} \\  \oplus \ZZ_l
y^{ab} \end{matrix}&
\stackrel{(\partial_i)_l^{ab}}{\to} 
& \begin{matrix} \oplus_{i=1}^3 \ZZ_l(N)x_{0i}^{ab}
 \\ \oplus \oplus_{1\leq i<j\leq 3} \ZZ_lx_{ij}^{ab} \end{matrix}
 \end{matrix}
\on{\ for\ }i\neq 0,\on {\ and\ } 
\begin{matrix} (F_2)_l &\stackrel{(\partial_0)_l}{\to} & P_4(\varphi_{3,N},l) 
\\ \scriptstyle{\on{ab}}
\downarrow && \downarrow\scriptstyle{\on{ab}} \\ 
\begin{matrix} \ZZ_l x^{ab} \\
 \oplus \ZZ_l y^{ab} \end{matrix} &
\stackrel{(\partial_0)_l^{ab}}{\to} 
& \begin{matrix} \oplus_{i=1}^3 \ZZ_l(N)x_{0i}^{ab} \\ 
\oplus \oplus_{1\leq i<j\leq 3}
\ZZ_l x_{ij}^{ab} \end{matrix}\end{matrix}.
$$ 
Set $\on{ab}(g)= b x^{ab} + c y^{ab}$, $b\in \ZZ_l(N)$, $c\in \ZZ_l$; in 
fact, $b = b_l(g)$. 
Then $(\partial_2)_l(g) \stackrel{\on{ab}}{\mapsto} b(x_{01}^{ab}+x_{02}^{ab})
+c(x_{13}^{ab}+x_{23}^{ab})$, $(\partial_4)_l(g) \stackrel{\on{ab}}{\mapsto} 
b x_{01}^{ab}+cx_{12}^{ab}$, $(\partial_3)_l(g)\stackrel{\on{ab}}{\mapsto} b
x_{01}^{ab}+c(x_{12}+x_{13})^{ab}$, and $(\partial_1)_{l}(g)
\stackrel{\on{ab}}{\mapsto} bx_{02}^{ab}+[b]x_{12}^{ab}+c x_{23}^{ab}$. 
Since $(\partial_0)_l(f) \stackrel{\on{ab}}{\mapsto} 0$, we get  
$\on{ab}((\partial_2)_l(g)) + \on{ab}( (\partial_4)_l(g) ) 
= \on{ab}( (\partial_1)_l(g) ) + \on{ab}( (\partial_3)_l(g) )$, which 
implies $b=0$, as wanted. 
\hfill \qed \medskip

If $N'|N$, the commuting diagrams $\begin{matrix} \ZZ & 
\to & \ZZ/N\ZZ \\  & \searrow  & \downarrow \\ & & \ZZ/N'\ZZ\end{matrix}$ and 
$\begin{matrix} F_2 & \stackrel{\varphi_N}{\to} & \ZZ/N\ZZ 
\\  & \scriptstyle{\varphi_{N'}}\searrow 
& \downarrow \\ & & \ZZ/N'\ZZ\end{matrix}$ induce morphisms 
$\pi_{NN'} : \ZZ_l(N)\to\ZZ_l(N')$, 
$F_2(\varphi_N,l)\to F_2(\varphi_{N'},l)$ which take place in commuting 
squares $\begin{matrix} F_2(\varphi_N,l) & \to & F_2(\varphi_{N'},l)
\\ \scriptstyle{\partial_{1,2,3}}\downarrow & &
\downarrow\scriptstyle{\partial_{1,2,3}}\\ 
P_4(\varphi_N,l)& \to & P_4(\varphi_{N'},l)
\end{matrix}$, etc. We derive from there morphisms 
$$
\pi_{NN'} : \ul{\on{GTM}}(N)_l \to \ul{\on{GTM}}(N')_l, 
$$
such that $\pi_{NN''}= \pi_{N'N''}\circ \pi_{NN'}$. 
The categorical interpretation of $\pi_{NN'}$ is that as 
$K_{n,N}\subset K_{n,N'}$, we have $\on{Cat}_{N'} \subset \on{Cat}_N$
where $\on{Cat}_N = \{$pairs $(\cC,\cM)$ where 
the image of $K_{n,N}$ in each $\on{Aut}_{\cM}(M\otimes X^{\otimes n})$
is a $l$-group$\}$; the action of 
$g\in \ul{\on{GTM}}(N)_l$ restricts to the action of $\pi_{NN'}(g)$
on $\on{Cat}_{N'}$. 

If now $N' = N/l^\alpha$, where $\alpha$ is the $l$-adic valuation of $N$, 
$\pi_{NN'}$ is an isomorphism as $\ZZ_l(N) \simeq \ZZ_l(N')$ and  
$F_2(\varphi_N,l) \simeq F_2(\varphi_{N'},l)$. So  
$$
\ul{\on{GTM}}(N)_l \simeq \ul{\on{GTM}}(N/l^\alpha)_l.
$$
This is consistent with the categorical interpretation
of $\ul{\on{GTM}}(N)_l$, since the image of $K_{n,N}$ is an $l$-group  
iff the image of $K_{n,N'}$ is, as $K_{n,N}$ is a normal subgroup of 
$K_{n,N'}$ and the quotient in an $l$-group, so $\on{Cat}_N =
\on{Cat}_{N/l^\alpha}$. 

\subsection{The morphism $\wh{\ul{\on{GTM}}} \to \ul{\on{GTM}}(N)_{l}$}

Recall that we have a group morphism $\wh F_2 \to F_2(\varphi_N,l)$
such that $\begin{matrix} \wh F_2 & \to & F_2(\varphi_N,l) \\ 
\searrow & & \swarrow \\ & \ZZ/N\ZZ& \end{matrix}$ commutes. 
In the same way, we have an additive group morphism $\wh\ZZ\to
\ZZ_l(N)$; this is actually a ring morphism, since we have a
ring isomorphism $\wh\ZZ\simeq \wh\ZZ(N)$, which we compose with the 
morphism $\wh\ZZ(N) \to \ZZ_l(N)$ induced by the projection
$\wh\ZZ\to\ZZ_l$. These induce maps $(1+2\wh\ZZ)\times\wh\ZZ
\times\wh F_2\times\wh F_2\to (1+2\ZZ_l)\times\ZZ_{l}(N)\times (F_2)_l
\times F_2(\varphi_N,l)$; one checks that this induces a semigroup 
morphism $\wh{\ul{\on{GT}}}\ltimes \wh\ZZ \to 
\wh{\ul{\on{GTM}}} \to \ul{\on{GTM}}(N)_{l}$, such that 
$\begin{matrix} \wh{\ul{\on{GT}}}\ltimes\wh\ZZ & \to & \ul{\on{GTM}}(N)_l \\ 
\downarrow & & \downarrow \\ \wh{\ul{\on{GT}}} & \to & \ul{\on{GT}}_l
\end{matrix}$ commutes. 

We have a semigroup morphism $\ul{\on{GTM}}(N)_{l}\to \ZZ_{l}(N)=
(\ZZ/N\ZZ) \times \ZZ_{l}$ (with multiplicative structure), $(\lambda,\lambda',f,g)
\mapsto \lambda'$. The composed map $\wh\ZZ\subset \wh{\ul{\on{GT}}}
\ltimes\wh\ZZ\to \wh\ZZ(N)$ is trivial. Let $\on{GTM}(N)_{l}:= 
\ul{\on{GTM}}(N)_{l}^{inv}$. The composed map $\on{Gal}(\bar\QQ/\QQ)
\to \wh{\on{GT}} \to \on{GTM}(N)_{l} \to (\ZZ/N\ZZ)^{\times} \times 
\ZZ_{l}^{\times}$ is the cyclotomic character. This implies that the group morphism 
$\on{Gal}(\bar\QQ/\QQ) \to \on{GTM}(N)_{l}$ restricts to 
$\on{Gal}(\bar\QQ/\QQ(\mu_{Nl^{\infty}})) \to \on{GTM}_{(\bar 1,1)}(N)_{l}$, 
where $\on{GTM}_{(\bar 1,1)}(N)_{l}$ is the preimage of the unit element $(\bar 1,1)$
by $ \on{GTM}(N)_{l} \to (\ZZ/N\ZZ)^{\times} \times \ZZ_{l}^{\times}$. 

\subsection{The semigroup $\ul{\on{GTM}}(N,\kk)$}
 
 Recall that $\ul{\on{GT}}(\kk)$ is the set of pairs $(\lambda,f)\in 
 \kk \times F_{2}(\kk)$, satisfying the field versions of the defining 
 conditions on $\ul{\on{GT}}$; the hexagon and duality conditions take place
 in $F_{2}(\kk)$ and the pentagon condition in $P_{4}(\kk)$. This is a 
 semigroup, which acts on braided monoidal categories $\cC$ over $\kk[[\hbar]]$, 
 such that  the images of the pure braid groups are of the form $\on{id}+O(\hbar)$; 
 and also on QTQBA's over $\kk[[\hbar]]$, such that $R = 1^{\otimes 2} + O(\hbar)$
 and $\Phi_{A}=1^{\otimes 3}+O(\hbar)$. We have a semigroup morphism 
 $\ul{\on{GT}}(\kk)\to\kk$, $(\lambda,f)\mapsto \lambda$, and 
 $\on{GT}(\kk):= \ul{\on{GT}}(\kk)^{inv}$ is the preimage of $\kk^{\times}$. 
 
 We define $\ul{\on{GTM}}(N,\kk)$ as the set of $(\lambda,\mu,f,g)
 \in \kk \times \kk(N)\times F_{2}(\kk) \times F_{2}(\varphi_{N},\kk)$, 
 satisfying the field octogon condition (in $F_{2}(\varphi_{N},\kk)$) and 
 mixed pentagon condition (in $P_{4}(\varphi_{3,N},\kk)$). As before, 
 if $(\lambda,\mu,f,g)\in\ul{\on{GTM}}(N,\kk)$, then $\lambda
 =[\mu]$ and $g\in \on{Ker}\varphi_{N}(\kk) = (\on{Ker}\varphi_{N})(\kk)$, where 
 $\varphi_{N}(\kk):F_{2}(\kk)\to\ZZ/N\ZZ$ is induced by $\varphi_{N}$. 
 
 We have a semigroup morphism $\ul{\on{GTM}}(N,\kk)\to \kk(N)\simeq 
 (\ZZ/N\ZZ)\times \kk$, $(\lambda,\mu,f,g)\mapsto \mu$. One checks that 
 the group $\on{GTM}(N,\kk):= \ul{\on{GTM}}(N,\kk)^{inv}$ coincides with the 
 preimage of $\kk(N)^{\times} = (\ZZ/N\ZZ)^{\times}\times \kk^{\times}$. 
 
 We also have a semigroup morphism $\ul{\on{GTM}}(N,\kk)\to\ul{\on{GT}}(\kk)$. 
 The semigroup $\ul{\on{GTM}}(N,\kk)$ acts on pairs $(\cC,\cM)$ defined over 
 $\kk[[\hbar]]$, such that the images of $K_{n,N}$ are of the form $\id + O(\hbar)$, 
 and also on pairs $(A,B)$, where $A$ is a QTQBA as above and  $B$ is a QRA
 over it, such that $E^{N} = 1\otimes 1 + O(\hbar)$. 
 
Let us prove that the map $F_{2}(\varphi_{N},l)\to F_{2}(\varphi_{N},\QQ_{l})$ is 
injective. For this, it suffices to prove that $(\on{Ker}\varphi_{N})_{l}\to 
(\on{Ker}\varphi_{N})(\QQ_{l})$ is injective. This follows from the fact that 
$\on{Ker}\varphi_{N}$ is the free group with generators $X$, 
$y(\alpha)$, $\alpha= 0,...,N-1$ (where $X:= x^{N}$ and 
$y(\alpha):= x^{\alpha}yx^{-\alpha}$). 

The inclusion $(2\ZZ_{l}+1)\times\ZZ_{l}(N)\times (F_{2})_{l}\times 
F_{2}(\varphi_{N},l) \subset \QQ_{l}\times\QQ_{l}(N)\times F_{2}(\QQ_{l})
\times F_{2}(\varphi_{N},\QQ_{l})$ then gives rise to a semigroup inclusion 
$\ul{\on{GTM}}(N,l)\subset \ul{\on{GTM}}(N,\QQ_{l})$ and an inclusion 
between the corresponding groups. 

\subsection{The Lie algebra $\gtm(N,\kk)$}

Define the Magnus semigroup $\ul{\on{Mag}}(N)_{l}:= \{(\lambda,\mu,f,g)
\in \ZZ_{l}\times \ZZ_{l}(N)\times (F_{2})_{l} \times (\on{Ker}\varphi_{N})_{l}|
\lambda = [\mu]\}$, equipped with the composition law of $\ul{\on{GTM}}(N)_{l}$. 
We have a morphism $\ul{\on{Mag}}(N)_{l}\to \on{End}((F_{2})_{l})^{op} 
\oplus \on{End}(F_{2}(\varphi_{N},l))^{op}$, $(\lambda,\mu,f,g)\mapsto 
(\theta_{(\lambda,f)},\theta_{(\lambda,\mu,g)})$. 

In the same way, one define $\ul{\on{Mag}}(N,\kk):= 
\{(\lambda,\mu,f,g)\in \kk\times\kk(N)
\times F_{2}(\kk)\times (\on{Ker}\varphi_{N})(\kk)|\lambda=[\mu]\}$ 
equipped with the composition law of $\ul{\on{GTM}}(N,\kk)$. We also have 
semigroup morphisms $\ul{\on{Mag}}(N,\kk)\to \on{End}(F_{2}(\kk))^{op}
\oplus \on{End}(F_{2}(\varphi_{N},\kk))$ and
$\ul{\on{Mag}}(N,\kk)\to\kk(N)$, $(\lambda,...)
\mapsto \mu$. Let $\on{Mag}(N,\kk):= \ul{\on{Mag}}(N,\kk)^{inv}$ be the preimage of
$\kk(N)^{\times}$. Let ${\on{Mag}}_{(\bar 1,1)}(N,\kk)$, 
(resp., 
${\on{Mag}}_{(\bar 1,1)}(N,\kk)$) be the preimage of $(\bar 1,1)$
under $\ul{\on{Mag}}(N,\kk)\to \kk(N)^{\times}$, resp., of $\bar 1$ under 
$\ul{\on{Mag}}(N,\kk)\to (\ZZ/N\ZZ)^{\times}$. We have an exact sequence 
$1\to \on{Mag}_{(\bar 1,1)}(N,\kk)\to \on{Mag}_{\bar 1}(N,\kk)\to\kk^{\times}
\to 1$. 

Let $\hat\f(u_{i},\in I)$ be the topologically free Lie algebra generated by 
$u_{i},i\in I$. 
As $\on{Ker}\varphi_{N}$ is the free group generated by $X := x^{N}$ and 
$y(\alpha):=x^{\alpha}yx^{-\alpha}$, $\alpha\in [0,N-1]$, 
$\on{Lie}(\on{Ker}\varphi_{N})(\kk)
= \hat\f(\Xi,\eta(\alpha),\alpha\in [0,N-1])$; here  $\Xi:= \on{log}X$ and 
$\eta(\alpha):= 
\on{log}y(\alpha)$; similarly $\on{Lie}F_{2}(\kk) = \hat\f_{2}(\xi,\eta)$, 
where $\xi:= \on{log}x$ and $\eta:= \on{log}y$. 

\begin{lemma} 
$\on{Mag}_{(1,\bar 1)}(N,\kk)$ is a prounipotent Lie group. The Lie algebra of 
$\on{Mag}_{\bar 1}(\kk)$ is $\kk\times \on{Lie}F_{2}(\kk)\times 
\on{Lie}(\on{Ker}\varphi_{N})(\kk) \simeq \kk\times \hat\f(\xi,\eta)\times 
\hat\f(\Xi,\eta(\alpha),\alpha\in[0,N-1])$, equipped with the Lie bracket
$$
[(s_1,\varphi_1,\psi_1),(s_2,\varphi_2,\psi_2)] =
(0,\langle\varphi_1,\varphi_2\rangle, \langle \psi_1,\psi_2\rangle), 
$$
where 
$$
\langle \varphi_1, \varphi_2\rangle  = [\varphi_1,\varphi_2] + s_2 D(\varphi_1)
- s_1 D(\varphi_2) + D_{\varphi_2}(\varphi_1) - D_{\varphi_1}(\varphi_2); 
$$
where $D,D_\varphi$ are the derivations of 
$\hat\f(\xi,\eta)$ defined by $D(\xi) = \xi$, 
$D(\eta) = \eta$ and $D_\varphi(\xi) = [\varphi,\xi]$, $D_{\varphi}(\eta)=0$,  
and
$$
\langle \psi_1, \psi_2\rangle  = 
 [\psi_1,\psi_2] + s_2 \overline D(\psi_1)
- s_1 \overline D(\psi_2) + \overline D_{\psi_2}(\psi_1) - 
\overline D_{\psi_1}(\psi_2); 
$$
here $\overline D, \overline D_\psi$ are the derivations of 
$\hat\f(\Xi,\eta(0),\ldots,\eta(N-1))$ defined by 
$$
\overline D : \Xi \mapsto \Xi, \; 
\eta(\alpha) \mapsto \eta(\alpha) + {\alpha\over N} [\Xi,\eta(\alpha)] \; 
(\alpha\in [0,N-1])
$$
$$
\overline D_\psi : \Xi \mapsto [\psi,\Xi], \; 
\eta(\alpha) \mapsto [\psi(\Xi | \eta(0),\ldots,\eta(N-1)) - 
\psi(\Xi | \eta(\alpha),\ldots,\eta(\alpha+N-1)) , \eta(\alpha)] \; (\alpha\in [0,N-1])
$$
(one checks that these formulas hold for any $\alpha\in \ZZ$). 
\end{lemma}

{\em Proof.} Identify $(\lambda,\mu,f,g)$ with $(a,k,f,g)$, where 
$\mu = (a,k)\in \kk(N)$; so $\lambda=\tilde a+Nk$. 

The composition in $\ul{\on{Mag}}(N,\kk)$ is given by 
$(a_{1},k_{1},f_{1},g_{1})(a_{2},k_{2},f_{2},g_{2})
=(a,k,f,g)$, where $a = a_{1}a_{2}$, $k$ is such that
$\tilde a+Nk = (\tilde a_{1}+Nk_{1})(\tilde a_{2}+Nk_{2})$, 
$f(x,y)$ is given by (\ref{comp:f}) and 
\begin{align} \label{new:semigroup} 
& g(X | y(0),\ldots, y(N-1)) = 
g_2(X | y(0),\ldots,y(N-1))  \cdot 
\\ & \nonumber 
g_1\Big( X^{\tilde a_2 + Nk_2} |  
\Ad(g_2(X | y(0),\ldots,y(N-1))^{-1})
(y(0)^{\tilde a_{2}+Nk_{2}}) , 
\\ & \nonumber 
\Ad \big( 
X^{k_2} g_2(X | y(\tilde a_2), \ldots, 
y(\tilde a_2 + N-1))^{-1} \big) (y(\tilde a_2)^{\tilde a_2+Nk_{2}}), 
\ldots, 
\\ & \nonumber
\Ad \big( X^{(N-1)k_2} g_2(X | y((N-1)\tilde a_2), \ldots, 
y((N-1)\tilde a_2 + N-1))^{-1} \big)
((y(N-1)\tilde a_2)^{\tilde a_2+Nk_{2}})
\Big) .  
\end{align}  
Here we extend $\alpha\mapsto y(\alpha)$ to $\alpha\in\ZZ$ by $y(\alpha+N)
=Xy(\alpha)X^{-1}$. It is then clear that $\ul{\on{Mag}}_{(\bar 1,1)}(N,\kk)$
is a prounipotent Lie group. The Lie algebra structure is obtained by linearization. 
\hfill \qed\medskip 

We will denote by $\on{Exp}:\hat\f(\xi,\eta)\times 
\hat\f(\Xi,\eta(0),...,\eta(N-1)) \to F_{2}(\kk)\times (\on{Ker}\varphi_{N})(\kk)$ 
the exponential map of $\on{Mag}_{(\bar 1,1)}(N,\kk)$. 

$\on{GTM}_{(\bar 1,1)}(N,\kk)$ has been defined as $\on{GTM}(N,\kk)\cap 
\on{Mag}_{(\bar 1,1)}(N,\kk)$; we similarly set 
$\on{GTM}_{\bar 1}(N,\kk):= \on{GTM}(N,\kk)\cap \on{Mag}_{\bar 1}(N,\kk)$. 
We then have an exact sequence $1\to \on{GTM}_{(\bar 1,1)}(N,\kk)\to 
\on{GTM}_{\bar 1}(N,\kk)\to \kk^{\times}$. We denote by 
$$
0\to \gtm_{(\bar 1,1)}(N,\kk)\to \gtm_{\bar 1}(N,\kk)\to \kk
$$
the corresponding exact sequence of Lie algebras. 

\begin{lemma} \label{prop:gtm}
The Lie subalgebra $\gtm_{\bar 1}(N,\kk)\subset \kk\times\hat\f(\xi,\eta)
\times \hat\f(\Xi,\eta(0),...,\eta(N-1))$ is the set of all $(s,\varphi,\psi)$ 
such that 
$$
\varphi(\eta,\xi) = - \varphi(\xi,\eta), \; 
\varphi(\xi,\eta) + \varphi(\eta,\zeta) + \varphi(\zeta,\xi) + {s\over 2}
(\xi+\eta + \zeta) = 0, 
$$ 
where $e^\xi e^\eta e^\zeta=1$, 
$$
\varphi(\xi^{12},\xi^{23} * \xi^{24}) + \varphi(\xi^{13} * \xi^{23},\xi^{34})
= \varphi(\xi^{23},\xi^{34}) + \varphi(\xi^{12}*\xi^{13},\xi^{24}*\xi^{34})
+ \varphi(\xi^{12},\xi^{23})
$$ 
(relation in $\on{Lie}(P_4(\kk))$, where $\xi^{ij} = \on{log}x_{ij}$),  
\begin{align} \label{octogon:psi}
& \psi(\Xi | \eta(0),\ldots,\eta(N-1)) - \psi(\Xi | \eta(1),\ldots,\eta(N))
\\ & \nonumber
+ \psi(\Theta | \eta'(1),\ldots,\eta'(2-N)) 
- \psi(\Theta | \eta'(0),\ldots,\eta'(1-N)) + {s\over 2}(\eta(0) + \eta(1))
+ {s\over N}(\Xi + \Theta) = 0, 
\end{align}
where $\eta(\alpha+N) = \Ad(e^{\Xi})(\eta(\alpha))$,
$e^\Theta = e^{-\eta(0)} \cdots e^{-\eta(1-N)} e^{-\Xi}$,  
$\eta'(1) = \eta(1)$, $\eta'(0) = \eta(0)$, 
$\eta'(-1) = \Ad(e^{-\eta(0)})(\eta(-1))$, 
$\eta'(-2) = \Ad(e^{-\eta(0)}e^{-\eta(-1)})(\eta(-2))$, etc., and
in general $\eta'(\alpha-N) = \Ad(e^\Theta)(\eta'(\alpha))$, and 
\begin{align*}
& \psi(\Xi^{01} | \xi(0)^{12} * \xi(0)^{13},\ldots,\xi(N-1)^{12} *
\xi(N-1)^{13}) 
\\ & + \psi(\Xi^{02} * \xi(0)^{12} * \cdots * \xi(N-1)^{12} | 
\tilde\xi(0),\ldots,\tilde\xi(N-1))
\\ & 
= \varphi(\xi(0)^{12},\xi(0)^{23}) 
+ \psi(\Xi^{01} |  \xi(0)^{12},\ldots,\xi(N-1)^{12})
\\ & + \psi(\Xi^{01} * \Xi^{02} * \xi(0)^{12} * \cdots * \xi(N-1)^{12} *
(-N\xi(0)^{12})| \tilde{\tilde\xi}(0),\ldots, \tilde{\tilde\xi}(N-1)) , 
\end{align*} 
(equality in $\on{Lie}(K_{3,N}(\kk))$)
where $\tilde\xi(\alpha):= \on{log} \on{Ad}[x_{12}(1)...x_{12}(\alpha)
u_{12}(\alpha,-\alpha)][x_{23}(\alpha)]$ and $\tilde{\tilde\xi}(\alpha):= \on{log}
\on{Ad}[x_{12}(0)^{-\alpha}x_{12}(1)...x_{12}(\alpha)][x_{13}(\alpha)]
\on{Ad}[x_{12}(0)^{-\alpha}x_{12}(-\alpha)...x_{12}(-1)][x_{23}(\alpha)]$; 
recall that $\on{log} : X_{0i},x_{i}(\alpha)\mapsto 
\Xi^{0i},\xi(\alpha)^{ij}$ and the $u_{ij}(\alpha,\beta)$ have been defined in 
Proposition \ref{beginning}. 

The Lie algebra $\gtm_{(\bar 1,1)}(N,\kk)$ is prounipotent and 
the map $\gtm_{\bar 1}(N,\kk)\to\kk$ is $(s,\varphi,\psi)\mapsto s$. 
\end{lemma}

{\em Proof.} The quadruple $(\lambda,\mu,f,g)\in \kk\times \kk(N)\times F_{2}(\kk)
\times (\on{Ker}\varphi_{N})(\kk)$ lies in $\ul{\on{GTM}}(N,\kk)$ iff it satisfies
$\lambda = [\mu]$, i.e., $\lambda = \tilde a+Nk$ if $\mu = (a,k)$, 
$(\lambda,f)\in\ul{\on{GT}}(\kk)$, 
\begin{align} \label{new:octogon}
& 
g \big( X | y(1), \ldots, y(N) \big)^{-1} y(1)^{(\lambda-1)/2} 
\\ & \nonumber
g \big( 
y(0)^{-1} \cdots y(1-N)^{-1} X^{-1}  | 
\\ &  \nonumber
y(1), y(0), \Ad(y(0)^{-1})(y(-1)), \ldots, 
\Ad(y(0)^{-1}\cdots y(3-N)^{-1})(y(2-N)) \big)
\\ & \nonumber
(y(0)^{-1} \cdots y(1-N)^{-1}X^{-1} )^k 
y(0)^{-1} y(-1)^{-1} \cdots y(2-\tilde a)^{-1} 
\\ & \nonumber
g \big( y(1-\tilde a)^{-1} \cdots y(2-\tilde a-N)^{-1} X^{-1} |  
\\ & \nonumber 
y(1-\tilde a),\Ad(y(1-\tilde a)^{-1})(y(-\tilde a)),\ldots, 
\Ad(y(1-\tilde a)^{-1}\cdots y(3-\tilde a-N)^{-1})(y(2-\tilde a-N)) \big)^{-1}
\\ & \nonumber
y(1-\tilde a)^{(\lambda-1)/2} g(X | y(1-\tilde a),\ldots,y(N-\tilde a))
X^k =1,  
\end{align}
and equation (\ref{mixed:pentagon}), where $\partial_{0}:F_{2}(\kk)\to 
K_{3,N}(\kk)$ is induced by $x\mapsto x_{12}(0)$, 
$y\mapsto x_{23}(0)$ and $\partial_{i} : (\on{Ker}\varphi_{N})(\kk)
\to K_{3,N}(\kk)$ are given by
$$
\partial_{1} : X \mapsto X_{02}x_{12}(0)...x_{12}(N-1), \quad y(\alpha)\mapsto 
\on{Ad}[x_{12}(1)...x_{12}(\alpha)u_{12}(\alpha,-\alpha)][x_{23}(\alpha)], 
$$
$$
\partial_{2} : X \mapsto X_{01}X_{02}x_{12}(0)...x_{12}(N-1)x_{12}(0)^{-N}, 
$$
$$ y(\alpha)\mapsto  
\Ad[x_{12}(0)^{-\alpha} x_{12}(1) \cdots x_{12}(\alpha)][ x_{13}(\alpha)] 
\Ad[ x_{12}(0)^{-\alpha} x_{12}(-\alpha) \cdots x_{12}(-1)][ x_{23}(\alpha)],
 $$
$$
\partial_{3} : X\mapsto X_{01}, \quad y(\alpha)\mapsto x_{12}(\alpha)x_{13}(\alpha),\quad 
\partial_{4} : X\mapsto X_{01}, \quad y(\alpha)\mapsto x_{12}(\alpha).  
$$
The formulas are then obtained by setting $a=\bar 1$, $\lambda = 1+\eps s$, 
$f(x,y)=e^{\eps\varphi(\xi,\eta)}$, $g(X|y(0),...,y(N-1)) = 
e^{\eps\psi(\Xi|\eta(0),...,\eta(N-1))}$ and linearizing in $\eps$. 
\hfill \qed \medskip 

\begin{remark} As in \cite{Dr2}, the surjectivity of $\gtm_{\bar 1}(N,\kk)\to \kk$
for $\kk = \QQ_l$, hence $\kk= \QQ$, follows from that of the cyclotomic 
character $\on{Gal}(\bar\QQ/\QQ(\mu_N))\to \ZZ_l^\times$, and another proof 
of this surjectivity can be given using cyclotomic associators (see next 
section). \end{remark}

\section{Torsor structure of $\Pseudo(N,\kk)$ and the group $\on{GRTM}(N,\kk)$}
\label{sect:7}

\subsection{Action of $\ul{\on{GTM}}(N,\kk)$ on $\ul{\Pseudo}(N,\kk)$}

The semigroup $\ul{\on{GTM}}(N,\kk)$ acts on $\ul{\Pseudo}(N,\kk)$ 
as follows: 
\begin{equation} \label{action:gtm:ps}
(\lambda,\mu,f,g) * (a',\lambda',\Phi',\Psi') = 
(\bar\mu a',[\mu]\lambda',\Phi'',\Psi''),  
\end{equation}
where 
$$
\Phi''(t^{12},t^{23}) := \Phi'(t^{12},t^{23}) f(e^{\lambda' t^{12}}, 
\Ad(\Phi'(t^{12},t^{23})^{-1})(e^{\lambda't^{23}})), 
$$
\begin{align*}
& \Psi''(t_0^{12} | t(0)^{23},\ldots,t(N-1)^{23}) := 
\Psi'(t_0^{12} | t(0)^{23},\ldots,t(N-1)^{23})
\\ & 
g \Big( e^{\lambda' t_0^{12}} | 
\on{Ad}\big( \Psi'(t_0^{12} | t(0)^{23},\ldots,t(N-1)^{23})^{-1} \big) 
(e^{\lambda' t(0)^{23}}), 
\\ & \Ad\big(e^{(\lambda'/N)t_0^{12}} \Psi'(t_0^{12} | 
t(a')^{23}, \ldots, t(a'+N-1)^{23})^{-1}\big)
(e^{\lambda' t(a')^{23}}), \ldots,
\\ & 
\Ad\big(e^{(N-1)(\lambda'/N)t_0^{12}} \Psi'(t_0^{12} | 
t((N-1)a')^{23}, \ldots, t((N-1)a'+N-1)^{23})^{-1}\big)
(e^{\lambda' t((N-1)a')^{23}}) \Big) 
\end{align*}
(recall that $\lambda=[\mu]$, so if $\mu=(a,k)$, then $\lambda = 
\tilde a+Nk$; also $\bar\mu = a$). 
This action is compatible with the action of $\ul{\on{GT}}(\kk)$ on
$\ul{\Assoc}(\kk)$ and the natural morphisms $\ul{\on{GTM}}(N,\kk)
\to \ul{\on{GT}}(\kk)$, $\ul{\Pseudo}(N,\kk) \to \ul{\Assoc}(\kk)$. 

This action restricts to an action of $\on{GTM}(N,\kk)$ on $\Pseudo(N,\kk)$, 
compatible with the action of $\on{GT}(\kk)$ on $\Assoc(\kk)$. 

\begin{theorem} \label{thm:torsor}
If $\Pseudo(N,\kk)$ is nonempty, then the 
action of $\on{GTM}(N,\kk)$ on $\Pseudo(N,\kk)$ is free and transitive, in
other words, $\Pseudo(N,\kk)$ is a torsor under $\on{GTM}(N,\kk)$. 
\end{theorem}

{\em Proof.} The proof is parallel to \cite{Dr2}, Proposition 5.1. 
One shows that formula (\ref{action:gtm:ps}) defines a free and 
transitive action of 
$\on{Mag}(N,\kk)$ on $\{(a,\lambda,f,g)\in (\ZZ/N\ZZ)^{\times}
\times \kk^{\times} \times F_{2}(\kk) \times (\on{Ker}\varphi_{N})(\kk)\}$. 
Hence given $(a'_i,\lambda'_i,\Phi_i,\Psi_i)$, $i=0,1$ in 
$\Pseudo(N,\kk)$, there is a unique $(\lambda,\mu,f,g)\in 
\on{Mag}(N,\kk)$, such that $(\lambda,\mu,f,g) * 
(a'_0,\lambda'_0,\Phi_0,\Psi_0) = (a'_1,\lambda'_1,\Phi_1,\Psi_1)$. 
It remains to show that $(a,\lambda,f,g)\in \on{GTM}(N,\kk)$. 
\cite{Dr2}, Proposition 5.1 implies that $(\lambda,f)$ satisfies 
 (\ref{cond:f:1}), (\ref{cond:f:2}). 
Let $\rho_{(a',\lambda',\Phi,\Psi)}^{n}:B_{n}^{1}(\varphi_{n,N},\kk)\to 
\on{exp}(\hat\t_{n,N}(\kk))\rtimes (\ZZ/N\ZZ)^{n}\rtimes \SG_{n}$
be the morphism induced by $(a',\lambda',\Phi,\Psi)\in 
\Pseudo(N,\kk)$. Replacing in the identities 
(\ref{cond:Psi:1}), (\ref{cond:Psi:2}) satisfied by 
$(\Phi_{1},\Psi_{1})$, $\Phi_{1}$ and $\Psi_{1}$ by their expressions in terms 
of $(f,g)$ and $(\Phi_{0},\Psi_{0})$, and using the fact that $(\Phi_{0},\Psi_{0})$
satisfies the same identities, we find that the images of relations 
(\ref{new:octogon}), (\ref{mixed:pentagon})
by $\rho^{n}_{(a'_{0},\lambda'_{0},\Phi_{0},\Psi_{0})}$ ($n=2,3$) 
are satisfied; since $\rho^{n}_{(a'_{0},\lambda'_{0},\Phi_{0},\Psi_{0})}$
is bijective, these relations are satisfied in $B^{1}_{n}(\varphi_{n,N},\kk)$ ($n=2,3$). 
\hfill \qed \medskip 

\subsection{Pseudotwists and splittings of exact sequences}

The group $(\ZZ/N\ZZ)^\times \times \kk^\times$ acts by automorphisms of 
$\t_{n,N}$ as follows: $(c,\gamma) \cdot t_0^{1i} = 
\gamma t_0^{1i}$, $(c,\gamma) \cdot t(a)^{ij} = \gamma t(c a)^{ij}$, and it acts
by automorphisms of $\t_n$ by $(c,\gamma) \cdot t^{ij} = \gamma t^{ij}$. 
This induces an action of $(\ZZ/N\ZZ)^\times \times \kk^\times$
on $\ul{\Pseudo}(N,\kk)$, such that the map $\ul{\Pseudo}(N,\kk)
\to \ZZ/N\ZZ \times \kk$ is equivariant ($\ZZ/N\ZZ \times \kk$
being equipped with the action $(c,\gamma) \cdot (a,\lambda) = 
(ca,\gamma\lambda)$). In particular, we get isomorphisms
$\Pseudo_{(a,\lambda)}(N,\kk) \simeq \Pseudo_{(\bar 1,1)}(N,\kk)$ for any 
$(a,\lambda)\in (\ZZ/N\ZZ)^\times \times \kk^\times$. 

\begin{lemma}
If $\Pseudo_{(\bar 1,1)}(N,\kk)\neq \emptyset$, then the group morphism 
$\on{GTM}(N,\kk) \to \kk(N)^\times \simeq (\ZZ/N\ZZ)^\times 
\times \kk^\times$ is onto. 
\end{lemma}

{\em Proof.}
Let $(a,\lambda)\in (\ZZ/N\ZZ)^\times \times \kk^\times$ and 
$\Psi\in \Pseudo_{(\bar 1,1)}(N,\kk)$, then $(a,\lambda) \cdot \Psi
\in \Pseudo_{(a,\lambda)}(N,\kk)$; according to Theorem 
\ref{thm:torsor}, there is a unique $g\in \on{GTM}(N,\kk)$
such that $g * \Psi = (a,\lambda) \cdot \Phi$. Then the image of 
$g$ under $\on{GTM}(N,\kk) \to \kk(N)^\times$ is $(a,\lambda)$. 
\hfill \qed \medskip 
 
It follows that if $\Pseudo_{(\bar 1,1)}(N,\kk)\neq \emptyset$, 
the morphisms $\on{GTM}_{\bar 1}(N,\kk) \to \kk^\times$ and 
$\gtm_{\bar 1}(N,\kk)\to \kk$ are also onto. 

The actions of $(\ZZ/N\ZZ)^\times \times \kk^\times$ and $\on{GTM}(N,\kk)$ 
on $\Pseudo(N,\kk)$ commute with each other. We have therefore an action of 
$\on{GTM}(N,\kk)$ on $\Pseudo_{(\bar 1,1)}(N,\kk)$, defined by 
$(\lambda,\mu,f,\overline g) \circ \Psi = (\bar\mu,[\mu])^{-1} \cdot 
(\lambda,\mu,f,\overline g)\Psi$. 

To each $\Psi\in \Pseudo_{(\bar 1,1)}(N,\kk)$, we associate the Lie algebra of the
stabilizer of $\Psi$. This defines a section of the exact sequence
\begin{equation} \label{exact:seq}
0\to \gtm_{(\bar 1,1)}(N,\kk) \to \gtm_{\bar 1}(N,\kk) \to \kk \to 0.
\end{equation} 

\begin{lemma} If $\Pseudo_{(\bar 1,1)}(N,\kk) \neq \emptyset$, then 
the map $\Pseudo_{(\bar 1,1)}(N,\kk) \to \{$splittings of
(\ref{exact:seq})$\}$ is a bijection. 
\end{lemma} 

{\em Proof.}  $\on{GTM}_{(\bar 1,1)}(N,\kk)$ acts on the target of this map by
conjugation, and the map is then equivariant. So it suffices to show that 
$\on{Spl}(\kk) := \{$splittings of $(\ref{exact:seq})\}$ is a torsor under 
$\on{GTM}_{(\bar 1,1)}(N,\kk)$. $\on{Spl}(\kk)$ 
identifies with $\{$elements of $\gtm_{\bar 1}(N,\kk)$ of the form 
$(1,\varphi,\psi)\}$. If we fix an origin $(1,\varphi_0,\psi_0)$, 
then $\on{Spl}(\kk)$ identifies with $\gtm_{(\bar 1,1)}(N,\kk)$, and the vector field 
at $(0,\varphi',\psi') \in \gtm_{(\bar 1,1)}(N,\kk) \simeq \on{Spl}(N,\kk)$ 
induced by the action of 
$(0,\varphi,\psi) \in \gtm_{(\bar 1,1)}(N,\kk)$ is $(0,[\varphi_0,\varphi']
+ D(\varphi_0) + D_{\varphi_0}(\varphi') - D_{\varphi'}(\varphi_0), 
[\psi_0,\psi'] + \overline D(\psi_0) + \overline D_{\psi_0}(\psi')
- \overline D_{\psi'}(\psi_0) )$; this defines a linear map 
$\vartheta_{(0,\varphi',\psi')} : \gtm_{(\bar 1,1)}(N,\kk)
\to \gtm_{(\bar 1,1)}(N,\kk)$. Now 
$\gtm_{(\bar 1,1)}(N,\kk)$ has a decreasing filtration $\g^{\geq k}
= \{(0,\varphi,\psi) \in \gtm_{(\bar 1,1)}(N,\kk)  | \psi$ has degree $\geq k\}$
(this condition implies that $\varphi$ has also degree $\geq k$), 
$\vartheta_{(0,\varphi',\psi')}$ is compatible with this filtration and
the associated graded map is multiplication by $k$ in degree $k$, so that 
$\vartheta_{(0,\varphi',\psi')}$ is bijective. This suffices to show that
the action of $\on{GTM}_{(\bar 1,1)}(N,\kk)$ is free, and using successive
approximations, one can also show that it is transitive. 
\hfill \qed \medskip 

\begin{proposition} \label{prop:isr}
If $\gtm_{\bar 1}(N,\kk) \to \kk$ is surjective, then 
$\Pseudo_{(\bar 1,1)}(N,\kk)\neq \emptyset$. 
\end{proposition} 

The proof follows \cite{Dr2}, Proposition 5.2. 

{\em Proof.} Let us assume that this map is surjective, i.e., there exists an element 
$(1,\varphi,\psi)$ in $\gtm_{\bar 1}(N,\kk)$. There exists a unique pair 
$(\Phi,\Psi)$, where $\Phi\in F_2(\kk)$ and $\Psi\in(\on{Ker}\varphi_{N})(\kk)$, 
such that 
$$
\Phi(U,V)^{-1} {{\on{d}}\over {\on{d}t}}_{|t=1} \Phi(tU,tV) = 
\varphi(NU, \Ad(\Phi(U,V)^{-1})(NV)) ,  
$$ 
\begin{align*}
& \Psi(A | b(0),\ldots,b(N-1))^{-1}
{{\on{d}}\over {\on{d}t}}_{|t=1} 
\Psi(tA | tb(0),\ldots,tb(N-1))
\\ & 
= \psi\Big(NA | \Ad( \Psi(A | b(0),\ldots,b(N-1)) )^{-1}(Nb(0)),
\Ad( e^{A}\Psi(A | b(1),\ldots,b(N)) )^{-1}(Nb(1)),
\\ & \ldots,
\Ad( e^{(N-1)A}\Psi(A | b(N-1),\ldots,b(2N-2)))^{-1}(Nb(N-1)) \Big). 
\end{align*} 

Proposition 5.2 in \cite{Dr2} implies that $\Phi\in \Assoc_N(\kk)$. 
Repeating the final argument of {\it loc. cit.,} we find that 
$(\Phi,\Psi)$ satisfies the mixed pentagon relation (\ref{cond:Psi:2}). 
Let us now show that $\Psi$ satisfies the octogon relation 
(\ref{cond:Psi:1}), where $(\alpha,\lambda) = (\bar 1,N)$. Let us set 
\begin{align*}
& Q(A | b(0),\ldots,b(N-1)) := 
\Psi(C | b(0),\ldots,b(1-N)) e^C  
\Psi(C | b(-1),\ldots,b(-N))^{-1}
e^{N b(-1)/2} 
\\ & 
\Psi(A | b(-1),\ldots,b(N-2)) e^A 
\Psi(A | b(0),\ldots,b(N-1))^{-1}  
e^{N b(0)/2} .  
\end{align*} 
We have $\Psi(A| b(0),\ldots,b(N-1)) = 1 + \psi(A|b(0),\ldots,b(N-1))$
modulo degree $\geq 2$. Since $\psi$ satisfies (\ref{octogon:psi}), we get 
$Q = 1$ modulo degree $\geq 2$. 
Assume that we have proved that $Q = 1$ modulo degree $\geq n$, 
where $n\geq 2$. 

We have 
\begin{align} \label{main:loc:id}
& Q(A | b(0),\ldots,b(N-1))^{-1} {{\on{d}}\over{\on{d}t}}_{|t=1} 
Q(tA | tb(0),\ldots,tb(N-1)) 
\\ & \nonumber 
= {{N b(0)}\over 2}
+ A' - \psi(NA' | Nb'(0),\ldots,Nb'(N-1)) 
+ \psi(NA'| Nb'(-1),\ldots,Nb'(N-2))
\\ & \nonumber
+ {{Nb'(-1)}\over 2}
+ C' + \psi(NC'| Nb''(-1),\ldots,Nb''(-N)) 
- \psi(NC'| Nb''(0),\ldots,Nb''(1-N)), 
\end{align}
where 
$$
A' = \Ad\big(e^{-Nb(0)/2}\Psi(A | b(0),\ldots,b(N-1))\big)(A),
$$ 
$$
b'(a) =  \Ad \big( e^{-Nb(0)/2} \Psi(A| b(0),\ldots, b(N-1)) e^{aA}
\Psi(A| b(a),\ldots, b(a+N-1))^{-1} \big) (b(a)), 
$$
\begin{align*}
C' = \Ad\big( & e^{-Nb(0)/2} \Psi(A | b(0),\ldots,b(N-1)) e^{-A}
\\ & \Psi(A | b(-1),\ldots,b(N-2))^{-1} e^{-Nb(-1)/2} 
\Psi(C|b(-1),\ldots,b(-N))
\big)(C), 
\end{align*}
\begin{align*}
b''(a) = \Ad\big( & e^{-Nb(0)/2} \Psi(A | b(0),\ldots,b(N-1)) e^{-A}
\Psi(A | b(-1),\ldots,b(N-2))^{-1} 
\\ & e^{-Nb(-1)/2} \Psi(C|b(-1),\ldots,b(-N))
e^{-(a+1)C} \Psi(C | b(a),\ldots,b(a-N+1))^{-1}
\big) (b(a)). 
\end{align*}

Define $q(A | b(0),\ldots, b(N-1))$ as $(Q - 1$ modulo degree 
$\geq n+1)$. So the l.h.s. of (\ref{main:loc:id}) is 
$1 + n q(A | b(0),\ldots,b(N-1)) + ($degree $\geq n+1)$.  

Let us compute the r.h.s. of (\ref{main:loc:id}). We have 
$$
e^{NC'} e^{Nb'(0)} \cdots e^{Nb'(N-1)} e^{NA'}
= 1 + \sum_{a=0}^{N-1} q(A | b(a),\ldots,b(a+N-1)) 
+ (\on{degree}\geq n+1),  
$$
and
\begin{align*}
e^{Nb''(-a-1)} & = \Ad\big(Q(A | b(0),\ldots,b(N-1))^{-1}\big)
\big( e^{-Nb'(-1)} \cdots e^{-Nb'(-a)} 
e^{Nb'(-a-1)} \cdots e^{Nb'(-1)} \big) 
\\ & + \on{degree}\geq n+1. 
\end{align*}

Define $\overline b(-a)$ by 
$$
e^{N \overline b(-a-1)} = e^{-Nb'(-1)} \cdots e^{-Nb'(-a)}
e^{Nb'(-a-1)} \cdots e^{Nb'(-1)}, 
$$
$$
C'' := \Ad(\Psi(C | b(0),\ldots,b(1-N)))(C), 
$$
$$
A'' := A' - {1\over N} \sum_{a=0}^{N-1} q(A | b(a),\ldots,b(a+N-1)). 
$$
We have $C' = \Ad(Q(A | b(0),\ldots,b(a+N-1))^{-1})(C'')$, therefore 
$C' = C'' +$ degree $\geq n+1$. Therefore $e^{NC''} e^{Nb'(0)} \cdots 
e^{Nb'(N-1)} e^{NA''} = 1 +$ (degree $\geq n+1$). So (\ref{octogon:psi})
implies 
\begin{align*}
& {{Nb(0)}\over 2} + A'' - \psi(NA'' | Nb'(0),\ldots,Nb'(N-1))
+ \psi(NA'' | Nb'(-1),\ldots,Nb'(N-2))
\\ & 
+ {{Nb'(-1)}\over 2} + C'' + \psi(NC'' | N \overline b(-1),\ldots, 
N \overline b(-N)) - \psi(NC'' | B \overline b(0),\ldots, \overline b(1-N)) 
= 0.
\end{align*}
So the r.h.s. of (\ref{main:loc:id}) is equal to 
$$
A' - A'' + (\on{degree}\geq n+1) = 
{1\over N} \sum_{a=0}^{N-1} q(A | b(a),\ldots,b(a+N-1)) + 
(\on{degree}\geq n+1). 
$$
Therefore (\ref{main:loc:id}) implies 
\begin{equation}
n q(A | b(0),\ldots,b(N-1)) = {1\over N} 
\sum_{a=0}^{N-1} q(A | b(a),\ldots,b(a+N-1)). 
\end{equation}
This equation implies that $q(A | b(0),\ldots,b(N-1))$ is invariant under the 
automorphism $A\mapsto A$, $b(a)\mapsto b(a+1)$. So its r.h.s. is equal to 
$q(A | b(0),\ldots,b(N-1))$. Therefore the equation is rewritten as 
$(n-1)q(A | b(0),\ldots,b(N-1)) = 0$, which implies (since $n\geq 2$)
$q(A | b(0),\ldots,b(N-1)) = 0$, so $Q = 1$ modulo degree $\geq n+1$. 
So $Q = 1$. So $(\Phi,\Psi)\in \Pseudo_{(\bar 1,N)}(N,\kk)$; hence 
$\Pseudo_{(\bar 1,1)}(N,\kk)\neq\emptyset$. \hfill \qed \medskip 

\begin{corollary} 
$\Pseudo_{(\bar 1,1)}(N,\QQ) \neq \emptyset$. 
\end{corollary}

{\em Proof.} As in \cite{Dr2}, $\Pseudo_{(\bar 1,1)}(N,\CC)\neq \emptyset$
implies that $\gtm_{\bar 1}(N,\kk) \to \kk$ is surjective for $\kk = \CC$. 
This implies the same statement for $\kk = \QQ$. Then Proposition 
\ref{prop:isr} implies that $\Pseudo_{(\bar 1,1)}(N,\QQ) \neq \emptyset$. 
\hfill \qed \medskip

\subsection{$\on{GRTM}(N,\kk)$ and $\grtm(N,\kk)$}

Define $\on{GRTM}_{(\bar 1,1)}(N,\kk)$ as the set of pairs 
$(h,k)$ with $h\in \exp(\hat\t_3^0)$ and $k\in \exp(\hat\t_{3,N}^0)$, such that 
(we recall that $U = t^{12}$, $V = t^{23}$, and $A = t_0^{01}$, 
$b(a) = t(a)^{12}$ for $a\in \ZZ/N\ZZ$)
$$
h(U,V) = h(V,U)^{-1}, \quad 
h(U,V)h(V,W) h(W,U) =1, 
$$
\begin{equation} \label{add:cond:GRT}
U + \Ad(h(U,V)^{-1})(V) + \Ad(h(U,W)^{-1})(W) = 0 \; \on{for} \; 
U+V+W=0, 
\end{equation}
$$
h^{1,2,34} h^{12,3,4} = h^{2,3,4} h^{1,23,4} h^{1,2,3}
$$
(equality in $\exp(\hat\t_{4})$), and 
\begin{align} \label{octogon:GRTM}
& k(A|b(1),\ldots,b(N))^{-1} k(C| b(1),b(0),\ldots,b(2-N))
\\ & \nonumber 
k(C | b(0),b(-1),\ldots,b(1-N))^{-1} k(A | b(0),\ldots,b(N-1)) = 1, 
\end{align}
\begin{align} \label{add:cond:GRTM}
& A + \sum_{a=1}^N \Ad(k(A|b(a),\ldots,b(a+N-1))^{-1})(b(a))
\\ & \nonumber 
+ \Ad\big( k(A|b(0),\ldots,b(N-1))^{-1}
k(C|b(0),\ldots,b(1-N)) \big)(C) = 0, 
\end{align}
where $A + \sum_{i=0}^n b(i) + C =0$,
$$ 
k^{0,1,23} k^{01,2,3} = h^{1,2,3} k^{0,12,3} k^{0,1,2} 
$$ 
(equality in $\exp(\wh \t_{3,N})$). 

$\on{GRTM}_{(\bar 1,1)}(N,\kk)$ is a group when equipped with the product 
$(h_1,k_1) * (h_2,k_2) = (h,k)$, where 
$$
h(U,V) = h_2(U,V) h_1(U,\Ad h_2(U,V)^{-1}(V)), 
$$
\begin{align*}
& k(A | b(0),\ldots,b(N-1)) = k_2(A | b(0),\ldots,b(N-1)) \cdot 
\\ & k_1\Big( A| \Ad(k_2(A| b(0),\ldots,b(N-1))^{-1})(b(0)),\ldots,
\Ad(k_2(A| b(N-1),\ldots,b(2N-2))^{-1})(b(N-1)) \Big) ;   
\end{align*}
so we have a group inclusion $\on{GRTM}_{(\bar 1,1)}(N,\kk) \subset 
\on{Mag}_{(\bar 1,1)}(N,\kk)$. 
The action of $(\ZZ/N\ZZ)^\times \times \kk^\times$ by automorphisms 
of $\t^0_{3,N}$ and $\t^0_3$ induces its action by automorphisms of 
$\on{GRTM}_{(\bar 1,1)}(N,\kk)$. We denote by $\on{GRTM}(N,\kk)$ 
the corresponding semidirect product. 

The Lie algebra $\grtm_{(\bar 1,1)}(N,\kk)$ of $\on{GRTM}_{(\bar 1,1)}(N,\kk)$ 
consists of all pairs $(\varphi,\psi)\in \hat\t_3^0\times \hat\t_{3,N}^0$, such that: 
\begin{equation} \label{cond:1:grt}
\varphi(V,U) = -\varphi(U,V), \; \varphi(U,V) + \varphi(V,W) + \varphi(W,U) = 0
\end{equation}
if $U+V+W=0$, 
\begin{equation} \label{add:cond:grt}
[V,\varphi(U,V)] + [W,\varphi(U,W)] = 0 \; \on{if}\; U+V+W=0, 
\end{equation}
\begin{equation} \label{cond:2:grt}
\varphi^{12,3,4} + \varphi^{1,2,34} = \varphi^{1,2,3} + \varphi^{1,23,4}
+ \varphi^{2,3,4},  
\end{equation}
\begin{align} \label{octogon:grtm}
& \psi(A | b(0),b(1),\ldots,b(N-1)) - \psi(A | b(1),b(2),\ldots,b(0))
\\ & \nonumber 
+ \psi(C | b(1),b(0),\ldots,b(2-N)) - \psi(C | b(0),b(-1),\ldots,b(1-N))=0,
\end{align}
where $C = - A - \sum_{i=0}^{N-1} b(i)$, 
\begin{equation} \label{add:cond:grtm}
[\psi(A|b(0),\ldots,b(N-1)) - \psi(C|b(0),\ldots,b(1-N)),C]
+ \sum_{a=0}^{N-1} [\psi(A|b(a),\ldots,b(a+N-1)),b(a)] =0
\end{equation}
where $C = - A - \sum_{i=0}^{N-1} b(i)$, 
\begin{equation} \label{pent:grtm}
\psi^{1,2,34} + \psi^{12,3,4} = \varphi^{1,2,3} + \psi^{1,23,4} + \psi^{1,2,3}.  
\end{equation}
The Lie bracket is given by 
\begin{equation} \label{magnus:0}
[(\varphi_1,\psi_1),(\varphi_2,\psi_2)] 
= (\langle\varphi_1,\varphi_2\rangle,\langle\psi_1,\psi_2\rangle),
\end{equation}
where 
\begin{equation} \label{magnus:1}
\langle \varphi_1,\varphi_2 \rangle = [\varphi_1,\varphi_2] 
+ D_{\varphi_2}(\varphi_1) - D_{\varphi_1}(\varphi_2)
\end{equation}
\begin{equation} \label{magnus:2}
\langle \psi_1,\psi_2\rangle = [\psi_1,\psi_2] 
+ \overline D_{\psi_2}(\psi_1) - \overline D_{\psi_1}(\psi_2).   
\end{equation}
Here $D_{\varphi}$, $\overline D_{\psi}$ are defined by 
\begin{equation} \label{D:phi}
D_\varphi(U) = [\varphi,U], \; 
D_\varphi(V) = 0
\end{equation} 
and 
\begin{equation} \label{D:psi}
\overline D_\psi(A) = [\psi,A], \; 
\overline D_\psi(b(a)) = [\psi(A | b(0),\ldots,b(N-1)) 
- \psi(A | b(a),\ldots,b(a+N-1)), b(a)]. 
\end{equation} 

The direct sum of its homogeneous components is a graded Lie subalgebra of 
$\grtm_{(\bar 1,1)}(N,\kk)$, and $\grtm_{(\bar 1,1)}(N,\kk)$ is then the corresponding 
graded completion.  

The Lie algebra $\grtm(N,\kk)$ is the semidirect product 
$\grtm_1(N,\kk) \rtimes \kk$, $1\in \kk$ acting by the 
continuous derivation taking an
element $(\varphi,\psi)$ of degree $n$ to $-n(\varphi,\psi)$. 
Then $\grtm(N,\kk)$ is equipped with an action of 
$(\ZZ/N\ZZ)^\times$. 

\begin{remark} 
The octogon condition (\ref{octogon:grtm}) means that 
$$
\psi(A | b(a),\ldots,b(a+N-1))
- \psi(C | b(a), \ldots,b(a+1-N))
$$ 
is independent of $a$. 
In the same way, the octogon condition (\ref{octogon:GRTM})
means that $k(A|b(a),\ldots,b(a+N-1))^{-1} 
k(C|b(a),b(a-1),\ldots,b(a+1-N))$ is independent of $a$. 
These facts are used when proving that $\on{GRTM}_{(\bar 1,1)}(N,\kk)$
is a group.  
\end{remark}

\begin{remark}
(\ref{add:cond:grtm}) implies (\ref{octogon:grtm}). Indeed,
(\ref{add:cond:grtm}) implies that 
$$
[\psi(A| b(0),\ldots,b(N-1)) - \psi(C|b(0),\ldots,b(1-N)),C]
$$ is invariant under the automorphism 
$\theta_1 : A\mapsto A$, $b(a)\mapsto b(a+1)$. So for some $\lambda\in\kk$, 
we have l.h.s. of (\ref{octogon:grtm})$ = \lambda C$. Applying 
$\id + \theta_1 + \cdots + \theta_1^{N-1}$ to this identity, 
we get $\lambda = 0$, whence (\ref{octogon:grtm}).  
\end{remark}

\begin{remark}
In \cite{Dr2}, 
Drinfeld showed that conditions (\ref{add:cond:grt}) and (\ref{add:cond:GRT})
can be removed from the definition of $\on{GRT}_1(\kk)$ and $\grt_1(\kk)$, so 
$\on{GRT}_1(\kk) = \ul{{\bf Assoc}}_0(\kk)$. On the other hand, we can only 
prove $\on{GRTM}_{(\bar 1,1)}(N,\kk) \hookrightarrow 
\ul{\Pseudo}_{(\bar 1,0)}(N,\kk)$.  \hfill \qed \medskip 
\end{remark}

As we noticed, $\gtm_{(\bar 1,1)}(N,\kk)$ is filtered by 
$\g^{\geq n} = \{(\varphi,\psi) | \psi$ has degree $\geq n\}$. 
 
\begin{lemma} \label{lemma:grading}
$\hat{\on{gr}}\gtm_{(\bar 1,1)}(N,\kk) \subset \grtm_{(\bar 1,1)}(N,\kk)$.
\end{lemma}

{\em Proof.}
Let $\Xi,\Theta,\eta(1),\ldots,\eta(N)$ be such that 
$e^\Theta e^{\eta(1)} \cdots e^{\eta(N)} e^\Xi = 1$. 
Let $\overline\psi\in \gtm_1(N,\kk)$ be of degree $\geq k$, 
and let $\psi$ be its degree $k$ component. Then for any $\alpha\in \ZZ$, we have 
\begin{align*}
& \overline \psi(\Xi | \eta(\alpha),\ldots, \eta(\alpha+N-1)) 
- \overline\psi(\Xi | \eta(\alpha+1),\ldots,\eta(\alpha+N))
\\ & 
= \Ad(e^{-\eta(\alpha)} \cdots e^{-\eta(1)}) 
\Big( \overline \psi(\Theta | \eta'(\alpha),\ldots,\eta'(\alpha+1-N)) 
- \overline\psi(\Theta | \eta'(\alpha+1),\ldots,\eta'(\alpha+2-N) )\Big)  
\end{align*}  
(we use the conventions of Lemma \ref{prop:gtm}).
Summing up these equalities for $\alpha=0,\ldots,N-1$, we get 
\begin{align*}
& (1 - \Ad(e^\Xi))\overline\psi(\Xi | \eta(0),\ldots,\eta(N-1))
\\ & = \sum_{\alpha=1}^{N-1} (\Ad(e^{-\eta(\alpha)}) -1) \Ad(e^{-\eta(\alpha-1)}\cdots 
e^{-\eta(1)}) \overline\psi(\Theta | \eta'(\alpha),\ldots,\eta'(\alpha+1-N))
\\ & + \big( 1 - \Ad(e^{\Xi} e^{\eta(0)})\big) 
\overline\psi(\Theta | \eta'(0),\ldots,\eta'(1-N)), 
\end{align*}
which implies that $\psi$ satisfies (\ref{add:cond:grtm}). 
The other conditions satisfied by $\overline\psi$ imply that 
$\psi \in \grtm_{(\bar 1,1)}(N,\kk)$. \hfill \qed \medskip 

\subsection{The action of $\on{GRTM}(N,\kk)$ on $\Pseudo(N,\kk)$}

$\on{GRTM}_{(\bar 1,1)}(N,\kk)$ acts on $\ul{\Pseudo}(N,\kk)$ from the right 
by $(\Phi,\Psi) * (h,k) = (\Phi',\Psi')$, 
where
\begin{equation} \label{act:grt}
\Phi'(U,V) = h(U,V)\Phi(U,\Ad(h(U,V))^{-1}(V)), 
\end{equation}
\begin{align} \label{act:grtm}
& \Psi'(A|b(0),\ldots,b(N-1)) = k(A|b(0),\ldots,b(N-1))
\\ & \nonumber
\Psi\Big(A | \Ad\big(k(A|b(0),\ldots,b(N-1))^{-1}\big)(b(0)), \ldots,
\Ad\big(k(A|b(N-1),\ldots,b(2(N-1)))^{-1}\big)(b(N-1)) \Big). 
\end{align}
This action preserves each $\ul{\Pseudo}_{(a,\lambda)}(N,\kk)$, and it
extends to an action of $\on{GRTM}(N,\kk)$ on 
$\ul{\Pseudo}(N,\kk)$, which is compatible with the action of $(\ZZ/N\ZZ)^{\times}
\times \kk^{\times}$ on $(\ZZ/N\ZZ)\times\kk$ and commutes with the left action of 
$\on{GTM}(N,\kk)$ on $\ul{\Pseudo}(N,\kk)$. 

\begin{theorem} 
$\Pseudo(N,\kk)$ (resp., $\Pseudo_{\bar 1}(N,\kk)$, $\Pseudo_{(\bar 1,1)}(N,\kk)$) 
is a torsor under $\on{GRTM}(N,\kk)$ (resp., $\on{GRTM}_{\bar 1}(N,\kk)$, 
$\on{GRTM}_{(\bar 1,1)}(N,\kk)$).  
\end{theorem}

{\em Proof.} It suffices to prove that last statement. Let 
$(\Phi_1,\Psi_1)$ and $(\Phi_2,\Psi_2)$
be elements of $\Pseudo_{(\bar 1,1)}(N,\kk)$. Then there exists a unique pair 
$(h,k) \in \exp(\hat\t_3^0) \times \exp(\hat\t^0_{3,N})$, such that 
$(\Phi_1,\Psi_1) * (h,k) = (\Phi_2,\Psi_2)$ (where the action is computed 
according to formulas (\ref{act:grt}), (\ref{act:grtm}). We have to show that 
$(h,k) \in \on{GRTM}_{(\bar 1,1)}(N,\kk)$.  

Assume that $(\Phi_1,\Psi_1)$ and $(\Phi_2,\Psi_2)$ coincide up to order
$\nu-1$, i.e., their difference belongs to 
${U(\t^0_{3})}_{\geq\nu}^{\wedge} \oplus {U(\t^0_{3,N})}_{\geq\nu}^{\wedge}$
(the index means the grading where each generator has degree $1$ and the $()^{\wedge}$
means degree completion). 
Set $\varphi := (\Phi_2 - \Phi_1$ mod ${U(\t^0_3)}_{\geq\nu+1}^{\wedge})$, 
$\psi := (\Psi_2 - \Psi_1$ mod ${U(\t^0_{3,N})}_{\geq\nu+1}^{\wedge})$. 
Since the $\Phi_i,\Psi_i$ are group-like and coincide up to order
$n-1$, we have $\varphi \in (\t^0_{3})_\nu$ and $\psi \in (\t^0_{3,N})_\nu$. 
Let us show that $(\varphi,\psi) \in \grtm_{(\bar 1,1)}(N,\kk)$. 

We first show that $\varphi\in \grt_1(\kk)$. Comparing identities 
(\ref{cond:Phi:1}) and (\ref{cond:Phi:2}) satisfied by $\Phi_1$
and $\Phi_2$, we find that $\varphi$ satisfies (\ref{cond:1:grt}), 
(\ref{cond:2:grt}). 
According to \cite{Dr2}, Proposition 5.7, this implies that 
$\varphi$ also satifies (\ref{add:cond:grt}) and therefore belongs 
to $\grt_1(\kk)$. 

Comparing now identities (\ref{cond:Psi:1}) and (\ref{cond:Psi:2}) 
satisfied by $(\Phi_1,\Psi_1)$
and $(\Phi_2,\Psi_2)$, we find that $\psi$ satisfies (\ref{octogon:grtm}) 
and  $(\varphi,\psi)$ 
satisfies (\ref{pent:grtm}). Let us 
show that $\psi$ also satisfies (\ref{add:cond:grtm}). If 
$\Psi(A | b(0),\ldots, b(N-1))$
satisfies (\ref{cond:Psi:1}), we have 
\begin{equation} \label{cyclic:id}
X(0)\cdots X(N-1) = 1, 
\end{equation} 
where
\begin{align} \label{def:X} 
& X(a) := 
\Psi(A | b(a),\ldots, b(a+N-1))^{-1} e^{Nb(a)/2} 
\Psi(C | b(a),\ldots, b(a+1-N)) e^C
\\ & \nonumber 
\Psi(C | b(a),\ldots, b(a+1-N))^{-1} e^{Nb(a)/2} 
\Psi(A | b(a),\ldots, b(a+N-1)) e^A. 
\end{align}  
Let $X_i(a)$ be the r.h.s. of (\ref{def:X}) for $\Psi = \Psi_i$
($i=1,2$). Then 
\begin{align*} 
& X_2(a) = X_1(a) + [\psi(A | b(A),\ldots,b(a+N-1))
- \psi(C | b(a),\ldots,b(a+1-N)),A] 
\\ & 
- N [\psi(C | b(a),\ldots,b(a+1-N)) ,b(a)] + \on{\ terms\ of\ degree}>\nu+1, 
\end{align*}
and $X_i(a) = 1+$ terms of degree $>0$. 

Therefore 
\begin{align*}
& (\on{l.h.s.\ of\ (\ref{cyclic:id})\ for\ }\Psi_2)
- (\on{l.h.s.\ of\ (\ref{cyclic:id})\ for\ }\Psi_1)
\\ & 
= \sum_{a=0}^{N-1} \big( 
[\psi(A|b(a),\ldots,b(a+N-1)) - \psi(C|b(a),\ldots,b(a+1-N)),A] 
\\ & - N[\psi(C | b(a),\ldots,b(a+1-N)),b(a)]\big) 
+ \on{terms\ of\ degree}>\nu+1. 
\end{align*}
Since $\Psi_1$ and $\Psi_2$ satisfy (\ref{cyclic:id}), we get 
\begin{align*}
\sum_{a=0}^{N-1} & \big( [\psi(A | b(a),\ldots,b(a+N-1)) - 
\psi(C | b(a),\ldots,b(a+1-N)),A] 
\\ & - N [\psi(C | b(a),\ldots,b(a+1-N)),b(a)]\big) = 0. 
\end{align*}

Since $\psi$ satisfies (\ref{octogon:grtm}), this identity is equivalent to 
\begin{align*}
& [\psi(A|b(0),\ldots,b(N-1)) - \psi(C | b(0),\ldots,b(1-N)) , A] 
\\ & - \sum_{a=0}^{N-1} [\psi(C | b(a),\ldots,b(a+1-N)),b(a)] = 0, 
\end{align*}
which is equivalent to (\ref{add:cond:grtm}), after the change of 
variables $A\mapsto C$, $C\mapsto A$, $b(a) \mapsto b(-a)$.  

Hence $(\varphi,\psi) \in \grtm_{(\bar 1,1)}(N,\kk)$. Set $(\Phi'_1,\Psi'_1) 
:= (\Phi_1,\Psi_1) * \on{Exp}(\varphi,\psi)$, where $\on{Exp} :
\grtm_{(\bar 1,1)}(N,\kk) \to \on{GRTM}_{(\bar 1,1)}(N,\kk)$ is the 
exponential map for the group structure of $\on{GRTM}_{(\bar 1,1)}(N,\kk)$
(it is a restriction of the exponential of $\on{Mag}_{(\bar 1,1)}(N,\kk)$). 
Then $(\Phi'_1,\Psi'_1)
\in \Pseudo_{(\bar 1,1)}(N,\kk)$ coincides with $(\Phi_2,\Psi_2)$ up to order $\nu$.
Using successive approximations, we then construct $(h,k)\in 
\on{GRTM}_{(\bar 1,1)}(N,\kk)$ such that $(\Phi_1,\Psi_1) * (h,k) 
= (\Phi_2,\Psi_2)$. 
\hfill \qed \medskip 

As in \cite{Dr2}, we get: 

\begin{corollary}
Any $(\Phi,\Psi)\in \Pseudo(N,\kk)$ gives rise to a group isomorphism 
$s_{(\Phi,\Psi)} : \on{GRTM}(N,\kk) \to \on{GTM}(N,\kk)$, 
defined by the condition that $(\Phi,\Psi) *g = s_{(\Phi,\Psi)}(g) * 
(\Phi,\Psi)$. The diagram 
$$
\begin{matrix}
\on{GRTM}(N,\kk) & \stackrel{s_{(\Phi,\Psi)}}{\longrightarrow} & 
\on{GTM}(N,\kk) \\ 
\downarrow       & & \downarrow  \\
 (\ZZ /N\ZZ)^\times \times \kk^\times &  \simeq & \kk(N)^{\times}          
\end{matrix}
$$
is commutative, so $s_{(\Phi,\Psi)}(\on{GRTM}_{(\bar 1,1)}(N,\kk)) 
= \on{GTM}_{(\bar 1,1)}(N,\kk)$. 
The resulting isomorphism $\grtm_{bar 1,1)}(N,\kk) \to 
\gtm_{bar 1,1)}(N,\kk)$ is filtered, and the inverse of the 
associated graded isomorphism coincides with the map defined 
in Lemma \ref{lemma:grading}. The splitting of the 
exact sequence (\ref{exact:seq}) corresponding to
$(\Phi,\Psi)$ is the morphism 
$\on{Lie}(s_{(\Phi,\Psi)}) \circ i : \kk \to \gtm_{\bar 1}(N,\kk)$, where $i$ 
is the canonical embedding $\kk\to \grtm_{\bar 1}(N,\kk)$. 
\end{corollary}

\subsection{$\ul{\Pseudo}(N,\kk)$ when $N=1$} 

\begin{lemma}
There is a unique isomorphism $\ul\Assoc(\kk) \times \kk \to \ul\Pseudo(1,\kk)$, 
taking $((\lambda,\Phi),\gamma)$ to $(\bar 1,\lambda,\Phi,\Psi)$, where 
$\Psi = e^{\gamma t^{23}} \Phi$. 
\end{lemma}

{\em Proof.} Let $(\bar 1,\lambda,\Phi,\Psi)$ belong to $\ul{\Pseudo}(1,\kk)$. 
Let $\eps^{(0)} : \t_3 \to \t_2$ be the Lie algebra morphism 
$t^{0i}\mapsto 0$, $t^{12}\mapsto t^{12}$. Then there exists $\gamma\in \kk$
such that $\eps^{(0)}(\Psi) = e^{\gamma t^{12}}$. Applying $\eps^{(0)}$
to (\ref{cond:Psi:2}) and taking into account that $\Phi,\Psi$ are 
group-like, we get $\Psi = e^{\gamma t^{12}} \Phi$. 
Conversely, if $\Phi \in \ul{\Assoc}_\lambda(\kk)$ and $\gamma\in\kk$, then 
$(\Phi,e^{\gamma t^{23}}\Phi)$ belongs to 
$\ul{\Pseudo}_{(\bar 1,\lambda)}(\kk)$. \hfill \qed \medskip 

\begin{lemma}
There exists a unique isomorphism $\grtm_{(\bar 1,1)}(1,\kk) \simeq 
\grt_{(\bar 1,1)}(\kk) \oplus \kk (0,t^{12})$, where $(0,t^{12})$ is 
central and has degree $1$.  
\end{lemma} 

{\em Proof.} Similarly to the previous lemma, one proves that $(\varphi,\psi) 
\in \grtm_{(\bar 1,1)}(1,\kk)$ iff $\varphi\in \grt_{(\bar 1,1)}(\kk)$
and there exists $\gamma\in\kk$ such that $\psi = \varphi + \gamma t^{12}$. 
One then checks that $(0,t^{12})$ is central. 

\begin{remark} {\it $\grtm_{(\bar 1,1)}(N,\kk)$ in degree $1$.}
One checks that the degree $1$ part of $\grtm_{(\bar 1,1)}(N,\kk)$ is spanned by the
$(0,b(a) + b(-a))$, where $a\in (\ZZ/N\ZZ) / \{\pm 1\}$. Moreover, 
$(0,b(0))$ is central in $\grtm_{(\bar 1,1)}(N,\kk)$.  
\end{remark}

\section{The equality $\on{GTK} = \wh{\on{GT}}$}\label{sec:ih}

Recall that $\varphi_N: F_2\to \ZZ/N\ZZ$ is the morphism $x\mapsto \bar 1$, 
$y\mapsto \bar 0$ and that $\on{Ker}\varphi_N \subset F_2$ is freely 
generated by $X,y(0),...,y(N-1)$, where $X:= x^N$, 
$y(\alpha):= x^\alpha yx^{-\alpha}$. Define a group morphism 
$\delta_N : \on{Ker}\varphi_N\to F_2$
by $X\mapsto x$, $y(0)\mapsto y$, $y(\alpha)\mapsto 1$ for $\alpha\in[1,N-1]$. 
Define also a group morphism $c_N : \on{Ker}\varphi_N \to \ZZ$
by $c_N:= c \circ \delta_N$, where $c:F_2\to \ZZ$ is $x\mapsto 0$, 
$y\mapsto 1$. 

Recall that for any $(\lambda,f)\in\wh{\ul{\on{GT}}}$, $f\in \wh{F}_2'
\subset (\on{Ker}\varphi_N)^\wedge$. In \cite{Ih2}, Ihara defined 
$\on{GTK}\subset \wh{\on{GT}}$ as the set of $(\lambda,f)$, 
such that 
\begin{equation} \label{distf}
\forall N\geq 1, \;  \delta_N(f)=y^{c_N(f)}f. 
\end{equation}

\begin{theorem} \label{thm:ihara}
$\on{GTK} = \wh{\on{GT}}$; more generally, if 
$(\lambda,f)\in\wh{\ul{\on{GT}}}$, then (\ref{distf}) holds. 
\end{theorem}

{\em Proof.} Set $W_n:= W_{n,1}= \{(Z_1,...,Z_n)\in (\CC^\times)^n |
\forall i\neq j, Z_i\neq Z_j\}$ and recall that $P_{n+1} = \pi_1(W_n,b)$
(see Subsection \ref{sec:phi:nN}). The morphism $\varphi_{n,N}: P_{n+1}\to (\ZZ/N\ZZ)^n$
is induced by the covering $W_{n,N}\to W_n$, $(z_1,...,z_n)\mapsto 
(z_1^N,...,z_n^N)$, where $W_{n,N}=\{(z_1,...,z_n)\in(\CC^\times)^n | 
\forall i\neq j,\forall \zeta\in \mu_N(\CC), z_i\neq \zeta z_j\}$.
Hence $\on{Ker}\varphi_{n,N} = \pi_1(W_{n,N},b)$ (it is denoted $K_{n,N}$). 
In addition to the inclusion $K_{n,N}\hookrightarrow P_{n+1}$, 
we have a morphism $\delta_N : K_{n,N}\to P_{n+1}$, induced by 
$W_{n,N}\to W_n$, $(z_1,...,z_n)\mapsto (z_1,...,z_n)$ by taking fundamental groups. 

The morphisms $\partial_i:P_n\to P_{n+1}$ ($i=0,...,n+1$) are given 
(in terms of the Artin generators $x_{jk}$, $0\leq j<k\leq n-1$ of $P_n$)
by $\partial_0:x_{jk}\mapsto x_{j+1,k+1}$, $\partial_{n+1}:x_{jk}\mapsto 
x_{jk}$ and $\partial_i:x_{jk}\mapsto x_{jk}$ if $k<i-1$, 
$x_{jk}\mapsto x_{j,k+1}$ if $k>i-1>j$, $x_{jk}\mapsto x_{j+1,k+1}$ if $j>i-1$, 
$x_{j,i-1}\mapsto x_{j,i-1}x_{ji}$ if
$j<i-1$, $x_{i-1,j}\mapsto x_{i-1,j+1}x_{i,j+1}$ if $j>i-1$. 

\begin{lemma} There exist morphisms $\tilde\partial_i : 
P_n \to P_{n+1}$ ($i=1,...,n$) given by the same formulas, except 
$\tilde\partial_{i}:x_{0,i-1}\mapsto x_{0,i-1}x_{0,i}x_{i-1,i}$
if $i=2,...,n$ and $\tilde\partial_{1}:x_{01}\mapsto x_{01}x_{02}x_{12}$. 
\end{lemma}

{\em Proof of Lemma.}
One can show that for $i=1,...,n$, the morphism $\partial_i: P_{n}\to P_{n+1}$ 
is induced by the sequence of maps $W_{n-1} \supset W_{n-1}^{i,\eps}
\to W_n$, where $\eps\in]0,1[$, and: (a) for $i=2,...,n$, 
$W_{n-1}^{i,\eps}:=\{(z_1,...,z_{n-1})\in W_{n-1} | 
\forall j\neq i, |z_j/z_{i-1} -1|>\eps\}$ and $W_{n-1}^{i,\eps}\to W_{n}$ is  
$(z_1,...,z_{n-1})\mapsto (z_1,...,z_{i-1},z_{i-1}+\eps|z_{i-1}|,z_i,...,z_{n-1})$, 
and (b) for $i=1$, $W_{n-1}^{1,\eps} = \{(z_{1},...,z_{n})\in W_{n} | 
\forall i, |z_{i}/z_{1}|>\eps\}$ and $W_{n-1}^{1,\eps}\to 
W_{n}$ is $(z_{1},...,z_{n-1})\mapsto (\eps |z_{1}|,z_{1},...,z_{n-1})$. 

$\tilde\partial_i : P_n\to P_{n+1}$ ($i=1,...,n$) is then defined as 
the morphism induced by the sequence of maps $W_{n-1} \supset W_{n-1}^{i,\eps}
\to W_n$, where $W_{n-1}^{i,\eps}\to W_n$ is  
$(z_1,...,z_{n-1})\mapsto (z_1,...,z_{i-1},z_{i-1}(1+\eps),z_{i},...,z_{n-1})$
for $i=2,...,n$ and $(z_{1},...,z_{n-1})\mapsto (\eps z_{1},z_{1},...,z_{n-1})$ 
for $i=1$. One checks that the additional possible winding of $i$ around 
$i+1$ corresponds to the above explicit formulas for $\tilde\partial_i$. 
\hfill \qed \medskip 

On the other hand, $\partial_0 : P_{n}\to P_{n+1}$ is induced by the inclusion 
$\{(z_1,...,z_{n})|\forall i,z_{i}\notin\RR_{-}$ and $\forall i\neq j, z_i\neq z_j\}
\subset W_{n}$ and  
$\partial_{n+1} : P_{n}\to P_{n+1}$ is induced by the sequence of maps 
$W_{n-1} \supset W_{n-1} \cap {\bf D}^{n-1} \to W_{n}$, where 
${\bf D}=\{z\in\CC|
|z|<1\}$ and $W_{n-1} \cap {\bf D}^{n-1} \to W_{n}$ is $(z_1,...,z_{n-1})
\mapsto (z_1,...,z_{n-1},1)$.

\begin{lemma} \label{lemma:rest}
The morphisms $\tilde\partial_1,...,\tilde\partial_n,
\partial_{n+1}$ restrict to morphisms $K_{n-1,N} \to K_{n,N}$, and
$\partial_0$ factors as $P_n\to K_{n,N}\subset P_{n+1}$. 
\end{lemma}
 
{\em Proof.} Recall that $\varphi_{n,N} : P_{n+1}\to (\ZZ/N\ZZ)^n$
is given by $x_{0i}\mapsto s_i$, the $i$th generator of 
$(\ZZ/N\ZZ)^n$. 
We have commutative diagrams $\begin{matrix}
P_n & \stackrel{\varphi_{n-1,N}}{\to}& (\ZZ/N\ZZ)^{n-1} \\
\scriptstyle{\tilde\partial_i,\partial_{n+1}}\downarrow & & \downarrow 
\scriptstyle{\tilde\alpha_i,\alpha_{n+1}}\\ 
P_{n+1} & \stackrel{\varphi_{n,N}}{\to}& (\ZZ/N\ZZ)^n \end{matrix}$, 
where for $i=2,...,n$, $\tilde\alpha_i(s_j)=s_j$ for $j<i-1$, 
$\tilde\alpha_i(s_j)=s_{j+1}$ for $j>i-1$, $\tilde\alpha_i(s_{i-1}) = 
s_{i-1}+s_i$, $\tilde\alpha_1 = \tilde\alpha_2$, and 
$\alpha_{n+1}(s_j) = s_j$ for any $j$. On the other hand, 
$\varphi_{n,N}\circ \partial_0=0$. This implies the result. 
\hfill \qed \medskip 

In view of the topological interpretation of $\partial_0,
\tilde\partial_{1,...,n},\partial_{n+1}$, the morphisms 
$\tilde\partial_{1,...,n},\partial_{n+1} : K_{n-1,N} \to K_{n,N}$
and $\partial_0 : P_n\to K_{n-1,N}$ constructed in Lemma \ref{lemma:rest}
can be interpreted as induced by the following sequences of 
maps: (a) in the case of $(\tilde\partial_i)_{i=2,...,n}$, the 
maps\footnote{Here 
$\eps>0$ is small enough for $\{u\in\CC| |u^N-1|<\eps\}$
to consist of disjoint contractible neighborhoods of the 
$\zeta\in\mu_N(\CC)$ in $\CC^{\times}$.} 
$W_{n-1,N} \supset W_{n-1,N}^{i,\eps}\to W_{n,N}$, 
$W_{n-1,N}^{i,\eps}:=\{(z_1,...,z_n)\in W_{n-1,N}|
\forall j\neq i-1,|(z_j/z_{i-1})^N-1|>\eps\}$ and 
$W_{n-1,N}^{i,\eps}\to W_{n,N}$ is $(z_1,...,z_{n-1})
\mapsto (z_1,...,z_{i-1}(1+\eps)^{1/N},...,z_{n-1})$; 
(b) in the case of $\tilde\partial_1$, the maps 
$W_{n-1,N} \supset W_{n-1,N}^{1,\eps}\to W_{n,N}$, 
where $0<\eps<1$, $W_{n-1,N}^{1,\eps} = 
\{(z_1,...,z_{n-1})\in W_{n-1,N}|\forall j\neq 1, |(z_j/z_1)^N|>\eps\}$
and $W_{n-1,N}^{1,\eps}\to W_{n,N}$ is $(z_1,...,z_{n-1})\mapsto 
(\eps^{1/N} z_1,z_1,...,z_{n-1})$; (c) in the case of $\partial_{n+1}$, 
the maps $W_{n-1,N}\supset  W_{n-1,N}\cap {\bf D}^{n-1} \to W_{n,N}$, 
where $W_{n-1,N}\cap {\bf D}^{n-1} \to W_{n,N}$ is $(z_1,...,z_{n-1})
\mapsto (z_1,...,z_{n-1},1)$; (d) in the case of $\partial_0$, the map 
$\{(z_1,...,z_n)\in W_n | \forall i, \on{arg}(z_i) < \pi/n\} \to W_{n,N}$, 
$(z_1,...,z_n)\mapsto (z_1,...,z_n)$. 

\begin{lemma}
The diagrams $\begin{matrix} P_n & \hookleftarrow & 
K_{n-1,N}& \stackrel{\delta_N}{\to}& P_{n} 
\\ \scriptstyle{\tilde\partial_i,\partial_{n+1}}\downarrow && 
\scriptstyle{\tilde\partial_i,\partial_{n+1}}\downarrow & &
\downarrow\scriptstyle{\tilde\partial_i,\partial_{n+1}} \\ 
P_{n+1} &\hookleftarrow & K_{n,N}& \stackrel{\delta_N}{\to}& P_{n+1}
\end{matrix}$  and 
$\begin{matrix} P_n &\stackrel{\on{id}}{\leftarrow} 
&P_n & \stackrel{\on{id}}{\to}& P_{n} 
\\ \scriptstyle{\partial_{0}}\downarrow &&
\scriptstyle{\partial_{0}}\downarrow & &
\downarrow\scriptstyle{\partial_{0}} \\ 
P_{n+1}&\hookleftarrow &K_{n,N}& \stackrel{\delta_N}{\to}& P_{n+1}
\end{matrix}$  
commute. 
\end{lemma}

{\em Proof.} We have already seen that the left part of each diagram 
commutes, so it remains to show the right parts. For this, we rely on the 
above topological interpretations. 

Let $i\in \{1,...,n\}$. If $i=2,...,n$, let $\eps>0$ be such that 
$\{u\in\CC| |u^N-1|<\eps\}$ consists of disjoint contractible neighborhoods 
of the $\zeta\in\mu_N(\CC)$; and let $0<\eps<1$ if $i=1$;  then the diagram 
$$
\begin{matrix}
W_{n-1,N} &\supset & W_{n-1,N}^{i,\eps} & \to & W_{n,N} \\ 
\downarrow && \downarrow && \downarrow \\ 
W_{n-1} &\supset & W_{n-1}^{i,\eps} &\to & W_n 
\end{matrix}
$$
commutes, where the vertical maps are all restrictions of
$(z_1,...,z_k)\mapsto (z_1^N,...,z_k^N)$ ($k=n-1,n$). 
Passing to fundamental groups, we see that the right diagram involving
$\tilde\partial_i$ commutes. 

When $i=n+1$, the diagram 
$$
\begin{matrix}
W_{n-1,N} &\supset & W_{n-1,N} \cap {\bf D}^{n-1} & \to & W_{n,N} \\ 
\downarrow && \downarrow && \downarrow \\ 
W_{n-1} &\supset & W_{n-1} \cap {\bf D}^{n-1} &\to & W_n 
\end{matrix}
$$
commutes, where the vertical maps are as above. 
This implies the wanted result involving $\partial_{n+1}$. 

The commutativity of the diagram involving $\partial_0$ similarly 
follows from that of

\noindent  
$\begin{matrix} \{(z_1,...,z_n)\in W_n | \forall i, 
\on{arg}(z_i) < \pi/n\} & \to & W_{n,N}
\\ \downarrow & & \downarrow \\ \{(Z_1,...,Z_n)\in W_n | \forall i, 
Z_i\notin \RR_-\} & \to & W_{n} \end{matrix}$. \hfill \qed \medskip 

\begin{lemma}
Let $f\in \wh{F}'_2$. Then for $i=1,2,3$, $\tilde\partial_i(f)
=\partial_i(f)$. 
\end{lemma}

{\em Proof.} For $f\in \wh F_2$, recall that $c(f)\in\wh\ZZ$ is defined by 
$f(1,y)=y^{c(f)}$ and define $b(f)\in\wh\ZZ$ by $f(x,1)=x^{b(f)}$. 
As $\tilde\partial_1(f)=f(x_{01}x_{02}x_{12},x_{23})$ and since $x_{01}$
commutes with both $x_{02}x_{12}$ and $x_{23}$, we get 
$\tilde\partial_1(f)=x_{01}^{b(f)}\partial_1(f)$. Similarly, 
$\tilde\partial_2(f) = f(x_{01}x_{02}x_{12},x_{13}x_{23})$, 
and since $x_{12}$ commutes with both $x_{01}x_{02}$ and $x_{13}x_{23}$, 
we get $\tilde\partial_2(f) = x_{12}^{b(f) }\partial_2(f)$; and
$\tilde\partial_3(f) = f(x_{01},x_{13}x_{23})$ implies 
$\tilde\partial_3(f) = \partial_3(f)$. If $f\in
\wh F_2'$, then $b(f)=c(f)=0$, which implies the result. 
\hfill \qed \medskip 

Let $(\lambda,f)\in\wh{\ul{\on{GT}}}$. Then $f\in \wh F_2'$, hence 
\begin{equation} \label{pent:ih}
\partial_0(f)\tilde\partial_2(f)\partial_4(f) = \tilde\partial_3(f)
\tilde\partial_1(f). 
\end{equation}
$\partial_0(f)$ belongs to $K_{4,N}$, and since 
$f\in (\on{Ker}\varphi_N)^\wedge\subset \wh K_{3,N}$, so do 
$\tilde\partial_{1,2,3}(f)$ and $\partial_4(f)$. 
So (\ref{pent:ih}) takes place in $\wh K_{3,N}$. Let us 
apply $\delta_N$ to this identity. In view of the above 
commutation relations, we get $\partial_0(f) \tilde\partial_2(h)
\partial_4(h) = \tilde\partial_3(h)\tilde\partial_1(h)$, where
$h=\delta_N(g)$; this equality takes place in $\wh P_4$.
This is rewritten $\partial_{0}(f)x_{12}^{b(h)}\partial_{2}(h)
\partial_{4}(h)=\partial_{3}(h)x_{01}^{b(h)}\partial_{1}(h)$. 
For a general $k\in \on{Ker}\varphi_{N}$, we have 
$b(k) = Nb(\delta_{N}(k))$ so $Nb(h)=b(f)=0$ as $f\in\wh F_{2}'$, so 
$b(h)=0$. Hence 
$$
\partial_{0}(f)\partial_{2}(h)
\partial_{4}(h)=\partial_{3}(h)\partial_{1}(h). 
$$
Applying now the morphism $\eps_0:\wh P_4\to \wh P_3$
of erasing of the $0$th strand (i.e., $x_{0i}\mapsto 1$, 
$x_{ij}\mapsto x_{ij}$ for $1\leq i<j\leq 3$), and using the fact 
that $\eps_0\circ \partial_0 = \eps_0\circ \partial_1 = \id$, 
$\eps_0\circ \partial_2(h) = (x_{13}x_{23})^{c(h)}$, 
$\eps_0\circ \partial_3(h) = (x_{12}x_{13})^{c(h)}$, 
$\eps_0\circ \partial_4(h) = x_{12}^{c(h)}$ for $h\in \wh F_2$, 
we get $f \cdot (x_{13}x_{23})^{c(h)} \cdot x_{12}^{c(h)}
= (x_{12}x_{13})^{c(h)}\cdot h$, hence $h = (x_{12}x_{13})^{-c(h)}
\cdot f \cdot (x_{13}x_{23})^{c(h)} x_{12}^{c(h)}= x_{23}^{c(h)}f$ as 
$x_{12}x_{13}x_{23}$ is central in $\wh P_{3}$, i.e., $h(x,y)=y^{c(h)}f(x,y)$, as wanted. 
\hfill \qed \medskip 

\begin{remark}
It is observed in \cite{Ih2} that Theorem \ref{thm:ihara} implies that 
$c_{N_{1}N_{2}}(f)=c_{N_{1}}(f)+c_{N_{2}}(f)$ for 
$(\lambda,f)\in \wh{\ul{\on{GT}}}$, and $N_{1},N_{2}\geq 1$; indeed, 
$y^{c_{N_{1}N_{2}}(f)}f=\delta_{N_{1}N_{2}}(f)
=\delta_{N_{1}}(\delta_{N_{2}}(f)) = 
\delta_{N_{1}}(y^{c_{N_{2}}(f)}f)=y^{c_{N_{1}}(f)+c_{N_{2}}(f)}f$. 
\end{remark}

\begin{remark}
The analogue of Theorem \ref{thm:ihara} for $\ul{\on{GT}}_{l}$ can be 
formulated as follows. For $\nu\geq 0$, one can show that $\on{Ker}[(F_{2})_{l}
\stackrel{(\varphi_{l^{\nu}})_{l}}{\to}\ZZ/l^{\nu}\ZZ] = 
(\on{Ker}\varphi_{l^{\nu}})_{l}$ (if $N\geq 1$, then $(\varphi_{N})_{l}=
(\varphi_{l^{\nu}})_{l}$, where $\nu$ is the $l$-adic valuation of $N$); as 
$\on{Ker}\varphi_{l^{\nu}}$ is freely generated by $X:= x^{l^{\nu}}$ and 
$y(\alpha):= x^{\alpha}yx^{-\alpha}$, $\alpha\in [0,l^{\nu}-1]$, we have a morphism 
$\delta_{l^{\nu}} : (\on{Ker}\varphi_{l^{\nu}})_{l}\to (F_{2})_{l}$, 
$X\mapsto x$, $y(0)\mapsto y$, $y(\alpha)\mapsto 1$, $\alpha\in [1,l^{\nu}-1]$. 
We then define $c_{l^{\nu}} : (\on{Ker}\varphi_{l^{\nu}})_{l}\to \ZZ_{l}$
as $c_{l^{\nu}}:= c \circ \delta_{l^{\nu}}$, where $c:(F_{2})_{l}\to \ZZ_{l}$
is the morphism given by $x\mapsto 1$, $y\mapsto 0$. The analogue of Theorem 
\ref{thm:ihara} then says that for any $(\lambda,f)\in \ul{\on{GT}}_{l}$, we have 
$\delta_{l^{\nu}}(f)=y^{c_{l^{\nu}}(f)}f$, which implies that for any $(\lambda,f)
\in \ul{\on{GT}}_{l}$, we have $c_{l^{\nu}}(f)=\nu c_{l}(f)$. 
\end{remark}

\section{The distribution subgroup $\on{GTMD}(N,\kk)\subset 
\on{GTM}(N,\kk)$}

\subsection{The endomorphisms $\delta_{N}$ of $\wh{\on{GTM}}$}
\label{sec:9:1}

For $\ul f = (\lambda,f)\in {\wh{\on{GT}}}$, set $\chi(\ul f):= \lambda$; 
$\chi : {\wh{\on{GT}}}\to \wh\ZZ^{\times}$ is the cyclotomic character. Set 
$c_{N}(\ul f):= c_{N}(f)$. Then $c_{N}$ satisfies the cocycle identity 
$c_{N}(\ul f_1\ul f_{2}) = c_{N}({\ul f}_{1})\chi({\ul f}_{2})
+c_{N}({\ul f}_{2})$
(more generally, this identity holds for ${\ul f}_{i}\in \ul{\wh{\on{GT}}}$, 
if $\chi({\ul f}_{2})$ mod $N$ is in $(\ZZ/N\ZZ)^{\times}$). 
It follows that we have a group endomorphism $\delta_{N}$ of
$\wh{\on{GT}}\ltimes\wh\ZZ$, $\delta_{N}(\ul f,s):= 
(\ul f,s+c_{N}(\ul f))$. 

According to Theorem \ref{thm:ihara}, this endomorphism corresponds, 
under the isomorphism $\wh{\on{GTM}}\simeq \wh{\on{GT}}\ltimes
\wh\ZZ$ (Proposition \ref{prop:str:GTM}), to the self-map $(\lambda,\lambda,f,g)\mapsto 
(\lambda,\lambda,f,\delta_{N}(g))$ of $\on{GTM}$. 
(According to Proposition \ref{prop:str:GTM}, $g$ has the form $y^{s}f$, so since 
$f\in \wh F_{2}'$, $g\in \on{Ker}\hat b\subset (\on{Ker}\varphi_{N})^{\wedge}$, 
where $g(x,1)=x^{\hat b(g)}$, so 
$\delta_{N}(g)$ is well-defined.) Of course, the diagram $\begin{matrix}
\wh{\on{GT}} & \to & \wh{\on{GTM}}\\
  & \searrow & \downarrow\scriptstyle{\delta_{N}}\\
   & & \wh{\on{GTM}}\end{matrix}$ commutes.

\subsection{The morphisms $\delta_{NN'} : \ul{\on{GTM}}(N)_{l}
\to \ul{\on{GTM}}(N')_{l}$}

Let us assume that $N'|N$; we set $d:= N/N'$. In this subsection, we construct 
 morphisms $\delta_{NN'} : \ul{\on{GTM}}(N)_{l}\to \ul{\on{GTM}}(N')_{l}$, 
 such that $\begin{matrix} \ul{\wh{\on{GTM}}}& \to & \ul{\on{GTM}}(N)_{l}\\
 \scriptstyle{\delta_{d}}\downarrow & & \downarrow
 \scriptstyle{\delta_{NN'}}\\
 \ul{\wh{\on{GTM}}}  & \to & \ul{\on{GTM}}(N')_{l}\end{matrix}$ commutes. 

We first define $\delta_{NN'} : K_{n,N}\to K_{n,N'}$ to be the group 
morphism induced by the map $W_{n,N}\to W_{n,N'}$, $(z_{1},...,z_{n})
\mapsto (z_{1}^{d},...,z_{n}^{d})$ by taking fundamental groups. 
One checks that this map restricts to $\delta_{NN'}:\on{Ker}\varphi_{N}
\to \on{Ker}\varphi_{N'}$ under $\on{Ker}\varphi_{\nu}\subset K_{2,\nu}$
($\nu = N,N'$). 

Recall the ring morphism $[[-]]:\ZZ_{l}(N)\to \ZZ_{l}(N')$ given by 
$(a,r)\mapsto (\bar{\bar a},\tilde{\tilde a}+dr)$, where $a\mapsto 
\bar{\bar a}$ is the morphism $\ZZ/N\ZZ\to \ZZ/d\ZZ$, $\bar 1\mapsto \bar 1$, 
$\tilde{\tilde a}\in [0,N'-1]$ is the integral part of $\tilde a/d$, where
$\tilde a\in [0,N-1]$ is the lift of $a$; here $d=N/N'$. 

Let $(\lambda,\mu,f,g)\in \ul{\on{GTM}}(N)_{l}$. Recall that 
$(\lambda,\mu,f,g)\in (2\ZZ_{l}+1)\times \ZZ_{l}(N)\times (F_{2})_{l}
\times (\on{Ker}\varphi_{N})_{l}$. We then set $\delta_{NN'}(\lambda,\mu,f,g)
:= (\lambda,[[\mu]],f,\delta_{NN'}(g))$. 

\begin{proposition}
This defines a map $\delta_{NN'}:\ul{\on{GTM}}(N)_{l}
\to \ul{\on{GTM}}(N')_{l}$, such that $\delta_{NN'}(\ul{g}_{1}
\ul{g}_{2}) = \delta_{NN'}(\ul{g}_{1})\delta_{NN'}(\ul{g}_{2})$
for $\ul{g}_{2}\in \ul{\on{GTM}}(N)_{l}\times_{\ZZ/d\ZZ}(\ZZ/d\ZZ)^{\times}$, 
which restricts to a group morphism $\delta_{NN'}:\on{GTM}(N)_{l}\to 
\on{GTM}(N')_{l}$. 
\end{proposition}

{\em Proof.} One checks that for a general $g\in F_2(\varphi_N)$, 
$\tilde\partial_1(g) = x_{01}^{b_N(g)}\partial_1(g)$, 
$\tilde\partial_2(g) = x_{12}^{b_N(g)}\partial_2(g)$, 
$\tilde\partial_3(g) = \partial_2(g)$, where 
$\partial_i,\tilde\partial_i$ is shorthand for the 
morphisms $F_2(\varphi_N,l) \to P_4(\varphi_{3,N},l)$, 
$(F_2)_l \to P_4(\varphi_{3,N},l)$ induced
by their homonyms and $b_N:F_2(\varphi_N,l)\to \ZZ_l(N)$ is the map from 
Lemma \ref{lemma:b0}.

Let $(\lambda,\mu,f,g)\in \ul{\on{GTM}}(N)_{l}$. Let us 
show that $(f,\delta_{NN'}(g))$ satisfies the mixed pentagon relation. 
By Lemma \ref{lemma:b0}, $b_N(g)=0$ so the mixed pentagon relation 
satisfied by $(f,g)$ can be rewritten 
\begin{equation} \label{mpn}
\tilde\partial_3(g)\tilde\partial_1(g) = \partial_0(f)\tilde\partial_2(g)
\partial_4(g). 
\end{equation}
As in Section \ref{sec:ih}, one can show that the diagrams
$\begin{matrix} K_{n-1,N} & \stackrel{\delta_{NN'}}{\to}& K_{n-1,N'}\\
\scriptstyle{\tilde\partial_{1,...,n},\partial_{n+1}}\downarrow &&
\downarrow \scriptstyle{\tilde\partial_{1,...,n},\partial_{n+1}}\\
K_{n,N} & \stackrel{\delta_{NN'}}{\to}& K_{n,N'} \end{matrix}$ and 
$\begin{matrix} P_{n}&\stackrel{\partial_{0}}{\to}& K_{n,N} \\
&\scriptstyle{\partial_{0}}\searrow & \downarrow\scriptstyle{\delta_{NN'}}
\\ && K_{n,N'}\end{matrix}$ commute. 

Applying $\delta_{NN'}$ to (\ref{mpn}), and with $h:= \delta_{NN'}(g)$, we get 
$\tilde\partial_3(h)\tilde\partial_1(h) = \partial_0(f)\tilde\partial_2(h)
\partial_4(h)$ (the $\tilde\partial_i,\partial_i$ are now morphisms 
$(F_2)_l \to P_4(\varphi_{3,N'},l)$, 
$F_2(\varphi_{N'},l) \to P_4(\varphi_{3,N'},l)$. 
The restriction of $b_N$ to $(\on{Ker}\varphi_N)_l$ factors as 
$(\on{Ker}\varphi_N)_l \stackrel{\tilde b_N}{\to} \ZZ_l \stackrel{\lambda
\mapsto (\bar 0,\lambda)}{\hookrightarrow} \ZZ_l(N)$, and the diagram 
$\begin{matrix}
(\on{Ker}\varphi_N)_l & \stackrel{\tilde b_N}{\to}& \ZZ_l & \stackrel{\lambda
\mapsto (\bar 0,\lambda)}{\hookrightarrow}& \ZZ_l(N) \\ 
\scriptstyle{\delta_{NN'}}\downarrow && \scriptstyle{\id}\downarrow
&& \\ 
(\on{Ker}\varphi_{N'})_l & \stackrel{\tilde b_{N'}}{\to}& 
\ZZ_l & \stackrel{\lambda
\mapsto (\bar 0,\lambda)}{\hookrightarrow}& \ZZ_l(N')   
\end{matrix}$ commutes; $b_N(g)=0$ then implies that $b_{N'}(h)=0$. 
It follows that $\tilde\partial_i(h) = \partial_i(h)$ for 
$i=1,2,3$, so $(f,h)$ satisfies the mixed pentagon equation 
$\partial_3(h)\partial_1(h) = \partial_0(f)\partial_2(h)\partial_4(h)$. 

Let us show that $([[\mu]],\delta_{NN'}(g))$ satisfies the octogon relations. 
The inclusion $(K_{n,N})_{l}\hookrightarrow P_{n+1}(\varphi_{N},l)$
factors through $(K_{n,N})_{l}\hookrightarrow K_{n,d}(\varphi_{N},l)
\hookrightarrow P_{n+1}(\varphi_{N},l)$, and we have a commuting square 
$\begin{matrix} K_{n,N}& \stackrel{\delta_{NN'}}{\to} & K_{n,N'}\\
\scriptstyle{\subset}\downarrow 
& & \downarrow\scriptstyle{\subset}\\
K_{n,d}  & \stackrel{\delta_{d}}{\to}& P_{n+1} \end{matrix}$, which gives rise 
to a commuting square 
$\begin{matrix} (K_{n,N})_{l}& \stackrel{\delta_{NN'}}{\to} & (K_{n,N'})_{l}\\
\scriptstyle{\subset}\downarrow & & \downarrow\scriptstyle{\subset}\\
K_{n,d}(\varphi_{N},l)  & \stackrel{\delta_{d}}{\to}& P_{n+1}(\varphi_{N'},l) 
\end{matrix}$. For $g\in F_{2}(\varphi_{N},l)$, we write $g(x,y)=g$, $g(x^{-1}y^{-1},y)
= \kappa(g)$, etc. 

Rewriting the octogon identity as 
$$
g(x,y)^{-1}y^{{{\lambda-1}\over 2}} g(x^{-1}y^{-1},y)(x^{-1}y^{-1})^{\mu}x
g(x^{-2}y^{-1}x,x^{-1}yx)^{-1}(x^{-1}yx)^{{{\lambda+1}\over 2}}g(x,x^{-1}yx)
= x^{1-\mu}
$$
(equality in $F_{2}(\varphi_{N},l)$) and multiplying this equation by its conjugates
by $x^{-1},...,x^{1-d}$, we get 
\begin{align} \label{pre:octogon}
& g(x,y)^{-1} y^{{{\lambda-1}\over 2}} g(x^{-1}y^{-1},y) 
\on{Ad}[(x^{-1}y^{-1})^{\mu}x g(x^{-2}y^{-1}x,x^{-1}yx)]
[(x^{-1}yx)^{\lambda}] \nonumber \\
 & \on{Ad}[(x^{-1}y^{-1})^{2\mu}x^{2} g(x^{-3}y^{-1}x^{2},x^{-2}yx^{2})]
[(x^{-2}yx^{2})^{\lambda}] ... \nonumber \\
 & \on{Ad}[(x^{-1}y^{-1})^{(d-1)\mu}x^{d-1} g(x^{-d}y^{-1}x^{d-1},x^{1-d}yx^{d-1})]
[(x^{1-d}yx^{d-1})^{\lambda}] \nonumber \\
 &  (x^{-1}y^{-1})^{d\mu}g(x^{-1}y^{-1},y)
y^{{\lambda+1}\over 2}g(x,y) = y^{-d\mu}. 
\end{align}

Now for $t\in \{(x^{-1}y^{-1})^{d},x^{d},x^{\alpha}yx^{-\alpha}\}$ 
$(\alpha\in\ZZ)$, the morphisms $\kappa_{t}:\ZZ\to K_{2,d}^{0}$, $a\mapsto t^{a}$, 
extend to morphisms $\ZZ_{l}(N)\stackrel{\kappa_{(x^{-1}y^{-1})^{d}},
\kappa_{x^{d}}}{\to} K_{2,d}^{0}(\varphi_{N},l)$ and $\ZZ_{l}
\stackrel{\kappa_{x^{\alpha}yx^{-\alpha}}}{\to} K_{2,d}^{0}(\varphi_{N},l)$. 
For $t\in\{x^{-1}y^{-1},x,y\}$, the morphisms $\tilde\kappa_{t}:\ZZ\to F_{2}$, 
$a\mapsto t^{a}$ extend to morphisms $\ZZ_{l}(N')\stackrel{\kappa_{x^{-1}y^{-1}},
\kappa_{x}}{\to} F_{2}(\varphi_{N},l)$ and $\ZZ_{l}
\stackrel{\kappa_{y}}{\to}F_{2}(\varphi_{N'},l)$. The diagrams 
$\begin{matrix}\ZZ_{l} & \stackrel{\kappa_{y}}{\to}& 
K_{2,d}^{0}(\varphi_{N},l)\\
 & \scriptstyle{\tilde\kappa_{y}}\searrow &
 \downarrow\scriptstyle{\delta_{d}} \\
  & & F_{2}(\varphi_{N'},l)\end{matrix}$, 
$\begin{matrix} \ZZ_{l} & \stackrel{\kappa_{x^{\alpha}yx^{-\alpha}}}{\to}& 
K_{2,d}^{0}(\varphi_{N},l)\\
 & \scriptstyle{1}\searrow & \downarrow\scriptstyle{\delta_{d}}\\
  & & F_{2}(\varphi_{N'},l)\end{matrix}$ (for $\alpha\in\ZZ-N\ZZ$)
and $\begin{matrix} \ZZ_{l}(N)& \stackrel{\kappa_{(x^{-1}y^{-1})^{d}},
\kappa_{x^{d}}}{\to}& K_{2,d}^{0}(\varphi_{N},l)\\
\scriptstyle{[[-]]}\downarrow  & & \downarrow\scriptstyle{\delta_{d}}\\
\ZZ_{l}(N')  & \stackrel{\tilde\kappa_{x^{-1}y^{-1}},\tilde\kappa_{x}}{\to}
& F_{2}(\varphi_{N'},l)\end{matrix}$ commute. 

\begin{lemma} The automorphism $\tilde\kappa$ of $F_{2}$ (given by $\kappa(g)= 
g(x^{-1}y^{-1},y)$) restricts to an automorphism $\kappa$ of $K_{2,d}^{0}
\subset F_{2}$; $\tilde\kappa$ and $\kappa$ induce automorphisms 
$\tilde\kappa\in\on{Aut}(F_{2}(\varphi_{N'},l))$ and $\kappa\in 
\on{Aut}(K_{2,d}^{0}(\varphi_{N},l))$, and the diagram 
$\begin{matrix} K_{2,d}^{0}(\varphi_{N},l)& \stackrel{\kappa}{\to}& 
K_{2,d}^{0}(\varphi_{N},l)\\
\scriptstyle{\delta_{d}}\downarrow  && \downarrow \scriptstyle{\delta_{d}}\\
F_{2}(\varphi_{N'},l)  &\stackrel{\tilde\kappa}{\to}& F_{2}(\varphi_{N'},l)
\end{matrix}$ commutes. 
\end{lemma}

{\em Proof of Lemma.} We have $\varphi_{d}\circ \kappa = -\varphi_{d}$, 
which implies that $\kappa$ restricts to an automorphism $\tilde\kappa$
of $\on{Ker}\varphi_{d}=
K_{2,d}^{0}$. Since we also have  $\varphi_{N}\circ\kappa =-\varphi_{N}$, 
$\varphi_{N'}\circ \tilde\kappa = -\varphi_{N'}$, $\kappa$ and $\tilde\kappa$
also induce automorphisms of the completions $K_{2,d}^{0}(\varphi_{N},l)$
and $F_{2}(\varphi_{N'},l)$. It remains to show that the diagram 
$\begin{matrix} K_{2,d}^{0}& \stackrel{\kappa}{\to} & K_{2,d}^{0}\\
\scriptstyle{\delta_{d}}\downarrow & & \downarrow\scriptstyle{\delta_{d}}\\
F_{2}  &\stackrel{\tilde\kappa}{\to} &F_{2} \end{matrix}$ commutes. It
suffices to check this for the generators $x^{d}$, $x^{\alpha}yx^{-\alpha}$ ($\alpha
\in [0,d-1]$) of $K_{2,d}^{0}$. 

We have $\delta_{d}\circ \kappa(x^{d})=\delta_{d}((x^{-1}y^{-1})^{d})=
\delta_{d}(x^{-d}(x^{d-1}y^{-1}x^{1-d})...(xy^{-1}x^{-1})y^{-1})=
x^{-1}1...1y^{-1}=x^{-1}y^{-1}=\tilde\kappa(x)=\tilde\kappa\circ
\delta_{d}(x^{d})$; $\delta_{d}\circ\kappa(y)=\delta_{d}(y)=y=\tilde\kappa(y)
=\tilde\kappa\circ\delta_{d}(y)$; and for $\alpha\in[1,d-1]$, 
$\delta_{d}\circ\kappa(x^{\alpha}yx^{-\alpha})=\delta_{d}
[(x^{-1}y^{-1})^{\alpha}y(x^{-1}y^{-1})^{-\alpha}]=
\delta_{d}[\on{Ad}[y(-1)^{-1}...y(-\alpha)^{-1}][y(-\alpha)]]=1=
\tilde\kappa(1)=\tilde\kappa\circ\delta_{d}(x^{\alpha}yx^{-\alpha})$, where
$y(\alpha) = x^{\alpha}yx^{-\alpha}$. \hfill \qed\medskip 

{\em End of proof of Proposition.} Identity (\ref{pre:octogon}) holds in 
$K_{2,d}(\varphi_{N},l)$ and may be rewritten 
\begin{align*}
& g^{-1}\kappa_{y}({{\lambda-1}\over 2})\kappa(g)
\on{Ad}[(x^{-1}y^{-1})^{\mu}g(x^{-1}y^{-1},y)^{-1}x][\kappa_{x^{-1}yx}(\lambda)]
...
\on{Ad}[(x^{-1}y^{-1})^{d\mu}g(x^{-1}y^{-1},y)^{-1}x^{d}]
\\
 & [\kappa_{x^{1-d}yx^{d-1}}(\lambda)]
\kappa_{(x^{-1}y^{-1})^{d}}(\mu)\kappa(g)\kappa_{y}({{\lambda+1}\over 2})g=
\kappa_{x^{d}}(-\mu). 
\end{align*}
Applying $\delta_{d}$ to this identity and using the above commutation relations, 
we get (with $h:= \delta_{NN'}(g)$)
$$
h^{-1}\tilde\kappa_{y}({{\lambda -1}\over 2})\tilde\kappa(h)
\tilde\kappa_{x^{-1}y^{-1}}([[\mu]])\tilde\kappa(h)
\tilde\kappa_{y}({{\lambda+1}\over 2})h=\tilde\kappa_{x}(-[[\mu]]), 
$$
i.e., $([[\mu]],h)$ satisfies the octogon relation. 

Let us now show that $\delta_{NN'}(\ul{g}_{1}\ul{g}_{2})=
\delta_{NN'}(\ul{g}_{1})\delta_{NN'}(\ul{g}_{2})$ for 
$\ul{g}_{i}=(\lambda_{i},\mu_{i},f_{i},g_{i})\in\ul{\on{GTM}}(N)_{l}$, 
such that $\bar{\bar\mu}_{2}\in(\ZZ/d\ZZ)^{\times}$ (where $\bar{\bar\mu}$
is the image of $\bar\mu\in\ZZ/N\ZZ$ under $\ZZ/N\ZZ\to \ZZ/d\ZZ$, $\bar 1\mapsto 
\bar 1$). For this, it suffices to show that 
$\delta_{NN'}(\theta_{(\lambda_{2},\mu_{2},g_{2})}(g_{1})) = 
\theta_{(\lambda_{2},[[\mu_{2}]],\delta_{NN'}(g_{2}))}(\delta_{NN'}(g_{1}))$. 
We will show that we have commutative diagrams  
\begin{equation} \label{double:square}
\begin{matrix}
F_{2}(\varphi_{N},l) & \hookleftarrow & K_{2,d}^{0}(\varphi_{N},l) & 
\stackrel{\delta_{d}}{\to}& F_{2}(\varphi_{N'},l)\\
\scriptstyle{\theta_{(\lambda,\mu,g)}}{\downarrow} & & \downarrow& & 
\downarrow\scriptstyle{\theta_{(\lambda,[[\mu]],\delta_{NN'}(g))}}\\
F_{2}(\varphi_{N},l) & \hookleftarrow & K_{2,d}^{0}(\varphi_{N},l) & 
\stackrel{\delta_{d}}{\to}& F_{2}(\varphi_{N'},l)
\end{matrix}
\end{equation}
for any $\lambda\in\ZZ_{l}$, $\mu\in\ZZ_{l}(N)$ such that $\bar{\bar\mu}
\in (\ZZ/d\ZZ)^{\times}$ and  $g\in (\on{Ker}\varphi_{N})_{l}$; as $g_{1}\in 
(\on{Ker}\varphi_{N})_{l}\subset K_{2,d}^{0}(\varphi_{N},l)$, this implies the 
wanted equality. 

If $F_{2}(\varphi_{N},l)\stackrel{(\varphi_{d})_{l}}{\to}\ZZ/d\ZZ$ is the completion 
of $\varphi_{d}$, we have a commutative diagram $\begin{matrix}
F_{2}(\varphi_{N},l) & \stackrel{(\varphi_{d})_{l}}{\to} & \ZZ/d\ZZ \\
\scriptstyle{\theta_{(\lambda,f,g)}}\downarrow& \nearrow
\scriptstyle{\bar{\bar\mu}}(\varphi_{d})_{l}& \\
F_{2}(\varphi_{N},l) & & \end{matrix}$; as $K_{2,d}^{0}(\varphi_{N},l)
=\on{Ker}[F_{2}(\varphi_{N},l)\stackrel{(\varphi_{N})_{l}}{\to}\ZZ/d\ZZ]$, 
this implies that $\theta_{(\lambda,\mu,g)}$ restricts to an automorphism of 
$K_{2,d}^{0}(\varphi_{N},l)$, i.e., the left square of (\ref{double:square}) commutes. 

Let us now show that the right square commutes. As $K_{2,d}^{0}(\varphi_{N},l)$
is topologically generated by $x^{d}$ and the $x^{\alpha}yx^{-\alpha}$, $\alpha\in
[0,d-1]$, it suffices to check this commutativity on these generators. Recall that 
$\theta_{(\lambda,\mu,g)}$ is defined by $x\mapsto g(x,y)x^{\mu}g(x,y)^{-1}$, 
$y\mapsto y^{\lambda}$. Set $h:= \delta_{NN'}(g)$. Then 
$\theta_{(\lambda,[[\mu]],h)}\circ\delta_{d}(x^{d})=\theta_{(\lambda,[[\mu]],h)}(x)
=h(x,y)x^{[[\mu]]}h(x,y)^{-1}$, while $\delta_{d}\circ \theta_{(\lambda,\mu,g)}
(x^{d}) = \delta_{d}(g(x,y)x^{d\mu}g(x,y)^{-1}) = 
\delta_{d}(g(x,y))\delta_{d}(x^{d\mu})\delta_{d}(g(x,y)^{-1})$
(as $x^{d\mu}$ and $g(x,y)$ both belong to $K^{0}_{2,d}(\varphi_{N},l)$)
so $\delta_{d}\circ \theta_{(\lambda,\mu,g)}(x^{d})=h(x,y)x^{[[\mu]]}h(x,y)^{-1}
= \theta_{(\lambda,[[\mu]],h)}\circ \delta_{d}(x^{d})$. 

We have $\theta_{(\lambda,[[\mu]],h)}\circ \delta_{d}(y) = 
\theta_{(\lambda,[[\mu]],h)}(y) = y^{\lambda}$, while 
$\delta_{d}\circ\theta_{(\lambda,\mu,g)}(y) = \delta_{d}(y^{\lambda})
= y^{\lambda}$; if now $\alpha\in [1,d-1]$, then $\theta_{(\lambda,[[\mu]],h)}
\circ\delta_{d}(x^{\alpha}yx^{-\alpha}) = 1$ while 
$\delta_{d}\circ \theta_{(\lambda,\mu,g)}(x^{\alpha}yx^{-\alpha}) = 
\delta_{d}[\on{Ad}\{g(x,y)x^{\alpha\mu}g(x,y)^{-1}\}\{y^{\lambda}\}]
 = \delta_{d}[\on{Ad}\{g(x,y)g(x,x^{\alpha\mu}y x^{-\alpha\mu})^{-1}\}
 \{(x^{\alpha\mu}yx^{-\alpha\mu})^{\lambda}\}] = \on{Ad}\{\delta_{d}
 [g(x,y)g(x,x^{\alpha\mu}y x^{-\alpha\mu})^{-1}]\}
 \{\delta_{d}  [(x^{\alpha\mu}yx^{-\alpha\mu})^{\lambda}]\}$
 as $g(x,y)g(x,x^{\alpha\mu}y x^{-\alpha\mu})^{-1}\in K_{2,d}^{0}(\varphi_{N},l)$; 
 now $\delta_{d}(x^{\alpha\mu}yx^{-\alpha\mu})=1$ as $\bar\alpha
 \bar{\bar\mu}\neq 0$ in $\ZZ/d\ZZ$, so $\delta_{d}\circ \theta_{(\lambda,
 [[\mu]],h)}(x^{\alpha}yx^{-\alpha})=1=\theta_{(\lambda,\mu,g)}\circ 
 \delta_{d}(x^{\alpha}yx^{-\alpha})$. This proves the commutativity of 
 (\ref{double:square}), as wanted. \hfill \qed\medskip 

\subsection{The subgroup $\on{GTMD}(N)_{l}\subset \on{GTM}(N)_{l}$}

Let $\ul y:= (1,1,1,y)\in \ul{\on{GTM}}(N)_{l}$; this is the image of  
$1\in \ZZ\subset \ZZ\rtimes \ZZ/2\ZZ = \ul{\on{GTM}}$. The group morphism 
$\ZZ\to F_{2}$, $s\mapsto y^{s}$ extends to a morphism $\ZZ_{l}\to 
F_{2}(\varphi_{N},l)$, denoted again $s\mapsto y^{s}$; for $s\in\ZZ_{l}$, 
we set $\ul y^{s}:= (1,1,1,y^{s})\in \ul{\on{GTM}}(N)_{l}$. 

For $\ul g = (\lambda,\mu,f,g)\in \ul{\on{GTM}}(N)_{l}$, set $\chi_{l}(\ul g)
:= \lambda\in\ZZ_{l}$. Then for any $\ul g\in  \ul{\on{GTM}}(N)_{l}$ 
and any $s\in \ZZ_{l}$, we have $\ul y^{s}*\ul g = \ul g * \ul y^{\chi_{l}(\ul g)s}$, 
where $*$ is the product on $\ul{\on{GTM}}(N)_{l}$. It follows that any cocycle 
$\rho : \on{GTM}(N)_{l}\to \ZZ_{l}$ (i.e., $\rho(\ul g_{1} * \ul g_{2}) = 
\rho(\ul g_{1})\chi_{l}(\ul g_{2}) + \rho(\ul g_{2})$) gives rise to an automorphism 
$a_{\rho}$ of $\on{GTM}(N)_{l}$, given by $a_{\rho}(\ul g):= 
\ul g * \ul y^{\rho(\ul g)}$, i.e., $a_{\rho}(\lambda,\mu,f,g) = 
(\lambda,\mu,f,y^{\rho(\ul g)}g)$. 

\begin{lemma}
There is a unique cocycle $c_{NN'} : \on{GTM}(N)_{l}\to \ZZ_{l}$, defined 
by $c_{NN'}(\lambda,\mu,f,g):= c'(\delta_{NN'}(g)) - c'(\pi_{NN'}(g)) 
= c'(\delta_{NN'}(g))
- c(g)$, where $c:F_{2}(\varphi_{N},l)\to \ZZ_{l}$ and $c' : 
F_{2}(\varphi_{N'},l)\to \ZZ_{l}$ are the morphisms induced by 
$x\mapsto 1$, $y\mapsto 0$. 
\end{lemma} 

{\em Proof.} We will prove that each of the maps $(\lambda,\mu,f,g)\mapsto 
c(g), c'(\delta_{NN'}(g))$ is a cocycle. If $(\lambda_{1},\mu_{1},f_{1},g_{1})
*(\lambda_{2},\mu_{2},f_{2},g_{2}) = (\lambda,\mu,f,g)$, then 
$y^{c(g)} = g(1,y) = g_{1}(1,y^{\lambda_{2}})g_{2}(1,y) = y^{\lambda_{2}c(g_{1})
+c(g_{2})}$, so $c(\ul g)=c(\ul g_{1})\chi(\ul g_{2}) + c(\ul g_{2})$. 

On the other hand, $g(x,y) = g_{1}(g_{2}(x,y)x^{\mu_{2}}g_{2}(x,y)^{-1},
y^{\lambda_{2}})g_{2}(x,y)$;  the two factors of this product belong to 
$\on{Ker}(F_{2}(\varphi_{N},l)\stackrel{\varphi_{d}}{\to}\ZZ/d\ZZ)$, 
and $c'\circ \delta_{NN'} : \on{Ker}\varphi_{d}\to \ZZ_{l}$ is a group 
morphism, so $c'\circ \delta_{NN'}(g) = c'\circ\delta_{NN'}[g_{1}
(g_{2}(x,y)x^{\mu_{2}}g_{2}(x,y)^{-1},y^{\lambda_{2}})] + c'\circ
\delta_{NN'}(g_{2})$, so it remains to show that  
$c'\circ\delta_{NN'}[g_{1}
(g_{2}(x,y)x^{\mu_{2}}g_{2}(x,y)^{-1},y^{\lambda_{2}})]  = \lambda_{2}
(c'\circ\delta_{NN'})(g_{1})$. If $\tilde g_{1}\in (F_{d+1})_{l}$ is defined 
by $g_{1}(x,y) = \tilde g_{1}(x^{d}|y,xyx^{-1},...,x^{d-1}yx^{1-d})$, then 
$\tilde g_{1}(1|y,1,...,1) = y^{c'\circ\delta_{NN'}(g_{1})}$. On the other hand, 
\begin{align*}
& g_{1}(g_{2}(x,y)x^{\mu_{2}}g_{2}(x,y)^{-1})
= \tilde g_{1}(g_{2}(x,y)x^{d\mu_{2}}g_{2}(x,y)^{-1}|
y^{\lambda_{2}},\on{Ad}[g_{2}(x,y)x^{\mu_{2}}g_{2}(x,y)^{-1}x^{-\mu_{2}}]
[x^{\mu_{2}}y^{\lambda_{2}}x^{-\mu_{2}}], 
\\
 & \on{Ad}[g_{2}(x,y)x^{2\mu_{2}}g_{2}(x,y)^{-1}x^{-2\mu_{2}}]
[x^{2\mu_{2}}y^{\lambda_{2}}x^{-2\mu_{2}}],...)
\stackrel{\delta_{NN'}}{\mapsto} \tilde g_{1}(1|y^{\lambda_{2}},1,...,1)
= y^{\lambda_{2}(c'\circ \delta_{NN'})(g_{1})}, 
\end{align*}
as $x^{d}\stackrel{\delta_{NN'}}{\mapsto}1$, $y\stackrel{\delta_{NN'}}
{\mapsto}y$ and $x^{\alpha}yx^{-\alpha}\stackrel{\delta_{NN'}}{\mapsto}1$
for $\alpha\in [1,d-1]$ as $\bar{\bar\mu_{2}}\in (\ZZ/d\ZZ)^{\times}$. 
So $c'\circ\delta_{NN'}[g_{1}
(g_{2}(x,y)x^{\mu_{2}}g_{2}(x,y)^{-1},y^{\lambda_{2}})]  = \lambda_{2}
(c'\circ\delta_{NN'})(g_{1})$, as wanted. 
\hfill \qed \medskip 
 
 We denote the resulting automorphism of $\on{GTM}(N)_{l}$ by $a_{NN'}$, 
 Since $\pi_{NN'},\delta_{NN'} : \on{GTM}(N)_{l}\to \on{GTM}(N')_{l}$
 are group morphisms, 
 \begin{align*}
 & \on{GTMD}(N)_{l}:= \{\ul g\in \on{GTM}(N)_{l} | \forall N'|N, 
 \delta_{NN'}(\ul g) = \pi_{NN'}\circ a_{NN'}(\ul g)\}
 \\
  & = \{(\lambda,\mu,f,g)\in \on{GTM}(N)_{l} | \forall N'|N, 
 \delta_{NN'}(g) = y^{c'(\delta_{NN'}(g)) - c(g)}\pi_{NN'}(g)\}
 \end{align*}
is a subgroup of $\on{GTM}(N)_{l}$. 

\begin{proposition}
The morphism $\wh{\on{GTM}}\to \on{GTM}(N)_{l}$ factors through 
$\wh{\on{GTM}} \to \on{GTMD}(N)_{l} \hookrightarrow 
\on{GTM}(N)_{l}$.  
\end{proposition}

{\em Proof.} Recall the cocycle $c_{N}$ of $\wh{\on{GTM}}$
constructed in Section \ref{sec:9:1}, then the diagram 
$\begin{matrix} \wh{\on{GTM}} & \stackrel{c_{d}}{\to}& \wh\ZZ\\
\downarrow && \downarrow \\
 \on{GTM}(N)_{l} &\stackrel{c_{NN'}}
{\to}& \ZZ_{l} \end{matrix}$ commutes. Recall that $c_{d}$ gives rise to 
an automorphism $\hat a_{d}\in \on {Aut}(\wh{\on{GTM}})$, 
$(\lambda,\mu,f,g)\mapsto (\lambda,\mu,f,y^{\hat c_{d}(f)}g)$. We 
have commutative diagrams 
$$
\begin{matrix} \wh{\on{GTM}}& \stackrel{a_{d}}{\to} & \wh{\on{GTM}}\\
\downarrow & & \downarrow\\
\on{GTM}(N)_{l}  & \stackrel{a_{NN'}}{\to}&\on{GTM}(N)_{l} \end{matrix},
\quad  
\begin{matrix} \wh{\on{GTM}}& & \\
\downarrow  & \searrow & \\
\on{GTM}(N)_{l}  & \stackrel{\pi_{NN'}}{\to}& \on{GTM}(N')_{l}\end{matrix}
\quad \on{and} \quad 
\begin{matrix} \wh{\on{GTM}}& \stackrel{\delta_{d}}{\to}& \wh{\on{GTM}}\\
\downarrow & & \downarrow \\
 \on{GTM}(N)_{l} & \stackrel{\delta_{NN'}}{\to}& \on{GTM}(N')_{l}\end{matrix}.
 $$ 
We have seen (Theorem \ref{thm:ihara}) that $\delta_{d}=a_{d}\in 
\on{Aut}(\wh{\on{GTM}})$; 
it follows that for any $\ul g\in \on{Im}(\wh{\on{GTM}}\to \on{GTM}(N)_{l})$, we have 
$\delta_{NN'}(\ul g)=\pi_{NN'}\circ a_{NN'}(\ul g)$. \hfill \qed \medskip   

Let $\on{Div}(N)$ be the partial semigroup (for the multiplication) of all the divisors of $N$. 

\begin{proposition}
1) We have a group morphism $\on{GTMD}(N)_{l}\to \on{Hom}(\on{Div}(N),\ZZ_{l})
\simeq \ZZ_{l}^{\nu}$ (here $\nu$ is the number of distinct prime powers in the decomposition 
of $N$), taking $(\lambda,\mu,f,g)$ to $[d\mapsto c_{NN'}(g)]$. 

2) The group morphisms $\delta_{NN'},\pi_{NN'} : \on{GTM}(N)_{l}\to 
\on{GTM}(N')_{l}$ restrict to morphisms $\on{GTMD}(N)_{l}\to 
\on{GTMD}(N')_{l}$. 
\end{proposition}

{\em Proof.} Both statements are consequences of the following equalities: if
$N''|N'|N$, then $\delta_{N'N''}\circ \delta_{NN'}=\delta_{NN''}$, 
$\pi_{N'N''}\circ \pi_{NN'}=\pi_{NN''}$, and if $d_{1},d_{2},d_{1}d_{2}
\in \on{Div}(N)$, $N_{i}=N/d_{i}$ ($i=1,2$), $N_{3}=N/(d_{1}d_{2})$, then 
$\pi_{N_{1}N_{3}}\circ \delta_{NN_{1}} = \delta_{N_{2}N_{3}}\circ 
\pi_{NN_{2}}$. \hfill \qed \medskip 

\subsection{The subgroup $\on{GTMD}(N,\kk)\subset \on{GTM}(N,\kk)$}

One defines as above morphisms $\pi_{NN'},\delta_{NN'}:\on{GTM}(N,\kk)\to 
\on{GTM}(N',\kk)$, a cocycle $c_{NN'}:\on{GTM}(N,\kk)\to \kk$ and the 
subgroup $\on{GTMD}(N,\kk)\subset \on{GTM}(N,\kk)$. We have as above a 
group morphism $\on{GTMD}(N,\kk)\to \on{Hom}(\on{Div}(N),\kk)\simeq \kk^{\nu}$. 
When $\kk = \QQ_{l}$, these morphisms are compatible with the morphisms $\on{GTM}(N)_{l}
\hookrightarrow \on{GTM}(N,\QQ_{l})$, in particular 
$\on{GTMD}(N)_{l}\hookrightarrow \on{GTMD}(N,\QQ_{l})$. 

\section{Distribution relations and the torsor $\Psdist(N,\kk)$} 
\label{sect:8}

Recall that $N,N'$ are integers $\geq 1$ with $N' | N$; we set $d:= N/N'$. 

\subsection{The morphisms
$\pi_{NN'}:\ul\Pseudo(N,\kk) \to \ul\Pseudo(N',\kk)$, $\on{GRTM}(N,\kk) \to 
\on{GRTM}(N',\kk)$}

\subsubsection{} 

The inclusion $K_{n,N}\stackrel{\pi_{NN'}}{\hookrightarrow} K_{n,N'}
(\subset P_{n+1})$ may be viewed as associated to $W_{n,N}\to W_{n,N'}$, 
$(w_{1},...,w_{n})\mapsto (w_{1}^{d},...,w_{n}^{d})$. This group morphism 
induces a Lie algebra morphism $\on{Lie}K_{n,N}\to 
\on{Lie}K_{n,N'}$, given by 
$\Xi_{0i} \mapsto d \Xi'_{0i}$, $\xi_{ij}(\alpha) \mapsto (q \Xi'_{0i}) 
* \xi'_{ij}(\alpha') * (-q\Xi'_{0i})$ for $\alpha=N'q+\alpha'$, $\alpha'\in[0,N'-1]$
(the primes denote the generators of $\on{Lie}K_{n,N'}$). 
This morphism is not injective, e.g.,  $\xi_{ij}(\alpha+N')$ and 
$(\Xi_{0i}/d) * \xi_{ij}(\alpha) * (-\Xi_{0i}/d)$ have the same image. 

The associated graded Lie algebra morphism is 
$\pi_{NN'} : \t_{n,N} \to \t_{n,N'}$ given by 
$$
t^{0i}_0 \mapsto d t^{\prime 0i}_0, \; 
t^{ij}(a) \mapsto t^{\prime ij}(\bar a); 
$$
here $t^{\prime 0i}_0$, $t^{\prime ij}(a)$ are the generators of 
$\t_{n,N'}$ and $\bar a$ is the image of $a$ under $\ZZ/N\ZZ \to \ZZ/N'\ZZ$, 
$\bar 1 \mapsto \bar 1$. 

\subsubsection{The morphisms $\pi_{NN'} : \on{GRTM}(N,\kk) \to 
\on{GRTM}(N',\kk)$}

The morphisms $\pi_{NN'} : \t_{n,N} \to \t_{n,N'}$ 
are compatible with the insertion-coproduct maps. Recall that  
$\pi_{NN'} : \t^0_{3,N} \to \t^0_{3,N'}$ is defined by 
$$
A\mapsto d A', \; 
b(a) \mapsto b'(\bar a). 
$$
Then
$\pi_{NN'}$ induces a group morphism $\on{GRTM}_{(\bar 1,1)}(N,\kk) \to 
\on{GRTM}_{(\bar 1,1)}(N',\kk)$, by $(h,k)\mapsto (h,\pi_{NN'}(k))$; 
together with the morphism $(\ZZ/N\ZZ)^\times \times \kk^\times 
\stackrel{\on{can} \times \on{id}}{\to}
(\ZZ/N'\ZZ)^\times \times \kk^\times$, this 
extends to a morphism $\on{GRTM}(N,\kk) \to \on{GRTM}(N',\kk)$. 
In the same way, $\pi_{NN'}$ induces a Lie algebra morphism 
$\grtm_{(\bar 1,1)}(N,\kk) \to \grtm_{(\bar 1,1)}(N',\kk)$, graded and 
compatible with the 
actions of $(\ZZ/N\ZZ)^\times$ and $(\ZZ/N'\ZZ)^\times$.  It extends 
to a morphism $\grtm_{\bar 1}(N,\kk) \to \grtm_{\bar 1}(N',\kk)$, 
commuting with the projections $\grtm_{\bar 1}(N'',\kk)\to\kk$, 
$N'' = N,N'$. 

\subsubsection{The morphisms $\pi_{NN'} : 
\ul\Pseudo(N,\kk) \to \ul\Pseudo(N',\kk)$}

There is a unique morphism $\ul\Pseudo(N,\kk)
\to \ul\Pseudo(N',\kk)$, taking $(a,\lambda,\Phi,\Psi)$ to 
$(\bar a,\lambda,\Phi,\pi_{NN'}(\Psi))$. 
This morphism is compatible with the morphism 
$\on{GRTM}(N,\kk) \stackrel{\pi_{NN'}}{\to}
\on{GRTM}(N',\kk)$, and it restricts to a morphism of torsors 
$\Pseudo(N,\kk) \to \Pseudo(N',\kk)$. 

\subsection{The morphisms 
$\delta_{NN'} : \Pseudo(N,\kk) \to \Pseudo(N',\kk)$, 
$\on{GRTM}(N,\kk) \to \on{GRTM}(N',\kk)$}

If $a\in \ZZ/N\ZZ$, we write $d|a$ iff $a\in d\ZZ/N\ZZ$. We then denote by 
$a/d$ the unique element $a'\in \ZZ/N'\ZZ$ such that $da' = a$. 

\subsubsection{} 

Recall that the morphism $\delta_{NN'} : K_{n,N}\to K_{n,N'}$ 
may be viewed as induced by the morphism $W_{n,N}\to W_{n,N'}$, 
$(w_{1},...,w_{n})\mapsto (w_{1},...,w_{n})$. 
It is given by $X_{0i}\mapsto X'_{0i}$,  $x_{ij}(\alpha)\mapsto 
x'_{ij}(\alpha/d)$ if $d|\alpha$, $x_{ij}(\alpha)\mapsto 1$ if $d\nmid \alpha$. 

The associated graded morphism of $\on{Lie}\delta_{NN'}$
is $\delta_{NN'} : \t_{n,N} \to \t_{n,N'}$, given by 
$$
t_0^{0i} \mapsto t^{\prime 0i}_0, \quad 
t^{ij}(a) \mapsto 
t^{\prime ij}(a/d) \on{\ if\ }d|a, \; t^{ij}(a) \mapsto 0
\on{\ if\ }d\nmid a. 
$$

\subsubsection{The morphisms $\delta_{NN'} : \on{GRTM}(N,\kk) \to
\on{GRTM}(N',\kk)$}

The morphisms $\delta_{NN'} : \t_{n,N} \to \t_{n,N'}$ are compatible with
the insertion-coproduct morphisms. Recall that $\delta_{NN'} : \t^0_{3,N}
\to \t^0_{3,N'}$ is defined by 
$$
A\mapsto A', \quad b(a) \mapsto b'(a/d) \on{\ if\ }d|a, \; 
b(a) \mapsto 0 \on{\ if\ }d\nmid a. 
$$ 
One checks that there is a unique group morphism $\delta_{NN'} : 
\on{GRTM}_{(\bar 1,1)}(N,\kk) \to \on{GRTM}_{(\bar 1,1)}(N',\kk)$, 
taking $(h,k)$
to $(h,\delta_{NN'}(k))$. In the same way as above, its extends to a 
group morphism $\delta_{NN'} : 
\on{GRTM}(N,\kk) \to \on{GRTM}(N',\kk)$. The Lie algebra
versions of these morphisms are $\delta_{NN'} : \grtm_{\bar 1}(N,\kk) \to 
\grtm_{\bar 1}(N',\kk)$ and $\delta_{NN'} : \grtm_{(\bar 1,1)}(N,\kk) \to 
\grtm_{(\bar 1,1)}(N',\kk)$ with the same properties as above. 

\subsubsection{The morphisms $\delta_{NN'} : \Pseudo(N,\kk) \to 
\Pseudo(N',\kk)$} 

There is a unique morphism $\ul\Pseudo(N,\kk) \to \ul\Pseudo(N,\kk)$
taking $(a,\lambda,\Phi,\Psi)$ to $(\bar a,\lambda,\Phi,\delta_{NN'}(\Psi))$, 
and it shares the properties of $\pi_{NN'} : \ul\Pseudo(N,\kk) 
\to \ul\Pseudo(N,\kk)$ stated above.  

\subsection{The case $N'=1$}

\begin{lemma}
The compositions of $\grtm_{(\bar 1,1)}(N,\kk) \stackrel{\pi_{N1}}{\to} 
\grtm_{(\bar 1,1)}(1,\kk) \simeq 
\grt_1(\kk) \oplus \kk(0,t^{23})$ and $\grtm_{(\bar 1,1)}(N,\kk) 
\stackrel{\delta_{N1}}{\to} 
\grtm_{(\bar 1,1)}(1,\kk) \simeq 
\grt_1(\kk) \oplus \kk(0,t^{23})$ with projection on the first component 
both coincide with the canonical morphism $\grtm_{(\bar 1,1)}(N,\kk) \to \grt_1(\kk)$.  
\end{lemma}

{\em Proof.} Denote by $\psi\mapsto \overline\psi$ the morphism 
$\pi_{N1} : \t^0_{3,N}\to \t^0_3$ defined by 
$A\mapsto N t^{12}$, $b(a) \mapsto 
t^{23}$ and $\eps_{0}:\t_{n,N}\to \t_{n}$ the morphism $t_{0}^{0i}\mapsto 0$, 
$t(a)^{ij}=\delta_{a,\bar 0}t^{ij}$. Then $\eps_{0}(\overline\psi)$ is proportional to 
$t^{12}$, so it is zero if $\psi$ has degree  $>1$. 
On the other hand, (\ref{pent:grtm}) implies that $\varphi^{1,2,3} + \overline
\psi^{0,12,3}
+ \overline \psi^{0,1,2,} = \overline \psi^{01,2,3} + \overline\psi^{0,1,23}$.
Applying $\eps_{0}$ to this identity, we get $\overline\psi = \varphi$. 
The case of $\delta_{N1}$ is similar. \hfill \qed \medskip 

$\pi_{N1}$ is related to the canonical projection $\ul\Pseudo(N,\kk) \to
\ul\Assoc(\kk)$, $\Psi\mapsto \Phi$ as follows: there exists 
a map $\gamma : \ul\Pseudo(N,\kk) 
\to \kk$, such that $\pi_{N1}(\Psi) = e^{\gamma(\Psi) t^{23}}\Phi$. 
Explicitly, we have $\gamma(\Psi) =  N \gamma_A(\Psi) + \sum_{a\in \ZZ/N\ZZ}
\gamma_a(\Psi)$, where $\gamma_A(\Psi)$, $\gamma_a(\Psi)$ are the 
coefficients of $A$, $b(a)$ in $\Psi$.  

In the same way, if we set $\gamma'(\Psi) := \gamma_A(\Psi) + \gamma_0(\Psi)$,
then $\delta_{N1}(\Psi) = e^{\gamma'(\Psi) t^{23}}\Phi$.

\subsection{The scheme $\ul\Psdist(N,\kk)$ and the groups 
$\on{GTMD}(N,\kk)$ and $\on{GRTMD}(N,\kk)$} 

\subsubsection{The scheme $\ul\Psdist(N,\kk)$}

Define $\ul\Psdist(N,\kk) \subset \ul\Pseudo(N,\kk)$ as the subset of all
$(a,\lambda,\Phi,\Psi)$ such that  
$$
\forall N'|N, 
 \; 
\pi_{NN'}(\Psi) = e^{\rho_{NN'}(\Psi) b(0)} \delta_{NN'}(\Psi). 
$$
Here $\rho_{NN'}(\Psi) = \rho(\pi_{NN'}(\Psi)) - \rho(\Psi)$, where
$\rho(\tilde\Psi)$ is the coefficient of $b(0)$ in $\tilde\Psi\in 
\on{exp}(\hat\t_{2,N''}^{0})$ ($N''=N,N'$). 

\begin{lemma} \label{lemma:eliminate}
For any $(a,\lambda,\Psi,\Psi)\in \ul\Psdist(N,\kk)$, the map 
$\on{Div}(N)\ni d \mapsto \rho_{N,N/d}(\Psi)$ is in $\on{Hom}(\on{Div}(N),\kk)$; 
this defines a map $\ul{\Psdist}(N,\kk)\to \on{Hom}(\on{Div}(N),\kk)$. 
\end{lemma}

{\em Proof.} This again follows from the identities $\pi_{N'N''}\circ \pi_{NN'}
=\pi_{NN''}$, $\delta_{N'N''}\circ \delta_{NN'} = \delta_{NN''}$ and 
if $d_{1},d_{2},d_{3}:= d_{1}d_{2}\in \on{Div}(N)$, $N_{i}=N/d_{i}$, then 
$\pi_{N_{1}N_{3}}\circ \delta_{NN_{1}} = \delta_{N_{2}N_{3}}\circ 
\pi_{NN_{2}}$ (morphisms between Lie algebras $\hat\t_{2,\nu}^{0}$ instead of 
groups $(K_{2,\nu})_{l}$).  \hfill \qed \medskip 


We will define ${\Psdist}(N,\kk) :=  \ul\Psdist(N,\kk)\cap 
\Pseudo(N,\kk)$, so this is the subset of $\ul\Psdist(N,\kk)$ 
defined by $(a,\lambda) \in (\ZZ/N\ZZ)^\times \times \kk^\times$. 

\subsubsection{The group $\on{GRTMD}(N,\kk)$} 

We have morphisms $\pi_{NN'},\delta_{NN'} : \on{GRTM}_{(\bar 1,1)}(N,\kk)
\to \on{GRTM}_{(\bar 1,1)}(N',\kk)$, given by $(h,k)\mapsto (h,\pi_{NN'}(k))$
and $(h,k)\mapsto (h,\delta_{NN'}(k))$. We also have $a_{NN'}\in 
\on{Aut}(\on{GRTM}_{(\bar 1,1)}(N,\kk))$, defined by $a_{NN'}(h,k):= 
(h,e^{\rho_{NN'}(k)b(0)}k)$. 

We then set 
\begin{align*}
\on{GRTMD}_{(\bar 1,1)}(N,\kk) & := \{(h,k)\in \on{GRTM}_{(\bar 1,1)}(N,\kk)|
\forall N'|N, \pi_{NN'}(h,k) = \delta_{NN'}\circ a_{NN'}(h,k)\} \\
& = \{(h,k)\in \on{GRTM}_{(\bar 1,1)}(N,\kk) | \pi_{NN'}(k) = e^{\rho_{NN'}(k)b(0)}
\delta_{NN'}(k)\}. 
\end{align*}
Then $\on{GRTMD}_{(\bar 1,1)}(N,\kk)$ is a subgroup of 
$\on{GRTM}_{(\bar 1,1)}(N,\kk)$, which is preserved by the action of 
$(\ZZ/N\ZZ)^{\times}\times \kk^{\times}$. Set 
$\on{GRTMD}(N,\kk):= \on{GRTMD}_{(\bar 1,1)}(N,\kk)\rtimes 
((\ZZ/N\ZZ)^{\times}\times\kk^{\times})$. 

\begin{lemma}
If $(h,k)\in \on{GRTMD}_{(\bar 1,1)}(N,\kk)$, then $\on{Div}(N)\ni d\mapsto 
\rho_{N,N/d}(k)$ belongs to $\on{Hom}(\on{Div}(N),\kk)\simeq \kk^{\nu}$. 
This defines a group morphism $\on{GRTMD}_{(\bar 1,1)}(N,\kk)\to\kk^{\nu}$. 
\end{lemma} 

{\em Proof.} The proof is the same as above; the group morphism property is checked
directly. \hfill \qed \medskip 

This group morphism is compatible with the projection $(\ZZ/N\ZZ)^{\times}
\times \kk^{\times}\to \kk^{\times}$ and the scaling action of $\kk^{\times}$
on $\kk^{\nu}$. We thus get a group morphism $\on{GRTMD}(N,\kk)\to 
(\ZZ/N\ZZ)^{\times} \times (\kk^{\nu}\rtimes \kk^{\times})$.



Denote by $\grtmd_{(\bar 1,1)}(N,\kk) \subset \grtmd_{\bar 1}(N,\kk)$ 
the Lie algebras of $\on{GRTMD}_{(\bar 1,1)}(N,\kk)\subset \on{GRTMD}(N,\kk)$. 
Let $\grtm_{\bar 1}(N,\kk) = \kk\oplus \grtm_{(\bar 1,1)}(N,\kk)_1 
\oplus \grtm_{(\bar 1,1)}(N,\kk)_2 \oplus \cdots$ be the decomposition of 
$\grtm_{\bar 1}(N,\kk)$ in homogeneous 
components. Then the decomposition of $\grtmd_{\bar 1}(N,\kk)$ is 
$ \kk\oplus \grtmd_{(\bar 1,1)}(N,\kk)_1 \oplus \grtmd_{(\bar 1,1)}(N,\kk)_2 
\oplus\cdots$, where 
$$
\grtmd_{(\bar 1,1)}(N,\kk)_1 = \{x\in \grtm_{(\bar 1,1)}(N,\kk)_1 
|  \forall N'|N, (\pi_{NN'}-\delta_{NN'})(x) = \rho_{NN'}(x)b(0)\}
$$
and for $n>1$, 
$$
\grtmd_{(\bar 1,1)}(N,\kk)_n = \{x\in \grtm_{(\bar 1,1)}(N,\kk)_n | 
\forall N'|N, 
\pi_{NN'}(x) = \delta_{NN'}(x)\}. 
$$  
The element $1\in \kk$ is such that $[1,\psi] = n\psi$ for $\psi$ of degree 
$n$, and $\grtmd_{(\bar 1,1)}(N,\kk)$ is the sum of homogeneous components of 
$\grtm_{\bar 1}(N,\kk)$ of degree $>0$. The action of 
$(\ZZ/N\ZZ)^{\times}$ on $\grtm_{(\bar 1,1)}(N,\kk)$ restricts to a 
degree-preserving action on $\grtmd_{(\bar 1,1)}(N,\kk)$. 
The Lie algebra morphism $\grtmd_{(\bar 1,1)}(N,\kk)\to \kk^{\nu}$
is the direct sum of $\grtmd_{(\bar 1,1)}(N,\kk)_{1}\to \kk^{\nu}$, 
$x\mapsto [d\mapsto \rho_{N,N/d}(x)]$, and zero on the other components; 
this morphism is invariant under the action of $(\ZZ/N\ZZ)^{\times}$. 

\subsubsection{Torsor structure of $\Psdist(N,\kk)$} \label{sec:torsor}

We have commutative diagrams 
$$\begin{matrix} \on{GTM}(N,\kk)\times \ul\Pseudo(N,\kk)\times\on{GRTM}(N,\kk)
& \to  & \ul\Pseudo(N,\kk)\\
\begin{matrix} \scriptstyle{\pi_{NN'}\times\pi_{NN'}\times\pi_{NN'}}\\
\scriptstyle{\on{(resp.,\ }
\delta_{NN'}\times\delta_{NN'}\times\delta_{NN'})}\end{matrix}
\downarrow & & \downarrow\scriptstyle{\pi_{NN'} \on{\ (resp.,\ } \delta_{NN'})}\\
\on{GTM}(N',\kk)\times \ul\Pseudo(N',\kk)\times\on{GRTM}(N',\kk)
& \to  & \ul\Pseudo(N',\kk) \end{matrix}$$
and 
$$\begin{matrix} \on{GTM}(N,\kk)\times \ul\Pseudo(N,\kk)\times\on{GRTM}(N,\kk)
& \to  & \ul\Pseudo(N,\kk)\\
\scriptstyle{a_{NN'}}\downarrow & 
& \downarrow\scriptstyle{a_{NN'}}\\
\on{GTM}(N,\kk)\times \ul\Pseudo(N,\kk)\times\on{GRTM}(N,\kk)
& \to  & \ul\Pseudo(N,\kk) \end{matrix}$$
which implies that the action of $\on{GRTM}(N,\kk)$ on $\ul{\Pseudo}(N,\kk)$ restricts 
to an action of $\on{GRTMD}(N,\kk)$ on $\ul{\Psdist}(N,\kk)$. 
One checks that that $\Psdist(N,\kk)$
is a torsor under the actions of $\on{GTMD}(N,\kk)$ and $\on{GRTMD}(N,\kk)$.  
The right action of $\on{GRTMD}(N,\kk)$ on $\ul{\Psdist}(N,\kk)$ maps to the 
right action of $(\ZZ/N\ZZ)^{\times}\times (\kk^{\nu}\rtimes\kk^{\times})$
on $(\ZZ/N\ZZ)\times (\kk^{\nu}\rtimes\kk)$. 


\subsection{Distributivity of $\Psi_{\on{KZ}}$} 

Let us denote by $\Psi^N_{\on{KZ}}$ the element of 
$\Pseudo_{(-\bar 1,2\pi\i)}(N,\CC)$ constructed in Section \ref{sect:PsiKZ}. 

\begin{lemma}
The morphism $\pi_{NN'} : \Pseudo_{(-\bar 1,2\pi\i)}(N,\CC)
\to \Pseudo_{(-\bar 1,2\pi\i)}(N',\CC)$ takes $(\Phi_{\on{KZ}},
\Psi^{N}_{\on{KZ}})$ to $(\Phi_{\on{KZ}},
e^{\log(d) b(0)}\Psi^{N'}_{\on{KZ}})$. 
Moreover $\delta_{NN'}(\Psi_{\on{KZ}}^N) = \Psi_{\on{KZ}}^{N'}$, 
so $(\Phi_{\on{KZ}},\Psi_{\on{KZ}}^N) \in \Psdist_{(-\bar 1,2\pi i)}(N,\CC)$; 
its image in $\on{Hom}(\on{Div}(N),\CC)$ is the restriction of $\log$  
to $\on{Div}(N)$. 
\end{lemma}

{\em Proof.} Let $H(z)$ be a solution of 
$$
{{\on{d} H}\over {\on{d}z}} = \big( {A\over z} + \sum_{a=0}^{N-1}
{{b(a)}\over{z-\zeta_N^a}}\big) H(z), 
$$
where $\zeta_N$ is a primitive $N$th root of $1$. 
Let $\overline H(z)$ be the image of $H(z)$ by $\pi_{NN'}$. 
Then it follows from the identity $\sum_{k=0}^{d-1} 1/(z-\zeta_d^k) = 
d z^{d-1}/(z^d-1)$ that $K(w) := \overline H(w^{1/d})$ satisfies 
$$
{{\on{d}\overline K}\over{\on{d}w}} =
\big( {{A'}\over{w}} 
+ \sum_{a'=0}^{N'-1} {{b'(a')}\over{w - \zeta_{N'}^{a'}}}\big)
\overline K(w) 
$$
(we set $\zeta_d := \zeta_N^{N'}$, $\zeta_{N'} := \zeta_N^d$).
The renormalized holonomy of $K(w)$ is $\Psi_{\on{KZ}}^{N'}$, 
hence the renormalized holonomy of $\overline H(z)$ is 
$e^{\on{log}(d) b(0)}\Psi_{\on{KZ}}^{N'}$ (taking into account the 
change of variables). On the other hand, this 
is the image of $\Psi^N_{\on{KZ}}$ by $\pi_{NN'}$.   
\hfill \qed \medskip 

\subsection{Surjectivity of $\Psdist(N,\kk) \to (\ZZ/N\ZZ)^\times 
\times (\kk^\nu \rtimes \kk^\times)$} \label{sect:rat}

One can use the torsor structure of $\Psdist(N,\kk)$ and 
its nonemptiness for $\kk = \CC$ to prove that 
$\Psdist(N,\kk)\to (\ZZ/N\ZZ)^{\times} \times \kk^{\times}$
is surjective. We will see in Theorem \ref{thm:estimate} that the morphism 
$\grtmd_{(\bar 1,1)}(N,\kk)\to \kk^{\nu}$ is surjective, which 
implies: 

\begin{proposition}
The map $\Psdist(N,\kk) \to (\ZZ/N\ZZ)^\times \times 
(\kk^{\nu}\rtimes \kk^\times)$ is surjective. 
\end{proposition}

\section{Generators of $\grtmd_{(\bar 1,1)}(N,\CC)$} 
\label{sect:drinfeld} \label{sect:9}

\subsection{} \label{sect:dr:conj}
The element $(\lambda=-1,f=1)$ of $\on{GT}$
takes a QTQBA $(A,m_{A},\Delta_{A},R,\Phi_{A})$ to 
$A' := (A,m_{A},\Delta_{A},(R_{A}^{21})^{-1},\Phi_{A})$; the element 
$(\lambda = -1,f=g=1)$ of $\on{GTM}$ takes the QRA
$(B,m_{B},\Delta_{B},E,\Psi_{B})$ over $A$ to the QRA 
$B' := (B,m_{B},\Delta_{B},E^{-1},\Psi_{B})$ over $A'$. 

It follows that there is a unique map $\ul{\Assoc}_\lambda(\kk)
\to \ul{\Assoc}_{-\lambda}(\kk)$ taking $\Phi$ to itself, and a 
unique morphism $\ul{\Psdist}_{(\alpha,\lambda)}(N,\kk) \to 
\ul{\Psdist}_{(-\alpha,-\lambda)}(N,\kk)$, taking 
$(\Phi,\Psi)$ to itself. Composing the first map with the 
action of $-1\in \kk^\times
\subset \on{GRT}(\kk)$, we get a permutation of 
$\ul{\Assoc}_\lambda(\kk)$, taking $\Phi(A,B)$ to $\Phi(-A,-B)$. 
Composing the second map with the action of $(-\bar 1,-1)\in 
(\ZZ/N\ZZ)^\times \times \kk^\times \subset 
\on{GRTMD}(N,\kk)$, we get a permutation of 
$\ul{\Psdist}_{(\alpha,\lambda)}(N,\kk)$, taking 
$$
\big(\Phi(U,V),\Psi(A | b(0),\ldots,b(N-1))\big) 
\; \on{to} \; 
\big(\Phi(-U,-V),\Psi(-A | -b(0),\ldots,-b(1-N))\big) . 
$$ 
In particular, there is a unique element $(g_{\on{KZ}}, 
h_{\on{KZ}}) \in \on{GRTMD}_{(\bar 1,1)}(N,\CC)$, such that 
$$
\big(\Phi_{\on{KZ}}(-U,-V),
\Psi_{\on{KZ}}(-A | -b(0),\ldots,-b(1-N))\big) 
= (\Phi_{\on{KZ}},\Psi_{\on{KZ}}) * (g_{\on{KZ}},h_{\on{KZ}}) . 
$$ 
Let $\on{Log} : \on{GRT}_1(\CC) \to \grt_1(\CC)$ be the logarithmic map. Set
$\varphi := \on{Log}(g_{\on{KZ}})$. According to \cite{Dr2}, for $n$ odd $\geq 3$,
the degree $n$ component $\varphi_n$ of $\varphi$ is a generator of 
$\grt_1(\CC)$ (i.e., its  class in 
$\grt_1(\CC)/[\grt_1(\CC),\grt_1(\CC)]$ is nonzero). 

Set $(\varphi,\psi) := \on{Log}(g_{\on{KZ}},h_{\on{KZ}})$, so $(\varphi,\psi)
\in \grtmd_{(\bar 1,1)}(N,\CC)$. When $n$ is odd $\geq 3$, the degree $n$ component 
$(\varphi_n,\psi_n)$ of $(\varphi,\psi)$ is a generator of 
$\grtmd_{(\bar 1,1)}(N,\CC)$, and is a preimage of $\varphi_n$. So the image of 
$\grtmd_{(\bar 1,1)}(N,\CC) \to \grt_1(\CC)$
contains the subalgebra generated by the $\varphi_n$. 

\subsection{} To obtain more information of $(\varphi_n,\psi_n) \in 
\grtmd_{(\bar 1,1)}(N,\CC)$, we 
study the expansion of $\on{log}(\Psi_{\on{KZ}})$. When $A$ and the 
$b[\zeta],\zeta\in \mu_N(\CC)$ commute, the renormalized holonomy of 
(\ref{basic:eqn}) is equal to 
$$
\exp(\sum_{\zeta\in \mu_N(\CC), \zeta\neq 1} \log(\zeta^{-1}-1) b[\zeta]).
$$ 
Let $\hat\f_{N+1}:= \hat\f(A,b[\zeta],\zeta\in\mu_{N}(\CC))$. 
Then $\on{Log}(\Psi_{\on{KZ}}) \in \hat\f_{N+1}$, and 
its image in $\hat\f_{N+1}/[\hat\f_{N+1},\hat\f_{N+1}]$ is 
$\sum_{\zeta\in \mu_N(\CC), \zeta\neq 1} \log(\zeta^{-1}-1) b[\zeta]$. 
It follows that $\on{Log}(\Psi_{\on{KZ}}) \in \hat\f_{N+1}^{A,b[1]}$, where the 
latter space is the topological sum of all finely homogeneous components of 
$\hat\f_{N+1}$, except $\CC A$ and $\CC b[1]$ (we define a fine degree in $\NN
\tilde\mu_N = \NN(\mu_{n}(\CC)\cup \{0\})$ by $\on{deg}(A)=0$, 
$\on{deg}(b[\zeta]) = \zeta$; see Appendix \ref{app:a}).  

By Lazard elimination (see \cite{Re}), 
$\q := \hat\f_{N+1}^{A,b[1]}$ is freely generated by 
the $U_{k,l}$ ($k,l\geq 1$) and the $V_{k,l,\zeta}$ ($k,l\geq 0$, $\zeta\in
\mu_N(\CC)$, $\zeta \neq 1$), where 
$$
U_{k,l} = \ad(A)^{k-1}\ad(b[1])^{l-1}([A,b[1]]), \quad 
V_{k,l,\zeta} = \ad(A)^k \ad(b[1])^l (b[\zeta]). 
$$
Let $V(z) := (1-x)^{-b[1]} x^{-A} H(z)$, where $H(z)$ is 
defined by (\ref{basic:eqn}). Then 
$$
{{\on{d}}\over{\on{d}z}} V(z) = 
\big( (1-z)^{-b[1]}z^{-A} 
(\sum_{\zeta\in \mu_N(\CC)} {{b[\zeta]}\over{z-\zeta}}) 
z^A (1-z)^{b[1]} - {{b[1]}\over{z-1}}\big) V(z).  
$$ 
Therefore the image of $\log(\Psi_{\on{KZ}})$ in $\q^{\on{ab}} := 
\q/[\q,\q]$ is 
$$
\sum_{k,l\geq 1} {1\over{k!(l-1)!}}
\int_0^1 (-\log z)^k(-\log(1-z))^{l-1} {{\on{d}z}\over{z-1}}
\overline{U}_{k,l} 
+ \sum_{k,l\geq 0} {1\over {k!l!}} \int_0^1 (-\log z)^k (-\log(1-z))^l 
{{\on{d}z}\over {z-\zeta}} \overline{V}_{k,l,\zeta} 
$$
($x\mapsto \overline x$ is the projection $\q\to \q^{\on{ab}}$). 

\subsection{} The infinitesimal action of $\grtmd_{(\bar 1,1)}(N,\kk)$ on 
$\ul\Psdist(N,\kk)$ is given by 
\begin{equation} \label{inf:act:grtm}
\delta_{(\varphi,\psi)}(\Phi,\Psi) = (\Phi\varphi + D_\varphi(\Phi), 
\Psi\psi + \overline D_\psi(\Psi)), 
\end{equation}
where $D_\varphi,\overline D_\psi$ are given by (\ref{D:phi}) and 
(\ref{D:psi}). 

Define a bracket on $\hat\f_2 \oplus \hat\f_{N+1}$ by formulas (\ref{magnus:0}),
(\ref{magnus:1}) and (\ref{magnus:2}); this is a Lie bracket. Then
$\grtmd_1(N,\kk) \subset \hat\f_2 \oplus \hat\f_{N+1}$ is an inclusion of Lie
algebras. As in \cite{Dr2}, let $\p := \hat\f_2^{A,B}\subset \hat\f_2$ be the
topological sum of all homogeneous components of degree $>1$.

If $\alpha\geq 1$, denote by $\grtm_{(\bar 1,1)}(N,\kk)_\alpha$ the degree $\alpha$
part of $\grtm_{(\bar 1,1)}(N,\kk)$. Set 
$$
\grtmd_{(\bar 1,1)}(N,\kk)_1^+ := \grtmd_{(\bar 1,1)}(N,\kk)_1 \cap 
\big( \bigoplus_{a\in \ZZ/N\ZZ, a\neq 0} \CC (0,b(a))\big) , 
$$  
and 
$$
\grtmd_{(\bar 1,1)}(N,\kk)^+ := \grtmd_{(\bar 1,1)}(N,\kk)_1^+ \oplus
\big(\wh\bigoplus_{\alpha\geq 2} \grtmd_{(\bar 1,1)}(N,\kk)_\alpha \big).  
$$
Then $\grtmd_{(\bar 1,1)}(N,\kk) = \grtmd_{(\bar 1,1)}(N,\kk)^+ \oplus 
\kk (0,b(0))$; the first summand is a Lie subalgebra and the second summand 
is central. 

\begin{lemma} \label{lemma:abelianize}
The linear map $\hat\f_{N+1} \otimes\hat\f_{N+1} \to\hat\f_{N+1}$, 
$\psi_1 \otimes \psi_2 \mapsto \overline D_{\psi_1}(\psi_2)$, restricts to a 
linear map $\q\otimes \q\to \q$. The latter map induces a map $\q^{\on{ab}}
\otimes \q^{\on{ab}} \to \q^{\on{ab}}$.  
\end{lemma} 

{\em Proof.} Set $\q_1 := \oplus_{a\in \ZZ/N\ZZ, a\neq 0} \kk b(a)$, then 
$\q = \q_1 \oplus (\wh\oplus_{\alpha\geq 2} (\f_{N+1})_\alpha)$. 
If $\psi_1$ and $\psi_2$ are homogeneous in $\f_{N+1}$, then 
$\on{deg}(\overline D_{\psi_1}(\psi_1)) = \on{deg}(\psi_1) + \on{deg}(\psi_2)
\geq 2$, hence $\overline D_{\psi_1}(\psi_2) \in \q$. If now $\psi_1 =
[\psi',\psi'']$, where $\psi',\psi''\in \q$, then 
$\overline D_{\psi_1}(\psi_2) = [\overline D_{\psi_1}(\psi'),\psi''] 
+ [\psi',\overline D_{\psi_1}(\psi'')] \in [\q,\q]$. Hence $\q\otimes [\q,\q]$
maps to $[\q,\q]$. 

Let us now show that $[\q,\q] \otimes \q$ maps to $[\q,\q]$. For this, we will
show that if $\psi\in [\q,\q]$ and $\psi\in \f_{N+1}$, then $\overline
D_\psi(\psi') \in [\q,\q]$. It suffices to show this when $\psi' = A$ or $\psi'
= b(a)$, $a\in \ZZ/N\ZZ$. 
We have $\overline D_\psi(b(0)) = 0$. 
$[\q,\q]$ is an ideal in $\hat\f_{N+1}$, hence
$\overline D_\psi(A) = [\psi,A]\in [\q,\q]$. For the same reason, 
$\overline D_\psi(b(a)) \in [\q,\q]$. Hence 
it remains to show that $[\psi(A | b(a),\ldots,b(a+N-1)),b(a)]\in [\q,\q]$
when $a\neq 0$. This follows from the facts that $b(a)\in \q$,  
$\psi$ has degree $\geq 2$, and for any $d\geq 2$, the degree $d$ part of 
$\q$ is stable under the automorphism $A\mapsto A$, $b(\alpha)\mapsto 
b(a+\alpha)$. 
\hfill \qed \medskip 

\begin{lemma} \label{lemma:formulas}
Let us denote by $\psi_1\otimes \psi_2\mapsto \overline D_{\psi_1}(\psi_2)$
the above map $\q^{\on{ab}} \otimes \q^{\on{ab}} \to \q^{\on{ab}}$. We have 
$$
\overline D_{\overline U_{kl}}(\overline U_{k'l'}) = 
\overline U_{k+k',l+l'}, \; 
\overline D_{\overline U_{kl}}(\overline V_{k',l',\zeta}) = \overline
V_{k+k',l+l',\zeta}, 
$$
$$
\overline D_{\overline V_{k,l,\zeta}}(\overline U_{k'l'}) = \overline
V_{k+k',l+l',\zeta}, \; 
\overline D_{\overline V_{k,l,\zeta}}(\overline V_{k',l',\zeta'}) = 0
$$
if $(k,l,\zeta') \neq (0,0,\zeta^{-1})$, and 
$\overline D_{\overline V_{0,0,\zeta}}(\overline V_{k',l',\zeta^{-1}}) = 
\overline V_{k',l'+1,\zeta^{-1}}$. 

In the same way, the map $\p^{\on{ab}} \otimes \p^{\on{ab}} \to \p^{\on{ab}}$,
which we denote by $\varphi_1\otimes \varphi_2\mapsto D_{\varphi_1}(\varphi_2)$
is given by 
$$
D_{\overline A_{kl}}(\overline A_{k'l'}) = \overline A_{k+k',l+l'}, 
$$
where $A_{kl} = \on{ad}(A)^{k-1}\on{ad}(B)^{l-1}([A,B])$, 
$k,l\geq 1$ are the free
generators of $\p$, and $\overline A_{kl}$ are their images in $\p^{\on{ab}}$. 
\end{lemma}

Let us set $\overline\Psi_{\on{KZ}} := \Psi_{\on{KZ}}(-A | -b(0),\ldots,
-b(1-N))$. Recall that $\psi$ is uniquely determined by $\overline\Psi_{\on{KZ}}=
\Psi_{\on{KZ}}*\on{Exp}(\psi)$. Set $\Psi(t):=\Psi_{\on{KZ}}*\on{Exp}(t\psi)$. 
Then $\on{Log}\Psi(t)\in \q$; we decompose its image in $\q^{\on{ab}}$ as follows 
$$
[\log\Psi(t)] = \sum_{k,l\geq 1} c_{kl}(t) \overline U_{kl} + \sum_{\zeta\in
\mu_N(\CC),\zeta\neq 1} \sum_{k,l\geq 0} d_{k,l,\zeta}(t) \overline
V_{k,l,\zeta}. 
$$
We also decompose the image of $\psi$ in $\q^{\on{ab}}$ as 
$$
[\psi] = \sum_{k,l\geq 1} \varphi_{kl} \overline U_{kl}
+ \sum_{\zeta\in \mu_N(\CC),\zeta\neq 1} \sum_{k,l\geq 0}
\psi_{k,l,\zeta} \overline V_{k,l,\zeta}. 
$$
Then $\psi_{0,0,\zeta} = -\log(\zeta-1)(\zeta^{-1}-1)$ when $\zeta\neq 1$, 
and $\psi_{0,0,1}=0$. 

We will set $\varphi(u,v) = \sum_{k,l\geq 1} \varphi_{kl} u^k v^l$, 
$\psi_{\zeta}(u,v) = \sum_{k,l\geq 0} \psi_{k,l,\zeta} u^kv^l$, and define
$c(u,v)$ and $d_\zeta(u,v)$ similarly. According to (\ref{inf:act:grtm}), we
have $\Psi'(t) = \Psi(t) \psi + \overline D_\psi(\Psi(t))$, which according to
Lemma \ref{lemma:formulas} implies 
$$
{{\partial}\over{\partial t}} c(u,v,t) = \varphi(u,v)\big(c(u,v,t)+1\big), 
$$
$$ 
{{\partial}\over{\partial t}} d_\zeta(u,v,t) = 
\psi_\zeta(u,v) \big(1 + c(u,v,t) \big) 
+ \big( \varphi(u,v) - \psi_{0,0,\zeta^{-1} v}\big) d_\zeta(u,v,t). 
$$

According to \cite{Dr2}, we have 
$$
1+c(u,v,0) = {{\Gamma(1-u)\Gamma(1-v)}\over{\Gamma(1-u-v)}}
= \exp(\sum_{n\geq 2} {{\zeta(n)}\over n} (u^n + v^n - (u+v)^n )), 
$$
$c(u,v,1) = c(-u,-v,0)$, and since $1+c(u,v,t) = 
(1+c(u,v,0))e^{t\varphi(u,v)}$, we get 
\begin{align*}
\varphi(u,v) & = 
\on{log}{{\Gamma(1+u)}\over{\Gamma(1-u)}}
+ \on{log}{{\Gamma(1+v)}\over{\Gamma(1-v)}}
- \on{log}{{\Gamma(1+u+v)}\over{\Gamma(1-u-v)}}
\\ & 
= -2\sum_{n\on{\ odd\ }\geq 3} {{\zeta(n)}\over n}
(u^n + v^n - (u+v)^n). 
\end{align*}

On the other hand, we have 
$$
d_\zeta(u,v,0) = \int_0^1 z^{-u} (1-z)^{-v} {{dz}\over{z-\zeta}}
$$
and $d_\zeta(u,v,1) = -d_{\zeta^{-1}}(-u,-v,0)$. 
Since 
$$
d_\zeta(u,v,t) e^{-t(\varphi(u,v) - \psi_{0,0,\zeta^{-1}}v)} 
= d_\zeta(u,v,0) + \psi_\zeta(u,v) 
\big( 1 + c(u,v,0)\big) {{e^{t\psi_{0,0,\zeta^{-1}}v} -1}
\over{\psi_{0,0,\zeta^{-1}}v }}, 
$$
we get 
$$
\psi_\zeta(u,v) = \alpha_\zeta(u,v) + \alpha_{\zeta^{-1}}(-u,-v), 
$$
where 
$$
\alpha_\zeta(u,v) = 
{{v \log((1-\zeta)(1-\zeta^{-1}))}
\over{((1-\zeta)(1-\zeta^{-1}))^{-v}-1}} 
{{\Gamma(1-u-v)}\over{\Gamma(1-u)\Gamma(1-v)}} 
\int_0^1 {{z^{-u}(1-z)^{-v}}\over{z-\zeta}} dz. 
$$

\subsection{} 

Recall that $\q^{\on{ab}}$ is spanned by the $\overline V_{k,l,\zeta}$, 
$k,l\geq 0$, $\zeta\in \mu_N(\CC)$, $\zeta\neq 1$. Let us denote by 
$\q^{\on{ab}}_{l>0}$ the span of all $\overline V_{k,l,\zeta}$ 
where $l>0$; we denote the quotient 
$\q^{\on{ab}} / \q^{\on{ab}}_{l>0}$ by $\q^{\on{ab}}_{l=0}$. 

\begin{lemma}
The image of the sequence of maps $$[\grtmd_{(\bar 1,1)}(N,\kk), 
\grtmd_{(\bar 1,1)}(N,\kk)] \subset \grtmd_{(\bar 1,1)}(N,\kk)^+ 
\subset \q \to \q^{\on{ab}}$$ is contained in $\q^{\on{ab}}_{l>0}$, 
so we have a map 
$$
\grtmd_{(\bar 1,1)}(N,\kk)^+ / [\grtmd_{(\bar 1,1)}(N,\kk),
\grtmd_{(\bar 1,1)}(N,\kk)] \to \q^{\on{ab}}_{l=0}. 
$$
\end{lemma} 

{\em Proof.} 
Let us denote by $\q_{l>0}$ the part of $\q$ of positive 
relative degree w.r.t. $b(0)$. Then $[\q,\q] + \q_{l>0}$ is an ideal 
of $\q$ (for the bracket $[-,-]$). We will show that 
\begin{equation} \label{inclusion}
\langle \q,\q \rangle \subset [\q,\q] + \q_{l>0}.
\end{equation} 

We will use the formula 
$$
\langle \psi_1,\psi_2 \rangle = -[\psi_1,\psi_2] 
+ \overline D'_{\psi_2}(\psi_1) - \overline D'_{\psi_1}(\psi_2), 
$$
where $\overline D'_\psi$ is the derivation such that 
$\overline D'_\psi(A) = 0$, $\overline D'_\psi(b(a)) = 
- [\psi(A | b(a),\ldots,b(a+N-1)),b(a)]$. 

If $\psi\in \q$ has degree $\geq 2$ and $a\in \ZZ/N\ZZ$, then 
$\overline D'_\psi(b(a)) \in [\q,\q] + \q_{l>0}$. Therefore, if $\psi'\in \q$
is arbitrary, then $\overline D'_\psi(\psi') \in [\q,\q] + \q_{l>0}$. 
If $\psi$ and 
$\psi'\in \q$ both have degree $\geq 2$, we then get 
$\langle \psi, \psi'\rangle \in [\q,\q] + \q_{l>0}$. 

Let again $\psi\in\q$ be of degree $\geq 2$ and $a\in \ZZ/N\ZZ - \{0\}$. 
If $a'\in\ZZ/N\ZZ$, then $\overline D'_{b(a)}(b(a')) \in [\q,\q]
+ \q_{l>0}$. Therefore if $\psi\in \q$, we have 
$\overline D'_{b(a)}(\psi) \in 
[\q,\q] + \q_{l>0}$. Since $\overline D'_\psi(b(a))\in [\q,\q]$, 
and $[b(a),\psi]\in [\q,\q]$, we get $\langle b(a),\psi\rangle
\in [\q,\q] + \q_{l>0}$. This implies that $\langle b(a) + b(-a),\psi
\rangle \in [\q,\q] + \q_{l>0}$. 

Finally, if $a,a'\neq 0$ and $a+a'\neq 0$, we have $\langle b(a),b(a')
\rangle = [b(a+a'),b(a')-b(a)] - [b(a),b(a')] \in [\q,\q]$, which 
implies that $\langle b(a)+b(-a), b(a') + b(-a')
\rangle \in [\q,\q]$. 

This proves (\ref{inclusion}) and therefore the lemma. 
\hfill \qed \medskip 

\subsection{Generators of $\grtmd_{(\bar 1,1)}(N,\CC)^+$}

We will set 
$$
\grtmd_{(\bar 1,1)}(N,\kk)^{+,\on{ab}} := \grtmd_{(\bar 1,1)}(N,\kk)^+ /
[\grtmd_{(\bar 1,1)}(N,\kk),\grtmd_{(\bar 1,1)}(N,\kk)]
$$ 
and denote by $\grtmd_{(\bar 1,1)}(N,\kk)^{+,\on{ab}}_n$ the degree $n$ part
of this space. 

We have a linear map $\grtmd_{(\bar 1,1)}(N,\kk)^{+,\on{ab}} \to \p^{\on{ab}} 
\oplus \q^{\on{ab}}_{l=0}$. The degree $n$ part of 
$\p^{\on{ab}} \oplus \q^{\on{ab}}_{l=0}$ is spanned by the 
$(\overline A_{kl},0)$, where $k,l\geq 1$, $k+l=n$, and the 
$(0,\overline V_{n-1,0,\zeta})$, $\zeta\in \mu_N(\CC)$, $\zeta\neq 1$.
The image of $(\varphi_n,\psi_n)$ in 
$\p^{\on{ab}} \oplus \q^{\on{ab}}_{l=0}$ is equal to 
\begin{equation} \label{vector:1}
\big(0,\sum_{\zeta\in \mu_N(\CC), \zeta\neq 1} 
-2\log | 1-\zeta |   \overline V_{0,0,\zeta}
\big)
\end{equation}
if $n=1$, and to 
\begin{equation} \label{vector:2}
\Big( 
{{\zeta(n)}\over n} (1+(-1)^{n+1}) 
\sum_{k,l\geq 1,k+l=n} \pmatrix n \\ k \endpmatrix \overline A_{kl}, 
\sum_{\zeta\in\mu_N(\CC),\zeta\neq 1}
\big( 
Z(n,\zeta^{-1}) + (-1)^{n+1}Z(n,\zeta)
\big) \overline V_{n-1,0,\zeta} \Big),  
\end{equation}
where $Z(n,\zeta) = \sum_{k\geq 1} \zeta^k/k^n$ (so $Z(n,\zeta) =
\on{Li}_n(\zeta)$) for $n\geq 2$. The images of $(\varphi_n,\psi_n)$
in $\q^{\on{ab}}_{l=0}$ can be written uniformly in $n$ since 
$Z(1,z) = - \log(1-z)$.  

Recall that $(\ZZ/N\ZZ)^\times$ acts on $\grtmd_{(\bar 1,1)}(N,\CC)$. 
If $k\in (\ZZ/N\ZZ)^\times$, then the vectors obtained from 
(\ref{vector:1}), (\ref{vector:2}) by replacing $\zeta$ by $\zeta^k$ 
in $\log|1-\zeta|$ and $Z(n,\zeta^{\pm 1})$, also belong to the image of 
$\grtmd_{(\bar 1,1)}(N,\CC)^{+,\on{ab}} \to \p^{\on{ab}} \oplus 
\q^{\on{ab}}_{l=0}$. 

Let $\CC\mu_N(\CC)$ be the complex vector space with basis $[\zeta]$, 
where $\zeta\in \mu_N(\CC)$. Set for $n\geq 2$
$$
V_n := \on{Span}_{\CC} \{\sum_{\zeta\in\mu_N(\CC)}
\big( (-1)^{n+1}Z(n,\zeta^k) + Z(n,\zeta^{-k})\big) [\zeta], 
\; k\in (\ZZ/N\ZZ)^\times \} \subset \CC\mu_N(\CC)
$$
and 
$$
V_1 := \on{Span}_{\CC} \{\sum_{\zeta\in \mu_N(\CC), \zeta\neq 1}
-2\log |1-\zeta^k| [\zeta], \; k\in (\ZZ/N\ZZ)^\times  \}
\subset \CC\mu_N(\CC). 
$$

Then for any $n\geq 1$, 
$\on{ch} \grtm_{(\bar 1,1)}(N,\CC)^{+,\on{ab}}_n \geq \on{ch}(V_n)$, where 
$\on{ch}$ means the character as $(\ZZ/N\ZZ)^\times$-module. 

Computing $\on{ch}(V_n)$ means decomposing this space in 
isotypic components under the action of $(\ZZ/N\ZZ)^\times$. 
We first decompose $\CC\mu_N(\CC)$. 

Let $\chi : (\ZZ/N\ZZ)^\times \to \CC^\times$ be a character. 
The set $D_\chi := \{d | d$ divides $N$ and $\chi$ factors through 
$(\ZZ/N\ZZ)^\times \to (\ZZ/d\ZZ)^\times\}$ is closed under the operation
of taking greatest common divisors. The smallest element $f_\chi$ in 
$D_\chi$ is called the conductor of $\chi$. The resulting character 
$\chi' : (\ZZ/f_\chi\ZZ)^\times \to \CC^\times$ has conductor $1$. 

Let us denote by $\on{Prim}_d$ the set of primitive $d$th roots of $1$ 
in $\CC$. Then $\on{Prim}_d$ is a torsor under the action of 
$(\ZZ/d\ZZ)^\times$. It follows that 
$$
\CC \on{Prim}_d = \bigoplus_{\chi \on{\ character\ of\ }(\ZZ/d\ZZ)^\times}
(\CC \on{Prim}_d)^\chi, 
$$
and each summand is $1$-dimensional, spanned by 
$\sum_{a\in (\ZZ/d\ZZ)^\times} \chi(a)^{-1} [\zeta_d^a]$. 
We therefore obtain 
\begin{equation} \label{decomp:CmuN}
\CC\mu_N(\CC) = \bigoplus_{\chi\on{\ character\ of\ }(\ZZ/N\ZZ)^\times} 
\CC\mu_N(\CC)^\chi, 
\end{equation}
where
$$
\CC\mu_N(\CC)^\chi = \bigoplus_{d\on{\ such\ that\ }f_\chi|d|N}
(\CC\on{Prim}_d)^{\chi' \circ \pi_{\chi,d}}, 
$$
$\chi'$ is the character of $(\ZZ/f_\chi\ZZ)^\times$ obtained from 
$\chi$ and $\pi_{\chi,d} : (\ZZ/d\ZZ)^\times \to (\ZZ/f_\chi\ZZ)^\times$ 
is the natural projection. 

Since $V_n$ is $(\ZZ/N\ZZ)^\times$-invariant, it is graded 
w.r.t. the decomposition (\ref{decomp:CmuN}). 
Let us compute its graded components.  

The isotypic component corresponding to the character $\chi$
is spanned by 
\begin{align*}
v_\chi & := \sum_{\zeta\in \mu_N(\CC)}
\Big(\sum_{k\in (\ZZ/N\ZZ)^\times} \chi(k) \big(
(-1)^{n+1} Z(n,\zeta^k) + Z(n,\zeta^{-k})
\big)\Big)
[\zeta]
\\ & 
= \big( \chi(-1) + (-1)^{n+1} \big) 
\sum_{\zeta\in \mu_N(\CC)}
\big(\sum_{k\in (\ZZ/N\ZZ)^\times} \chi(k) Z(n,\zeta^k) \big)
[\zeta];   
\end{align*}
here we adopt the convention that $Z(1,\zeta) = -\log(1-\zeta)$ 
if $\zeta\neq 1$ and $Z(1,1)=0$.  
A character is called even (resp., odd) if $\chi(-1) = 1$
(resp., $-1$), so odd (resp., even) characters do not contribute if 
$n$ is odd (resp., even). 

\begin{proposition} \label{nonvan}
If $n$ and $\chi$ have opposite parity, then $v_\chi\neq 0$. 
\end{proposition}

{\em Proof.} Assume that $n>1$. 
The coefficient of $[\zeta_N]$ is equal to 
\begin{equation} \label{coord:zetaN}
2(-1)^{n+1}\sum_{k\in (\ZZ/N\ZZ)^\times} \chi(k) Z(n,\zeta_N^k) 
= 
2(-1)^{n+1}\sum_{a\geq 1} {1\over {a^n}} \big( \sum_{k\in (\ZZ/N\ZZ)^\times} 
\chi(k)\zeta_N^{ka} \big) . 
\end{equation}

\begin{lemma}
$\sum_{k\in (\ZZ/N\ZZ)^\times} \chi(k)\zeta_N^{ka}$
vanishes if $(N/f_\chi) \nmid a$, or if $a = (N/f_\chi) a'$, 
where $a'$ is an integer such that $(a',f_\chi) \neq 1$. 
If $a = (N/f_\chi)a'$ where $a'$ is an integer such that 
$(a',f_\chi) = 1$, then 
$$
\sum_{k\in (\ZZ/N\ZZ)^\times} \chi(k)\zeta_N^{ka}
= (N/f_\chi)\tau(\chi) \chi^{-1}([a']), 
$$
where $[a']$ is the class of $a'$ in $(\ZZ/f_\chi\ZZ)^\times$, and 
$\tau(\chi)$ is the Gauss sum $\tau(\chi) 
= \sum_{l\in (\ZZ/f_\chi\ZZ)^\times} \chi(l) \zeta_f^l$.  
\end{lemma}

{\em Proof of Lemma.}
We have 
$$
\sum_{k\in (\ZZ/N\ZZ)^\times} \chi(k)\zeta_N^{ka}
= 
\sum_{l \in [0,f_\chi-1], (l,f_\chi)=1} \chi(l) 
\zeta_N^{la} \big( 1 + \zeta_N^{f_\chi a} + \cdots 
+ \zeta_N^{({N\over f_\chi} - 1) f_\chi a} \big) , 
$$
so the l.h.s. is zero if $N \nmid f_\chi a$, and is equal to  
$(N/f_\chi) \sum_{l \in [0,f_\chi-1], (l,f_\chi)=1} \chi(l) 
\zeta_N^{la}$ otherwise. 

Let us place ourselves in the last situation, and let 
$a = (N/f_\chi)a'$. We must compute 
$\sum_{l \in [0,f_\chi-1], (l,f_\chi)=1} \chi(l) 
\zeta_N^{(N/f_\chi)a'l}
= \sum_{l \in (\ZZ/f_\chi\ZZ)^\times} \chi(l) 
\zeta_{f_\chi}^{a'l}$. When $(a',f_\chi) =1$, the change of variables
$l' := a'l$ and the multiplicativity of $\chi$ yield the result. 
If $(a',f_\chi)\neq 1$, let $f'$ be the order of $\zeta_{f_\chi}^{a'}$; 
$f'$ divides strictly $f_\chi$.  The last sum is rewritten as follows
$$
\sum_{l'\in (\ZZ/f'\ZZ)^\times} (\zeta_{f_\chi}^{a'})^{l'}
(\sum_{l\in (\ZZ/f_\chi\ZZ)^\times, l\equiv l'(f_\chi)} \chi(l)). 
$$
Since the restriction of $\chi$ to the kernel of 
$(\ZZ/f_\chi\ZZ)^\times \to (\ZZ/f'\ZZ)^\times$ is nontrivial, 
we get 
$$
\sum_{l\in (\ZZ/f_\chi\ZZ)^\times, l\equiv 1(f_\chi)} \chi(l)=0,
$$ 
and since $\{ l\in (\ZZ/f_\chi\ZZ)^\times, l\equiv l'(f_\chi) \}$ 
is a coset space under this kernel, 
$\sum_{l\in (\ZZ/f_\chi\ZZ)^\times, l\equiv l'(f_\chi)} \chi(l) = 0$. 
So $\sum_{l \in [0,f_\chi-1], (l,f_\chi)=1} \chi(l) 
\zeta_N^{(N/f_\chi)a'l} = 0$. 
\hfill \qed \medskip 

Now (\ref{coord:zetaN}) is equal to $2(-N/f_\chi)^{1-n} \tau(\chi)
L(n,\chi^{-1})$, where $L(s,\chi)$ is the Dirichlet 
$L$-function $\sum_{k\geq 1, (k,f_\chi)=1} \chi(k) / k^s$. 

The Euler expansion of $L(s,\chi)$ implies that when $n>1$, 
$L(n,\chi) \neq 0$ for any $\chi$. 
According to \cite{W}, Lemma 4.8, the Gauss sum $\tau(\chi)$ 
is also nonzero. It follows that $v_\chi\neq 0$. 

Let us assume now that $n=1$ and $\chi$ is even. The coefficient of 
$[\zeta]$ in $v_\chi$ is equal to $-2 \sum_{k\in(\ZZ/N\ZZ)^\times}
\chi(k) \log |1-\zeta^k|$. Assume that $\chi$ is nontrivial. If
$\zeta = \zeta_{f_\chi}$, then 
this coefficient is equal to 
$$
-2 {{\phi(N)}\over{\phi(f_\chi)}}\sum_{k\in(\ZZ/f_\chi\ZZ)^\times}
\chi(k) \log |1-\zeta_{f_\chi}^k|,  
$$
where $\phi$ is the Euler function. According to \cite{W}, Theorem 4.9, 
this is equal to 
$2 {{\phi(N)}\over{\phi(f_\chi)}} {{f_\chi}\over{\tau(\chi^{-1})}} 
L(1,\chi^{-1})$ (which can be rewritten 
$2 {{\phi(N)}\over{\phi(f_\chi)}} \tau(\chi) L(1,\chi^{-1})$
according to \cite{W}, Lemma 4.8), and according to \cite{W}, Corollary 4.4, 
$L(1,\chi^{-1}) \neq 0$ if $\chi\neq 1$, hence $v_\chi$ is nonzero. 

Assume now that $n=1$ and $\chi=1$. Let $p$ be a prime divisor of $N$. 
The coefficient of $[\zeta_p]$ in $v_1$ is equal to 
$$
-2 {{\phi(N)}\over{\phi(p)}} \sum_{k\in(\ZZ/p\ZZ)^\times}
\log|1-\zeta_p^k| = -2 {{\phi(N)}\over{\phi(p)}} \log(p), 
$$
hence $v_1$ is also $\neq 0$. 
\hfill \qed \medskip

So 
\begin{theorem} \label{thm:estimate}
Let $N\geq 2$ and $\nu$ be the number of distinct prime factors of $N$. 
Then: 

$\bullet$ $\grtmd_{(\bar 1,1)}(N,\CC)^{+,\on{ab}}_1 = \nu[1] + \sum_{\chi
\on{\ even\ character\ of\ }(\ZZ/N\ZZ)^\times} [\chi]$ as 
$(\ZZ/N\ZZ)^\times$-modules. The map $\grtm_{(\bar 1,1)}(N,\CC)_{1}^{+}
\to \on{Hom}(\on{Div}(N),\CC)$ is surjective; it is a bijection when restricted to the 
$(\ZZ/N\ZZ)^{\times}$-invariant part of the source of this map. 

$\bullet$ If $n\geq 3$ is odd, then
$\grtmd_{(\bar 1,1)}(N,\CC)^{+,\on{ab}}_n \geq \sum_{
\begin{matrix} \scriptstyle{\chi \on{\ even\ 
character}} \\
\scriptstyle{ \on{of\ }(\ZZ/N\ZZ)^\times} \end{matrix}}[\chi]$ as 
$(\ZZ/N\ZZ)^\times$-modules. 

$\bullet$ if $n\geq 2$ is even, then 
$\grtmd_{(\bar 1,1)}(N,\CC)^{+,\on{ab}}_n \geq
\sum_{\begin{matrix} \scriptstyle{\chi \on{\ odd\ 
character}} \\ \scriptstyle{\on{of\ }(\ZZ/N\ZZ)^\times}\end{matrix}}[\chi]$ as 
$(\ZZ/N\ZZ)^\times$-modules. 

In particular, set $d(n) := \on{dim}_{\QQ} \grtmd_{(\bar 1,1)}(N,\QQ)^{+,\on{ab}}_n$, 
$D(t) = \sum_{n\geq 1} d(n) t^n$. 
Then:  

$\bullet$ If $N=2$, $D(t) \geq t/(1-t^2)$; 

$\bullet$ if $N\geq 3$, $D(t)\geq {{\phi(N)}\over 2} {t\over {1-t}}
+ (\nu-1)t$. 
\end{theorem}

Note that $D(t)$ coincides with $f(t)$ from \cite{DG}, proof of Cor. 5.25. 

{\em Proof.} Proposition \ref{nonvan} implies the two last statements. 
Let us compute the character of $\grtmd_{(\bar 1,1)}(N,\CC)^{+,\on{ab}}_1$. 
Recall that $\grtmd_{{(\bar 1,1)}}(N,\QQ)^{+,\on{ab}}_1$ is the space of pairs 
$$
\big( \sum_{a\in\ZZ/N\ZZ} x(a) b(a), \rho \big) 
$$ 
where $\rho\in \QQ^\nu = \on{Hom}(\on{Div}(N),\QQ)$ and $x(a)\in \QQ$, 
such that: 

a) $x(0) = 0$,  

b) for any $a\in \ZZ/N\ZZ$, $x(a) = x(-a)$, and 

c) if $N = dN'$, then for any $a'\in \ZZ/N'\ZZ - \{0\}$, 
we have $x(da') = \sum_{k=0}^{d-1} x(a'+kN')$, and 

d) if $N = dN'$, then $\sum_{k=0}^{d-1} x(kN') = \rho(d)$.  

As in Lemma \ref{lemma:eliminate}, one shows that 
if $(x(a))_{a\in \ZZ/N\ZZ}$
satisfies a), b), c), then $d\mapsto \rho(d)$ defined by d) 
automatically belongs to $\on{Hom}(\on{Div}(N),\QQ)$. So 
$\grtmd_{(\bar 1,1)}(N,\QQ)^{+,\on{ab}}_1$ identifies with the space of 
$(x(a))_{a\in \ZZ/N\ZZ}$ satisfying a), b), c). According to \cite{W}, 
Theorem 12.18, we get 
\begin{equation} \label{eq:dim}
\on{dim}_{\QQ}
\grtmd_{(\bar 1,1)}(N,\QQ)^{+,\on{ab}}_1 = {{\phi(N)}\over 2} + \nu - 1.
\end{equation}  

The multiplicity of an odd $\chi$ in $\grtmd_{(\bar 1,1)}(N,\CC)^{+,\on{ab}}_1$
is $0$, and the multiplicity of an even $\chi\neq 1$ is $\geq 1$.
 
Let us compute the multiplicity of the trivial character $\chi = 1$, i.e., 
$\on{dim}_{\QQ} (\grtmd_{(\bar 1,1)}
(N,\QQ)^{+,\on{ab}}_1)^{(\ZZ/N\ZZ)^\times}$. 
We have a morphism $(\grtmd_{(\bar 1,1)}(N,\QQ)^{+,\on{ab}}_1)^* \to \RR$, 
taking $x(a)$ to $\log|1-\zeta_N^a|$. The image of the invariant part of 
$(\grtmd_{(\bar 1,1)}(N,\QQ)^{+,\on{ab}}_1)^*$ is the $\QQ$-linear span of all 
$f(\zeta) = \sum_{k\in (\ZZ/N\ZZ)^\times} \log|1-\zeta^k|$, 
where $\zeta\in \mu_N(\CC)$. 
Now either the order of $\zeta$ is a prime power $p^a$ dividing $N$, and
$f(\zeta) = {{\phi(N)}\over {\phi(p^a)}}\log(p)$; 
or the order of $\zeta$ is not a prime power, and $f(\zeta)=0$. Hence the 
invariant part of $(\grtmd_{(\bar 1,1)}(N,\QQ)^{+,\on{ab}}_1)^*$ 
maps surjectively to $\on{Span}_{\QQ} \{\log(p),\ p$ prime dividing $N\}$. 
By uniqueness of the prime powers decomposition, the $\log(p)$
form a free family over $\QQ$, hence $\on{dim}_\QQ \on{Span}_{\QQ} 
\{\log(p),\ p$ prime dividing $N\} = \nu$. So $\on{dim}_{\QQ} 
(\grtmd_{(\bar 1,1)}(N,\QQ)^{+,\on{ab}}_1)^{(\ZZ/N\ZZ)^\times} 
\geq \nu$.

Let us now show that 
$(\grtmd_{(\bar 1,1)}(N,\QQ)^{+,\on{ab}}_1)^{(\ZZ/N\ZZ)^\times} 
\to \on{Hom}(\on{Div}(N),\QQ)$ is bijective. According to what we have seen, 
it suffices to show that the composed map 
$\on{Hom}(\on{Div}(N),\QQ)^{*} \to 
(\grtmd_{(\bar 1,1)}(N,\QQ)^{+}_{1})^{*}
\to  \on{Span}_{\QQ}\{\on{log}(p)$, 
$p$ prime dividing $N\}$ is bijective. The first map 
takes $p$ to the linear map $(\sum_{a}x(a)b(a),\rho)
\mapsto \rho(p)$, which can be identified in view of d) with the 
linear map $\sum_{a}x(a)b(a)\mapsto \sum_{k=0}^{p-1}
x(kN/p)$. The image of this element by the second map is then 
$\sum_{k=0}^{p-1}\on{log}|1-\zeta_{N}^{kN/p}| = \on{log}(p)$. 
So the composed map is $\on{Hom}(\on{Div}(N),\QQ)^{*} \to
\on{Span}_{\QQ}\{\on{log}(p)$, $p$ prime dividing $N\}$, 
$p\mapsto \on{log}(p)$, which is bijective, as wanted.  

Comparing (\ref{eq:dim}) with the lower bounds obtained for the 
multiplicities, we obtained the wanted formula for the 
character of $\grtmd_{(\bar 1,1)}(N,\CC)^{+,\on{ab}}_1$. \hfill \qed \medskip

\section{$\grtm_{(\bar 1,1)}(N,\kk)$, $\grtmd_{(\bar 1,1)}(N,\kk)$ and special 
derivations} \label{sect:Ihara}

In this section, we relate $\grtm_{(\bar 1,1)}(N,\kk)$ and $\grtmd_{(\bar 1,1)}(N,\kk)$ 
with analogues of Ihara's algebra and Drinfeld's Hamiltonian interpretation 
(see \cite{Dr2}, Section 6). 

Recall that $\f_{N+1} = \f(A,b(a),a\in\ZZ/N\ZZ)$ and $C$ is given by 
$A + C + \sum_{a\in \ZZ/N\ZZ} b(a) = 0$. 

The dihedral group $D_N := (\ZZ/N\ZZ) \rtimes (\ZZ/2\ZZ)$ acts on $\f_{N+1}$
as follows: 

$\bullet$ the element $a'\in \ZZ/N\ZZ$ acts by the automorphism 
$\theta_{a'} : A\mapsto A$, $b(a) \mapsto b(a+a')$; 

$\bullet$ the element $\overline 1\in \ZZ/2\ZZ$ acts by the involutive 
automorphism $\theta : A\mapsto C$, $b(a) \mapsto b(-a)$. 

As in \cite{Ih1,Dr2}, let us say that derivation $\partial$ of $\f_{N+1}$
is special iff there exist $R_A,R_C,R_a \in \f_{N+1}$, such that 
$\partial(A) = [R_A,A]$, $\partial(C) = [R_C,C]$, $\partial(b(a)) = 
[R_a,b(a)]$. The special derivations form a Lie algebra 
$\on{SDer}(\f_{N+1})$. The inner derivations of $\f_{N+1}$ form a 
Lie ideal $\on{Inn}(\f_{N+1}) \subset \on{SDer}(\f_{N+1})$. 
The group $D_N$ acts by automorphisms of $\on{SDer}(\f_{N+1})$ and 
preserves $\on{Inn}(\f_{N+1})$; so does $(\ZZ/N\ZZ)^\times$.  

Let us denote by $\on{Ih}(N,\kk)$ the set of $\psi\in
\t^0_{3,N} = \f_{N+1}$, satisfying relations (\ref{octogon:grtm}), 
(\ref{add:cond:grtm}) (but not necessarily (\ref{pent:grtm})). 
This is a Lie algebra when equipped with the bracket 
(\ref{magnus:2}), (\ref{D:psi}). The elements $A,b(0)$ are 
central in $\on{Ih}(N,\kk)$; $\on{Ih}(N,\kk)$ is graded, 
we set $\on{Ih}^+(N,\kk)_1 := \on{Ih}^+(N,\kk) \cap 
\on{Span}_\kk(b(a),a\neq 0)$, $\on{Ih}^+(N,\kk) := \on{Ih}^+(N,\kk)_1 \oplus 
\on{Ih}(N,\kk)_{2} \oplus \cdots$. So $\on{Ih}(N,\kk) := \on{Ih}^+(N,\kk) 
\oplus (\kk A \oplus \kk b(0))$, where the first summand in a subalgebra
and the second summand is central. 

\begin{proposition}
We have an isomorphism of graded Lie algebras with $(\ZZ/N\ZZ)^\times$-action
$$
\on{Ih}^+(N,\kk) \simeq \big( \on{SDer}(\f_{N+1})/\on{Inn}(\f_{N+1}) \big)^{D_N},  
$$ 
taking $\psi$ to the class of the special derivation $\partial_\psi : A\mapsto 
[\psi(A| b(0),\ldots,b(N-1)) - \psi(C|b(0),\ldots,b(1-N)) , A]$, 
$b(a) \mapsto [\psi(A | b(a),\ldots,b(a+N-1)),b(a)]$. 
\end{proposition}

{\em Proof.} $\on{SDer}(\f_{N+1})/\on{Inn}(\f_{N+1})$ identifies 
with the quotient $\{$special derivations $\partial$ such that $R_C=0\} / \kk
\on{ad}(C)$. The statement is carried out directly in degree $1$. 
If now $\partial$ is such a special derivation of degree $>1$
(meaning that $R_A,R_a$ have degree $d>1$), then its class 
is $(\ZZ/N\ZZ)$-invariant iff $R_A$ is $(\ZZ/N\ZZ)$-invariant and  
$R_a = \theta_{a}(R_0)$. So $R_0,R_A$ are such that 
\begin{equation} \label{int:eq}
[R_A,A] + \sum_{a\in \ZZ/N\ZZ} [\theta_{a}(R_0),b(a)] = 0
\end{equation}
and $\partial$ is given by $A\mapsto [R_A,A]$, $b(a)\mapsto 
[\theta_{a}(R_0),b(a)]$. Now the special derivation representing
$\theta\circ \partial\circ \theta^{-1}$ is 
$A\mapsto [-\theta(R_A),A]$, $b(a) \mapsto [\theta\theta_{-a}(R_0) - \theta(R_A),
b(a)]$. We have $\partial = \theta\circ \partial \circ \theta^{-1}$, 
hence $R_A + \theta(R_A) =0$; and $R_A = \theta_{-a}(R_0) - \theta\theta_{a}(R_0)
= \theta_{-a}(\id - \theta)(R_0)$ for any $a$. Hence $R_A = (\id-\theta)(R_0)$, 
so $\partial = \partial_{R_0}$; 
and $R_A$ is $(\ZZ/N\ZZ)$-invariant, which means that $R_0$ satisfies
(\ref{octogon:grtm}); finally, (\ref{int:eq}) implies that $R_0$ 
satisfies (\ref{add:cond:grtm}). Hence we have an isomorphism of graded vector
spaces. The proof that it is a Lie morphism is as in \cite{Dr2}. 
\hfill \qed \medskip 

The Hamiltonian interpretation of $\on{Ih}(\kk)$ from \cite{Dr2} can be 
generalized as follows. We first recall some notation. Let $\cF(\f_{N+1})$
be the quotient of $\f_{N+1}^{\otimes 2}$ by the subspace generated by 
the elements $x\otimes y - y\otimes x$, $[x,y] \otimes z - x\otimes [y,z]$, 
$x,y,z\in \f_{N+1}$. As in \cite{Dr2}, any $f\in \cF(\f_{N+1})$ 
yields a universal function $f_\g : (\g^*)^{N+1} \to \kk$, where $\g$
is a finite dimensional Lie algebra with nondegenerate invariant 
scalar product. Using the Kostant-Kirillov Poisson bracket on $(\g^*)^n$, this 
interpretation allows to define a Lie bracket on $\cF(\f_{N+1})$. 

\begin{proposition} 
1) $\cF(\f_{N+1})^{D_N}\subset \cF(\f_{N+1})$ is a Lie subalgebra. 

2) There is an isomorphism of graded Lie algebras with 
$(\ZZ/N\ZZ)^\times$-action
$\on{Ih}(N,\kk) \simeq \cF(\f_{N+1})^{D_N}$. 
\end{proposition}

{\em Proof.} Similar to \cite{Dr2}; the map $\psi\mapsto f$ is given by 
the formulas 
$$
{{\partial f}\over {\partial A}} = (\id - \theta)(\psi), \; 
{{\partial f}\over{\partial b(a)}} = \theta(a)(\psi); 
$$
the fact that $f$ is $D_N$-invariant follows from the properties of 
$\psi$ and the identities $(\partial/\partial A) \circ \theta_1 = 
\theta_1 \circ (\partial/\partial A)$, $(\partial/\partial b(a)) \circ 
\theta_1 = \theta_1 \circ (\partial/\partial b(a-1))$, 
and $(\partial/\partial A) \circ \theta = 
- \theta \circ (\partial/\partial A)$, $(\partial/\partial b(a)) \circ 
\theta = \theta \circ (\partial/\partial b(-a))$. \hfill \qed \medskip 

In particular, the element $b(0)\in \on{Ih}(N,\kk)$ corresponds to 
${1\over 2} \sum_{a\in \ZZ/N\ZZ} (b(a),b(a))$, and $A$ corresponds to 
$-(A,C)$. 

We also define a Lie algebra $\on{Ihdist}(N,\kk)$ by forgetting the 
pentagon conditions in the definition of $\grtmd_{(\bar 1,1)}(M,\kk)$. 

If $N|N'$, we have Lie algebra morphisms $\pi_{NN'},\delta_{NN'} : 
\cF(\f_{N+1})^{D_N} \to \cF(\f_{N'+1})^{D_{N'}}$, respectively obtained 
by replacing $b(a)$ by $b'(a\on{\ mod\ }N')$ in a given expression, or 
by replacing it by $b'(a/d)$ if $d|a$, and by $0$ otherwise. 

Then we define a Lie algebra 
$\cF(N,\kk) \subset \cF(\f_{N+1})^{D_N} \times \kk^\nu$
as the set of pairs $(f,\rho)$ such that for any $N'|N$, 
$\pi_{NN'}(f) = \delta_{NN'}(f) + {{\rho(d)}\over 2}  
\sum_{a\in \ZZ/N\ZZ} (b(a),b(a))$. 

\begin{proposition}
$\on{Ih}(N,\kk) \simeq \cF(\f_{N+1})^{D_N}$ restricts to an isomorphism 
of graded Lie algebras with $(\ZZ/N\ZZ)^\times$-action 
$\on{Ihdist}(N,\kk) \simeq \cF(N,\kk)$.  
\end{proposition} 

{\em Proof.} Straightforward. \hfill \qed \medskip 

\section{Actions of $\on{Gal}(\bar\QQ/\QQ)$} \label{sect:galois}

In this section, we show that the natural action (\cite{Gr}) of $G_{\QQ}:= 
\on{Gal}(\bar\QQ/\QQ)$ on the orbifold fundamental 
group\footnote{If the group $\Gamma$ acts on the topological space $X$,
and $x_{0}\in X$, then $\pi_{1}^{orb}(X/\Gamma,[x_{0}])$ is the set of pairs
$(\gamma,\ell)$, where $\gamma\in \Gamma$ and $\ell$ is the isotopy class of a path
from $x_{0}$ to $\gamma x_{0}$; the composition is $(\gamma,\ell)(\gamma',\ell')
= (\gamma\gamma',\gamma\ell'\circ\ell)$. We have an exact sequence $1\to 
\pi_{1}(X,x_{0})\to \pi_{1}^{orb}(X/\Gamma,[x_{0}])\to\Gamma\to 1$.}  
of $(\CC^{\times} - \mu_{N}(\CC)) / D_{N}$
factors through $\wh{\on{GT}}$ in the profinite setup, 
and through $\on{GTMD}(N)_{l}$ in the pro-$l$ setup. 
Here $\ZZ/N\ZZ\subset D_{N}$ acts by multiplication (using the isomorphism 
$\ZZ/N\ZZ\simeq \mu_{N}(\CC)$ induced by $\zeta_{N}$) and 
$\bar 1\in\ZZ/2\ZZ$ acts by $z\mapsto 1/z$. 

We set $M_{N}:= \CC^{\times}-\mu_{N}(\CC)$, $X_{N}:= M_{N}/D_{N}$. 
We first describe the exact sequence $1\to \pi_{1}(M_{N})\simeq F_{N+1}
\to \pi_{1}^{orb}(X_{N}) \to D_{N}\to 1$
(the base-point $x_{0}$ is close to $0$ and $>0$). Recall that $B_{2}^{1}=
\langle \tau,\sigma \rangle/((\sigma\tau)^{2}=(\tau\sigma)^{2})$; the center
of $B_{2}^{1}$ is $Z(B_{2}^{1})\simeq \ZZ = \langle (\tau\sigma)^{2}\rangle$. 
We have an exact seqence $1\to K_{2,N}\to B_{2}^{1}\to (\ZZ/N\ZZ)^{2}\rtimes 
\SG_{2}\to 1$, where the second map is $\tau\mapsto (\bar 1,\bar 0)\in 
(\ZZ/N\ZZ)^{2}\subset (\ZZ/N\ZZ)^{2}\rtimes \SG_{2}$ and $\sigma\mapsto 
(21)\in \SG_{2}\subset (\ZZ/N\ZZ)^{2}\rtimes\SG_{2}$; the kernel is presented by 
$K_{2,N} = \langle X_{01},X_{02},x_{12}(0),...,x_{12}(N-1) 
\rangle / (X_{02}x_{12}(0)...x_{12}(N-1)X_{01}$ is central$)$, where $X_{01}=
\tau^{N}$, $x_{12}(\alpha)=\tau^{\alpha}\sigma^{2}\tau^{-\alpha}$, and 
$X_{02} = (\sigma\tau\sigma^{-1})^{N}$. Then $Z(B_{2}^{1})\cap K_{2,N}$
is generated by $X_{02}x_{12}(0)...x_{12}(N-1)X_{01} = (\tau\sigma)^{2N}$. The 
image of $(\tau\sigma)^{2}$ in $(\ZZ/N\ZZ)^{2}\rtimes \SG_{2}$ is 
$(\bar 1,\bar 1)\in (\ZZ/N\ZZ)^{2}$, which generates a central $\ZZ/N\ZZ$. 
Factoring the above exact sequence by $0\to Z(B_{2}^{1})\cap K_{2,N}
\to Z(B_{2}^{1})\to \ZZ/N\ZZ\to 0$ (which coincides with $0\to \ZZ
\stackrel{N\times -}{\to}\ZZ\to \ZZ/N\ZZ\to 0$), we get an exact sequence 
$1\to F_{N+1}\to B_{2}^{1}/Z(B_{2}^{1})\to D_{N}\to 1$, where 
$F_{N+1}:= K_{2,N}/(X_{02}x_{12}(0)...x_{12}(N-1)X_{01} =1)$ is freely 
generated by $X_{01},x_{12}(0),...,x_{12}(N-1)$. We then have 
\begin{equation} \label{isom:orb}
\pi_{1}^{orb}(X_{N},[x_{0}]) \simeq B_{2}^{1}/Z(B_{2}^{1})\simeq 
\langle \tau,\sigma \rangle/((\tau\sigma)^{2}=1);  
\end{equation}
the projection $\pi_{1}^{orb}(X_{N},[x_{0}]) \stackrel{\tilde\varphi_{N}}{\to} 
D_{N}= (\ZZ/N\ZZ)\rtimes 
(\ZZ/2\ZZ)$ coincides with the above map, so is given by $\tau\mapsto (\bar 1,\bar 0)$, 
$\sigma\mapsto (\bar 0,\bar 0)$; this is the factorization by the relations 
$\tau^{N}=\sigma^{2}=1$. 

The isomorphism (\ref{isom:orb}) is described geometrically by $\tau\mapsto$
the counterclockwise path from $x_{0}$ to $\zeta_{N}x_{0}$, and $\sigma\mapsto$
the path from $x_{0}$ to $1/x_{0}$, contained in $\{z | 0\leq \on{arg}(z)\leq 2\pi/N\}
- \{0,1,\zeta_{N}\}$. The prop-$l$ (resp., algebraic) fundamental group of $X_{N}$
is $\pi_{1}^{orb}(X_{N})(\tilde\varphi_{N},l)$ (resp., $\pi_{1}^{orb}
(X_{N})^{\wedge}$); this is $\on{Gal}(F_{l}/E)$, where $E = \bar\QQ(z)^{D_{N}}$, 
$F_{l}$ is the maximal pro-$l$ algebraic extension of  $E$ in $\bar\QQ(z^{1/n},n>0)$, 
unramified outside $0,\infty,\mu_{N}(\bar\QQ)$ (resp., $\on{Gal}(F/E)$, where
the pro-$l$ condition is dropped in the definition of $F$). 

The action of $G_{\QQ}$ on $\pi_{1}^{orb}(X_{N})(\tilde\varphi_{N},l)$
factors through $\on{GTMD}(N)_{l}$, and the action of the latter group 
takes $(\lambda,\mu,f,g)$ to the automorphism $\tau\mapsto \tau^{\mu}$, 
$\sigma\to g^{-1}(\tau,\sigma^{2})\sigma^{\lambda}g(\tau,\sigma^{2})$
(one can prove that the octogon identity is equivalent to the condition that this
is an automorphism of $\pi_{1}^{orb}(X_{N})(\tilde\varphi_{N},l)$).
In the same way, the action of $G_{\QQ}$ on  $\pi_{1}^{orb}(X_{N})^{\wedge}$ 
factors through $\wh{\on{GT}}$ and it given by the same formulas. 

When $N=1,2,4$, the subgroup $\Gamma_N$ of $\on{PSL}_{2}(\CC)$ of 
automorphisms of $\CC^{\times} - \mu_{N}(\CC)$ is strictly larger than $D_{N}$. 
Taking into account these additional symmetries, one defines subgroups 
$\on{GTMD}'(N)_{l} \subset \on{GTMD}(N)_{l}$, such that we 
have morphisms $G_{\QQ}\hookrightarrow \wh{\on{GT}}\to \on{GTMD}'(N)_{l}
\to \on{Aut}(\pi_{1}^{orb}(\CC^{\times} - \mu_{N}(\CC)/\Gamma_{N})^{(l)})$. 
We will return to this question elsewhere. 

For completeness, we recall the action of $G_{\QQ}$ on $\pi_{1}^{orb}(\CC^{\times}
- \{1\} / \SG_{3})^{(l)}$ (when $N=1$, $\Gamma_{N}=\SG_{3}$), described in 
\cite{Dr2}, and show its compatibility with its action on 
$\pi_{1}^{orb}(\CC^{\times} - \{1\} / D_{1})^{(l)}$ described above. 
The exact sequence $1\to F_{2}\to \pi_{1}^{orb}(\CC^{\times} - \{1\}/\SG_{3},
[x_{0}])\to \SG_{3}\to 1$ coincides with $1\to F_{2}\to B_{3}/Z(B_{3})
\to \SG_{3}\to 1$, where $B_{3} = \langle \sigma_{1},\sigma_{2} \rangle
/(\sigma_{1}\sigma_{2}\sigma_{1}=\sigma_{2}\sigma_{1}\sigma_{2})$; 
$Z(B_{3}) = \langle (\sigma_{1}^{2}\sigma_{2})^{2}\rangle\simeq\ZZ$, 
and the morphism $B_{3}/Z(B_{3})\to \SG_{3}$ is $\overline{\sigma}_{1}
\mapsto (12)$, $\overline\sigma_{2}\mapsto (23)$; it coincides with the 
factorization by the relations $\overline\sigma_{1}^{2} 
= \overline\sigma_{2}^{2}=1$. The identification $\pi_{1}^{orb}(\CC^{\times}
- \{1\}/\SG_{3})\simeq B_{3}/Z(B_{3})$ takes $\overline\sigma_{1}$ to the path
$x_{0}\mapsto 1-x_{0}$ contained in $]0,1[$, and $\overline\sigma_{2}$ to the 
path $x_{0}\mapsto 1/x_{0}$ contained in $\{z|\Im(z)\geq 0\}-\{0,1\}$. 

We then have an inclusion $\pi_{1}^{orb}(\CC^{\times}-\{1\}/D_{1},[x_{0}])
\subset \pi_{1}^{orb}(\CC^{\times}-\{1\}/\SG_{3},[x_{0}])$, corresponding to 
$B_{2}^{1}/Z(B_{2}^{1})\to B_{3}/Z(B_{3})$, $\bar\tau\mapsto 
\bar\sigma_{1}^{2}$, $\bar\sigma\mapsto\bar\sigma_{1}$. The action of 
$\on{GT}_{l}$ on $\pi_{1}(\CC^{\times}-\{1\}/\SG_{3})^{(l)}$ takes $(\lambda,f)$
to $\overline\sigma_{1}\mapsto \overline\sigma_{1}^{\lambda}$,   
$\overline\sigma_{1}\overline\sigma_{2}\overline\sigma_{1}
\mapsto \overline\sigma_{1}\overline\sigma_{2}\overline\sigma_{1}
f(\overline\sigma_{1}^{2},\overline\sigma_{2}^{2})$; this restricts to the action of 
$\on{GT}_{l}$ on $\pi_{1}(\CC^{\times}-\{1\}/D_{1})^{(l)}$ as
$\overline\sigma_{1}\overline\sigma_{2}\overline\sigma_{1}
f(\overline\sigma_{1}^{2},\overline\sigma_{2}^{2}) = 
\overline\sigma_{1}^{\lambda} f^{-1}(\overline\sigma_{1}^{2},\overline
\sigma_{2}^{2})\overline\sigma_{2}^{\lambda}f(\overline\sigma_{1}^{2},
\overline\sigma_{2}^{2})f(\overline\sigma_{1}^{2},\overline\sigma_{2}^{2})
\overline\sigma_{1}^{\lambda}$ (which is used implicitly in 
\cite{Dr2} and follows from the hexagon identity).

\begin{appendix} 

\section{The coefficients of $\Psi_{\on{KZ}}$} \label{app:a}

Let $A$, $b[\zeta]$, $\zeta\in \mu_N(\CC)$ be free variables 
(here $A$ and $b[\zeta]$ have degree $1$; 
the identification is $A = t_0^{01}$, $b[\zeta] = t[\zeta]^{12}$). 
Then $\Psi_{\on{KZ}}$ may be identified with the renormalized 
holonomy from $0$ to $1$ of  
\begin{equation} \label{basic:eqn}
{ {\on{d}} \over {\on{d}z} } H(z) = \big( {A\over z}
+ \sum_{\zeta\in \mu_N(\CC)} {{b[\zeta]}\over{z-\zeta}}\big) H(z) ,  
\end{equation}
i.e. $\Psi_{\on{KZ}} = H_1^{-1}H_0$, where $H_0(z) \sim z^A$ and $H_1(z) \sim 
(1-z)^{b[1]}$. 

Le and Murakami's formula (\cite{LM}) for the KZ associator can 
be generalized to $\Psi_{\on{KZ}}$, showing that $\Psi_{\on{KZ}}$ 
may be viewed as a generating series for multiple polylogarithms (MPL's) 
at $N$th roots of unity.  

We first define these numbers. If $\zeta_1,\ldots,\zeta_r\in \mu_N(\CC)$, 
and $s_1,\ldots,s_r$ are positive integers such that $(s_r,\zeta_r)\neq 
(1,1)$, we set
$$
\on{Li}_{s_1,\ldots,s_r}(\zeta_1,\ldots,\zeta_r) = 
\sum_{1\leq n_1 < \ldots < n_r}
{{\zeta_1^{n_1} \cdots \zeta_r^{n_r}}\over
{n_1^{s_1} \cdots n_r^{s_r}}}. 
$$ 
It is well-known that these numbers correspond bijectively to 
iterated integrals. 
Set $\tilde\mu_N = \{0\} \cup \mu_N(\CC)$. If $(\tilde\zeta_1,\ldots,\tilde\zeta_m)$
is a sequence with values in $\tilde\mu_N$, such that $\tilde\zeta_1\neq 0$
and $\tilde\zeta_m\neq 1$, we set 
$$
I(\tilde\zeta_1,\ldots,\tilde\zeta_m) = \int_{0\leq t_1\leq \ldots \leq t_r \leq 1}
\omega_{\tilde\zeta_1}(t_1) \wedge \cdots \wedge 
 \omega_{\tilde\zeta_r}(t_r) , \quad \omega_0 = \on{d}t/t, \omega_{\zeta} = 
 \on{d}t/(\zeta-t). 
$$
Then 
$$
\on{Li}_{s_1,\ldots,s_r}(\zeta_1,\ldots,\zeta_r) = 
I((\zeta_1\cdots \zeta_r)^{-1},0^{(s_1-1)},(\zeta_2\cdots \zeta_r)^{-1},
0^{(s_2-1)},\ldots,\zeta_r^{-1},0^{(s_r-1)}) 
$$
(here $0^{(i)}$ means $0$ repeated $i$ times). 

\begin{proposition}
Set $b[0]:=A$. If $(\tilde\zeta_1,\ldots,\tilde\zeta_r)$ is a sequence in 
$(\tilde\mu_N)^r$, we set $w(\tilde\zeta_1,\ldots,\tilde\zeta_r)
:= b[\tilde\zeta_r] \cdots b[\tilde\zeta_1]$. If $I \subset \{i|\tilde\zeta_i=0\}$, 
$J\subset \{i|\tilde\zeta_i=1\}$, set $w(\tilde\zeta_1,\ldots,\tilde\zeta_r)^{I,J}
:= \prod_{i\in [1,r] \setminus(I\cup J)} b[\tilde\zeta_i]$ (the product is in
decreasing order) and 
$$
\tilde w(\tilde\zeta_1,\ldots,\tilde\zeta_r) := 
\sum_{I\subset \{i|\tilde\zeta_i=1\}, 
J\subset \{i|\tilde\zeta_i=0\}} (-1)^{\on{card}(I) + \on{card}(J)}
b[1]^{\on{card}(I)} w(\tilde\zeta_1,\ldots,\tilde\zeta_r)^{I,J}
b[0]^{\on{card}(J)}. 
$$
Then 
$$
\Psi_{\on{KZ}} = \sum_{r\geq 0} \sum_{(\tilde\zeta_1,\ldots,\tilde\zeta_r)
\in (\tilde\mu_N)^r | \tilde\zeta_1\neq 0, \tilde\zeta_r\neq 1}
(-1)^{\on{card}\{i\in [1,r] | \tilde\zeta_i\neq 0\}}
I(\tilde\zeta_1,\ldots,\tilde\zeta_r) \tilde w(\tilde\zeta_1,\ldots,\tilde\zeta_r). 
$$
\end{proposition}

{\em Proof.} We denote $\Psi_{\on{KZ}}$ by $\Psi$ in the course of 
this proof, which is parallel to \cite{LM}. Let
$\CC\langle\!\langle b[\tilde\zeta],\tilde\zeta\in \tilde\mu_N\rangle\!\rangle$ 
be the degree completion of the free algebra with generators 
$A,b[\zeta],\zeta\in\mu_N(\CC)$. Let $\alpha,\beta$ be additional 
indeterminates of degree $1$ commuting with each other and with the 
$b[\tilde\zeta]$ ($\tilde\zeta\in\tilde\mu_N$). Let $\CC\langle\!\langle
b[\tilde\zeta]\rangle\!\rangle_0 \subset \CC\langle\!\langle
b[\tilde\zeta]\rangle\!\rangle$ be the complete subspace 
spanned by the monomials
not starting with $b[1]$ and not ending with $b[0]$. 
We have a direct sum decomposition 
$\CC\langle\!\langle b[\tilde\zeta] \rangle\!\rangle 
= \CC\langle\!\langle b[\tilde\zeta] \rangle\!\rangle_0 \oplus  
\big( b[1] \CC\langle\!\langle b[\tilde\zeta] \rangle\!\rangle + 
\CC\langle\!\langle b[\tilde\zeta] \rangle\!\rangle b[0] \big)$. 
We denote by $\pi : \CC\langle\!\langle b[\tilde\zeta] \rangle\!\rangle 
\to \CC\langle\!\langle b[\tilde\zeta] \rangle\!\rangle_0$ the corresponding 
projection map. We also have two maps
$$
\CC\langle\!\langle b[\tilde\zeta] \rangle\!\rangle 
\stackrel{u}{\to} 
\CC\langle\!\langle b[\tilde\zeta] \rangle\!\rangle [[\alpha,\beta]]
\stackrel{v}{\to} 
\CC\langle\!\langle b[\tilde\zeta] \rangle\!\rangle,  
$$
defined by: $u$ is an algebra morphism, $u(b[0]) = b[0]-\alpha$, 
$u(b[1]) = b[1]-\beta$, $u(b[\zeta]) = b[\zeta]$ for $\zeta\in 
\mu_N(\CC) -\{1\}$; $v$ takes $w\alpha^p \beta^q$ to 
$b[1]^q w b[0]^p$, where $w\in \CC\langle\!\langle b[\tilde\zeta] 
\rangle\!\rangle$. 

Then $v\circ u(b[1] \CC\langle\!\langle b[\tilde\zeta] \rangle\!\rangle
) = v\circ u(\CC\langle\!\langle b[\tilde\zeta] \rangle\!\rangle
b[0])=0$. Therefore $v\circ u = v\circ u \circ \pi$. 

Analogues of $G_0,G_1$ for the equation (\ref{basic:eqn}), where 
$A,b[1]$ are replaced by $A - \alpha$, $b[1] - \beta$ are  
$G_0(z) z^{\alpha} (1-z)^{\beta}$, $G_1(z) z^{\alpha}
(1-z)^{\beta}$, hence the holonomy of this equation 
is $\Psi$. So $u(\Psi) = \Psi$ (a constant series in $\alpha,\beta$), 
so $v\circ u(\Psi) = \Psi$.  It follows that 
\begin{equation} \label{Psi:Psi}
\Psi = v\circ u \circ \pi(\Psi).
\end{equation} 
We have 
$$
\Psi = \on{lim}_{\eps\to 0}
\eps^{-b[1]}
\big( \sum_{(\tilde\zeta_1,\ldots,\tilde\zeta_r) \in (\tilde\mu_N)^r}
(-1)^{\on{card}\{i\in[1,r]|\tilde\zeta_i\neq 0\}}
I_\eps(\tilde\zeta_1,\ldots,\tilde\zeta_r)
w(\tilde\zeta_1,\ldots,\tilde\zeta_r) \big) 
\eps^{b[0]}, 
$$
where $I_\eps(\tilde\zeta_1,\ldots,\tilde\zeta_r)$ is defined as
$I(\tilde\zeta_1,\ldots,\zeta_r)$ with integration domain $\eps \leq
t_1\ldots\leq t_r \leq 1-\eps$. 
Applying $\pi$ to this identity, we get 
\begin{equation} \label{Psi:reg}
\pi(\Psi) = \sum_{(\tilde\zeta_1,\ldots,\tilde\zeta_r)
\in (\tilde\mu_N)^r | \tilde\zeta_r\neq 1, \tilde\zeta_1\neq 0} 
(-1)^{\on{card}\{i\in[1,r]|\tilde\zeta_i\neq 0\}}
I(\tilde\zeta_1,\ldots,\tilde\zeta_r)
w(\tilde\zeta_1,\ldots,\tilde\zeta_r). 
\end{equation}  
The result now follows from (\ref{Psi:Psi}), (\ref{Psi:reg}) and the fact that 
$v\circ u(w(\tilde\zeta_1,\ldots,\tilde\zeta_r)) 
= \tilde w(\tilde\zeta_1,\ldots,\tilde\zeta_r)$.  
\hfill \qed \medskip 

\begin{remark}Using this proposition, one may use \cite{EG} to give
an explicit formula for $\on{log}(\Psi_{\on{KZ}})$. \end{remark}


\end{appendix}

\end{document}